\documentclass[a4paper,12pt]{book}

\usepackage{amscd}
\usepackage{amsfonts}
\usepackage{amsmath}
\usepackage{amssymb}
\usepackage{amstext}
\usepackage{amsthm}
\usepackage{graphicx}
\usepackage{makeidx}
\usepackage{pifont}
\usepackage{psfrag}

\makeindex

\newcommand{\C}{{\cal C}}
\newcommand{\dd}{{\rm d}}
\newcommand{\E}{{\cal E}}
\newcommand{\N}{\mathbb{N}}
\newcommand{\I}{1\!{\rm I}}
\newcommand{\NN}{{\cal N}}
\newcommand{\PPP}{{\cal P}}
\newcommand{\PP}{\mathbb{P}}
\newcommand{\R}{\mathbb{R}}
\newcommand{\X}{{\cal X}}
\newcommand{\Z}{\mathbb{Z}}

\newcommand{\cadlag}{c\`adl\`ag}

\newtheorem{theo}{Theorem}[chapter]
\newtheorem{lemma}{Lemma}[chapter]
\newtheorem{prop}{Proposition}[chapter]
\newtheorem{cor}{Corollary}[chapter]
\newtheorem{claim}{Claim}[chapter]

\begin{document}

\frontmatter

\parskip 0pt

\cleardoublepage

\thispagestyle{empty}
\begin{center}
{\large INSTITUTO DE MATEM\'ATICA PURA E APLICADA}
{\rule{130mm}{0.3mm}}
\end{center}

\vspace{20mm}
\begin{center}
{\bf Doctoral Thesis}
\end{center}
\vspace{10mm}

\baselineskip 2.3em

\begin{center}
{\LARGE
\textbf{\sc Generalized Hammersley Process \\ and \\ Phase Transition for Activated Random Walk Models}
}
\end{center}

\baselineskip 1.3em

\begin{center}

\normalsize \vspace{10mm}

\vspace{15mm}

{\large {\bf Leonardo Trivellato Rolla}}

\vspace{15mm}
\end{center}

\begin{center}

{\bf Advisor: {Vladas Sidoravicius}}
\end{center}

\vfill

\begin{center}
\textsc{March 20, 2008}
\end{center}

\clearpage \thispagestyle{empty}

\parskip 6pt

\cleardoublepage

\thispagestyle{empty}

\ \vfill
\begin{flushright}
\textit{%
this thesis is dedicated to my good friends
}
\end{flushright}

\clearpage \thispagestyle{empty}

\cleardoublepage

\thispagestyle{empty}

\ \vfill

{\it
\begin{flushright}
\begin{tabular}{rl}
\ \ \ \
&
\begin{tabular}{l}
{\sc Ah! Os Rel\'ogios}
\\
\
\\
Amigos, n\~ao consultem os rel\'ogios
\\
quando um dia eu me for de vossas vidas
\\
em seus f\'uteis problemas t\~ao perdidas
\\
que at\'e parecem mais uns necrol\'ogios...
\\
\\

Porque o tempo \'e uma inven\c c\~ao da morte:
\\
n\~ao o conhece a vida - a verdadeira -
\\
em que basta um momento de poesia
\\
para nos dar a eternidade inteira.
\\
\\

Inteira, sim, porque essa vida eterna
\\
somente por si mesma \'e dividida:
\\
n\~ao cabe, a cada qual, uma por\c c\~ao.
\\
\\

E os Anjos entreolham-se espantados
\\
quando algu\'em - ao voltar a si da vida -
\\
acaso lhes indaga que horas s\~ao...
\\
\\
\\
{\hfill \small \rm M\'ario Quintana - A Cor do Invis\'\i vel}
\end{tabular}
\end{tabular}
\end{flushright}
}

\clearpage \thispagestyle{empty}

\chapter{Abstract}

This thesis consists of two parts that are independent of each other.

In the first part I report on a joint work with V.~Sidoravicius concerning the Activated Random Walk Model.
This is a conservative particle system on the lattice, with a Markovian continuous-time evolution.
Active particles perform random walks without interaction, and they may as well change their state to passive, then stopping to jump.
When particles of both types occupy the same site, they all become active.
This model exhibits phase transition in the sense that for low initial densities the system locally fixates and for high densities it keeps active.
Though extensively studied in the physics literature, the matter of giving a mathematical proof of such phase transition remained as an open problem for several years.
In this work we identify some variables that are sufficient to characterize fixation and at the same time are stochastically monotone in the model's parameters.
We employ an explicit graphical representation in order to obtain monotonicity.
This representation has a very useful commutativity property that allows direct, constructive approaches.
With this method we prove that there is a unique phase transition for the one-dimensional finite-range random walk.

In the second part I report on a joint work with V.~Sidoravicius, D.~Surgailis, and M.~E.~Vares about the Broken Line Process.
In this work we introduce the broken line process and derive some of its properties.
Its discrete version is presented first and a natural generalization to the continuum
is then proposed and studied.
The broken lines are related to the Young diagram and the Hammersley process and are useful for computing last passage percolation values and finding maximal oriented paths.
For a class of passage time distributions there is a family of boundary conditions that make the process stationary and reversible. For such distributions there is a law of large numbers and the process extends to the infinite lattice.
One application is a simple proof of the explicit law of large numbers for last passage percolation with exponential and geometric distributions.

\clearpage \thispagestyle{empty}

\chapter{Acknowledgments}

So many people and institutions have been crucial for my work during the last four years that I couldn't
mention them all.
I apologize to all those that are omitted here.

I thank CNPq and FAPERJ for financial support. 

I would also like to thank all the staff of IMPA for their efficiency, dedication, and availability that have provided me a pleasant working environment.

D. Jane  helped me considerably to improve the quality of the English presentation. The mistakes that remain are, of course, all my responsibility.

I am indebted to Profs.~A.~Ram\'\i rez,
J.~van~den~Berg,
and
W.~Werner,
respectively for their hospitality during my visits to
PUC in Santiago,
CWI in Amsterdam,
and ENS in Paris.

I am deeply grateful to A.~Q.~Teixeira and to Prof.~V.~Beffara\- for enjoyable mathematical discussions that undoubtedly have inspired me very much.
I also thank Prof.~A.~Ram\'\i rez for fruitful discussions.

I thank Profs.~R.~Dickman, G.~Kozma, B.~N.~B.~de~Lima, E.~Pujals, and M.~Viana for accepting the invitation to compose the thesis committee.
Special thanks to B.~N.~B.~de~Lima for his suggestions and careful reading.

I find it really hard to express properly all my gratitude to my friend and advisor, Prof.~V.~Sidoravicius.
The words that follow are a poor attempt to this task.
Vladas certainly had a huge influence on my formation as a person and as a professional.
First of all he taught me how to live in the scientific world with grace and joy, which cannot be found in any book.
His criticisms were slowly assimilated and have contributed to make me start developing my own critical sense.
Vladas' way of thinking of mathematics and doing mathematics is quite singular: I have partly inherited his ways and I benefit a lot from this.
Much more than playing with epsilons, deltas and fancy theories, he showed to me a particular way of doing mathematics using intuition.
It was mostly because of him that I became passioned by the beauty of Probability.
Beyond all that, Vladas has been an endless source of support.
Whenever I thought I couldn't make it to the end, he would immediately stop whatever he was doing and spend as long as it would take to convince me otherwise.
His encouragement and enthusiasm have pushed me to move on and on.
It is a great privilege to finish this doctorate under his guidance.

\clearpage \thispagestyle{empty}

\parskip 0pt
\tableofcontents

\clearpage \thispagestyle{empty}

\parskip 6pt

\cleardoublepage
\chaptermark{Introduction}

\mainmatter

\chapter*{Introduction}
\addcontentsline{toc}{chapter}{Introduction}

In Chapter~\ref{chap:arw} we study a system of activated random walks (ARW) on the infinite
one-dimensional lattice $\Z$.
We assume that there
are infinitely many particles in the system, each of which can be in one of two
states: $B$ (active) or $A$ (passive or sleeping).
Each $B$-particle performs an independent, continuous-time, simple random walk on $\Z$,
with the same jump rate, which we assume, without loss of generality, to be
equal to $1$.
We will be treating both the symmetric and asymmetric walk cases.
When a $B$-particle jumps to a site with one or more $A$-particles, any such particle at this site is immediately activated (i.e.,
switches to state $B$).
Each isolated $B$-particle goes to sleep (switches to state $A$), at a rate $\lambda >0$.
From this rule it follows that if two or more particles occupy the same site, then they are all of type $B$ or all of type $A$ (the latter situation can only arise in the initial condition).
According to this rule, at most one $B$-particle can go to sleep per site, and, if undisturbed, remain in state $A$ forever after.

This model belongs to a broad class of {interacting particle systems with conservation}, which have attracted great interest in physics, probability, and allied fields, in part because they afford simple
examples of phase transitions in systems maintained far from equilibrium.
In this {class of models the particles exist in two states that may be termed \emph{active} and \emph{passive}.
Such that activation of a passive particle requires the intervention of one or more active ones.}
This class includes the so-called {conserved lattice gases}~\cite{lubeck02b,lubeck02a,lubeck03b,lubeck03a,rossi00} and {stochastic sandpile models}~\cite{dickman02b,dickman08,dickman01,manna90,manna91}.
Such models exhibit {self-organized criticality}~\cite{bak87,bak88,dhar99,grinstein95} when coupled to a suitable control mechanism~\cite{dickman02a,dickman00}.

In fact, the stochastic conserved sandpile, generally known as \emph{Manna's Model}~\cite{manna90,manna91}, served as the primary motivation for the present study.
In infinite volume, this model is defined as follows.
Initially there are infinitely many particles, distributed in such a way that at each site of $\Z^d$ we have a Poisson number of particles with mean $\mu>0$.
Each site is equipped with an exponential rate $1$ clock, and each time a clock rings at a site bearing $2d$
or more particles, $2d$ particles jump to randomly chosen nearest neighbors.
Differently from the deterministic Bak-Tang-Wiesenfeld sandpile model \cite{bak87,bak88}, each particle chooses its direction among the $2d$ possibilities with probability $\frac1{2d}$, independent of any other particle.
In contrast to the deterministic sandpile \cite{dhar99}, very little is known rigorously about this system, and the ARW model is a reasonable caricature that seems to capture some essential aspects of Manna's model.

The ARW model may also be viewed as a special case of a {diffusive epidemic process}.
In this process, an infected particle performs a simple symmetric random walk with jump rate $D_B$, and recuperates at a given rate, while a healthy particle performs a simple symmetric random walk with jump rate $D_A$; healthy particles are infected on contact with infected ones.
The ARW model corresponds to $D_A = 0$. The model was introduced to probabilistic community (in the case $D_A = D_B$) 
in the late 1970's by F.~Spitzer, but due to its tremendous technical difficulties and complexity, remained unsolved until recently, when
 it was studied in detail in
\cite{kesten03,kesten05,kesten06,kesten08}.
The diffusive epidemic process has also been studied via renormalization group techniques and numerical simulation \cite{freitas00,fulco01b,fulco01a,janssen01,kree89,oerding89,wijland98}.
A general conclusion from these
studies is that there are three distinct regimes of critical behavior, for
$D_A < D_B$, $D_A = D_B$ and $D_A > D_B$.  
It is not yet clear whether the ARW model falls within the remit of the first regime, or, alternatively, that $D_A =0$ marks a special case.

Here we suppose that each site $x \in \Z$ initially contains a certain number $\eta^B_0(x)$ of $B$-particles and no $A$-particles.
The $\eta^B_0(x)$ are i.i.d., mean $\mu$, Poisson random variables.

We are mainly interested in the phenomenon of fixation.
By fixation we mean that almost surely, for any finite volume $\Lambda$, there is a finite time $t_{\Lambda}$ such that after this time there are no $B$-particles within $\Lambda$.
When there is no fixation we expect that there is a limiting density of active particles for long times.

Numerical analysis and some general theoretical arguments suggest that the ARW model exhibits a phase transition in the parameters $\lambda$ and $\mu$, and that there should be two distinct regimes:

Low particle density.
There is a phase transition in $\lambda$ in this case.
Namely, there is $0<\lambda_c<\infty$ such that the system fixates for $\lambda>\lambda_c$ and does not fixate for $\lambda<\lambda_c$.
The threshold $\lambda_c$ should vanish when $\mu$ tends to zero and it should diverge as $\mu$ approaches the high-density regime.

High particle density. In this case there is no phase transition and the system does not fixate for any value of $\lambda$.

The second regime was studied in~\cite{kesten06}, where it was shown that indeed if the initial density is large enough there is no phase transition.

On the other hand, the behavior predicted for low particle densities, in spite of being nearly obvious, resisted many attempts to be established rigorously. 

In this work we estabilish part of the predictions for in the one-dimensional case.
The main result is that $\lambda_c<\infty$ for $\mu<1$ and $\lambda_c\to0$ as $\mu\downarrow0$.

We also prove that $\lambda_c=\infty$ for $\mu\geqslant 1$, characterizing what ranges of $\mu$ pertain to each regime.
Yet two conjectures remain open: that $\lambda_c>0$ for all $\mu>0$ and $\lambda_c\to\infty$ as $\mu\uparrow1$.

\bigskip

{\centering {\rule{40mm}{0.15mm}}

}

\bigskip

In Chapter~\ref{chap:brokenline} we study a \emph{continuous} variant of so-called  \emph{Broken Line Process}.
The meaning in which we use
word ``continuous'' requires more detailed explanation, which will be given below.

Speaking in more general terms, this process might be viewed as a generalization of well known Hammersley Process, introduced by Aldous and Diaconis, and independently by Rost. 

The Hammersley Process belongs to the broad class of Polygonal Markov Fields, which have been promoted by A. N. Kolmogorov in his later years.
However, due to technical difficulties arising from the interplay between stochastic geometry and an intention to construct and preserve a spatial analog of the Markov property, this topic became infamous for not only being difficult technically, but also conceptually.

As a consequence, until now it remains poorly understood.
A few fundamental works were produced by Arak and Surgailis~\cite{arak89,arak89b,arak89a,arak89c,arak90,arak91}, which had fair success among the community of statisticians and statistical physicists studying random patterns.
We find in writings of J.~Hammersley some indication that he had thought about particular cases of such fields, however without developing much of research in this direction.

Later on, in the context of applications to the First and Last Passage Percolation, the process which we call nowadays \emph{Hammersley Process}, was introduced by Aldous and Diaconis~\cite{aldous95}, and Rost~\cite{rostXX}, and on more general basis studied by Sidoravicius, Surgailis, and Vares~\cite{sidoravicius99}.
Results of~\cite{sidoravicius99} were extended to the last passage percolation model on the square lattice with geometric weights~\cite{sidoraviciusXX}.

In this work we generalize and extend~\cite{sidoraviciusXX} to the case of lattice models with exponential weights (and in this sense we use the word ``continuous'').
We derive basic properties of this random field, first of all its existence.
Finally we identify all possible stationary broken line processes related to both geometric or exponential weights to be given by a one-parameter family.
As a consequence we obtain the law of large numbers for the asymptotic velocity of last passage percolation by relatively simple and transparent arguments -- results that were previously known and are summarized in~\cite{johansson00}.
Fluctuations, however, are still beyond of the reach for these geometric techniques.

\chapter
{Phase Transition for Activated Random Walk Models}
\label{chap:arw}

\chaptermark{Phase Transition for ARW Models}

The model of activated random walks evolves as follows.
On $\Z^d$, at time zero there is an i.i.d. number of active particles whose expectation is $\mu$.
Each particle performs a continuous-time random walk with rate one.
When a particle is found alone at some site, it will change its state from active to passive when an exponential clock of rate $\lambda$ rings.
Once a particle is passive, it no longer jumps.
The particle becomes active again when some other particle jumps into the same site, if that ever happens.

The phenomenon we consider is that of local fixation.
By local fixation we mean that at any finite box a.s. there will be only passive particles or empty sites for large enough times.
Absence of fixation means that at every site there are at least two particles for arbitrarily large times.
A natural conjecture is that for low $\mu$ the system locally fixates and for large $\mu$ there is no fixation.

The first difficulty for the study of this model lies on the nature of interactions, that even though sound very simple, are quite severe to analyze and forbid the use of the most common available tools.
Namely the evolution of this model does not preserve any type of domination, or, in the language of particle systems, this model is not attractive.
We identify an observable that suffices to characterize fixation and we introduce an explicit graphical representation that recovers the monotonicity of this variable.
Such representation has also a very useful commutativity property that strongly favors direct, constructive approaches.
With this method we prove that there is a unique phase transition for the one-dimensional finite-range random walk.

This chapter is based on a joint work with
V.~Sidoravicius.

\section{Introduction}
\label{sec:sa_introduction}

We study a particle system that is informally described as follows.
The system starts with an i.i.d.~Poisson$(\mu)$ number of active particles
at each site of the lattice $S=\Z^d$.
As time passes these particles perform continuous-time random walks.
The particles do not interact except that when there is only one particle at a site, this particle may become passive, or fall asleep, at rate $\lambda>0$.
Once the particle is passive it will stop moving until some other particle jumps into that site.
We say that particles which are active are particles of type $B$ and particles which are passive are of type $A$. In the dynamics the $B$-particles jump at rate $D_B=1$ and $A$-particles jump at rate $D_A=0$;
a $B$-particle becomes type $A$ at rate $\lambda$ and an $A$-particle becomes type $B$ immediately if there is some other $B$-particle at the same site.

The state of the system is described by an element of
$\Sigma = (\Z_+\times\Z_+)^S$ and is denoted
$\eta=\bigl(\eta^A(x),\eta^B(x)\bigr)_{x\in S}$.
In this setting $\eta^A_t(x)$ denotes the number of $A$-particles (that is, particles which are passive) at position $x$ at time $t$;
$\eta^B_t(x)$ denotes the number of $B$-particles;
and $\eta^{AB}_t(x)=\eta^A_t(x)+\eta^B_t(x)$ denotes the total number of particles.
\index{eta@$\eta$}

It is known that there is an underlying probability space where a \cadlag\ strong Markov process taking values on an appropriate subset $\Sigma'\subseteq\Sigma$ and corresponding to the dynamics described above is defined. This technical issue is discussed in Section~\ref{sec:existence}.

A natural question that arises is whether the system fixates, which amounts to say that almost surely $\eta_t^B(x)=0$ for all $x$ in a given finite box $\Lambda$ and for all $t$ large enough.
In particular one wonders if there is phase transition, which means that there is fixation when the initial density is small enough and there is no fixation when it is large enough.
More precisely, phase transition means that for each fixed $\lambda$ there exists a $\mu_c\in(0,\infty)$ such that the systems fixates for $\mu<\mu_c$ and does not fixate for $\mu>\mu_c$.

A severe difficulty one has to plunge into when approaching this question lies on the nature of the interactions that take place.
To be more precise, most attempts to attack this problem using modern tools and techniques have failed for the following reasons.
Good estimates that locally solve the problem do not suffice; once the evolution halts at a finite region, there is no guarantee that in the future more particles from outside will not enter the controlled region.
This restarts the process and invalidates the argument.
A renormalization approach is the standard technique to overcome the shortcomings with the local argument.
Typically one gets better and better estimates for all scales, and finally combines them to obtain a global result from local estimates.
This type of approach breaks down because they rely on attractiveness, a property that our model unfortunately lacks.
(Another consequence of non-attractiveness is that it is not even entirely obvious that there is only one transition $\mu_c\in[0,\infty]$.)

\index{phase transition}
This chapter is mainly devoted to the proof of
\begin{theo}
\label{theo:fixation}
Consider the activated random walk model on the lattice $\Z^d$ with fixed parameter $\lambda>0$.
Then there exists $0 \leqslant \mu_c \leqslant \infty$ such that the system fixates for $\mu<\mu_c$ and does not fixate for $\mu>\mu_c$.

If $d=1$ and the random walk consists of finite-range jumps then $0<\mu_c \leqslant 1$.
Moreover, if the jumps are only to nearest neighbors then $\mu_c \geqslant \frac\lambda{1+\lambda}$.
\end{theo}

To overcome the difficulties that arise from the lack of attractiveness, we consider an explicit graphical construction that has strong commutativity and monotonicity properties.
This construction will form the core of the proof.

The fact that $\mu_c<\infty$ is intuitively obvious, since for $\mu>1$ there are `more particles than sites'.
Indeed this proof is quite simple once the framework discussed above is set up.
Less obvious is the conjecture that $\mu_c$ should be strictly less than one.
Another open question is whether there is fixation at $\mu=\mu_c$.
For the totally asymmetric, nearest-neighbor walk, a negative answer was given by Hoffman and Sidoravicius~\cite{hoffman07}.

Section~\ref{sec:existence} concerns technical questions such as the existence of the process and the approximation of the probability of local events on the original model by a large, finite one.
We present the graphical construction in Section~\ref{sec:flag} and study its properties.
In Section~\ref{sec:phasetransition} we prove Theorem~\ref{theo:fixation} using the background developed in the previous sections.
Finally, in Section~\ref{sec:concluding} we make some general comments on generalizations of the proofs presented here as well as some open questions.

\section
[Existence and approximation by finite boxes]
{Existence of the process and approximation by finite boxes}
\label{sec:existence}

In this section we discuss the construction of a strong Markov process with \cadlag\ trajectories whose dynamics correspond to the infinitesimal description given in the Introduction.

On the lattice $S=\Z^d$,
consider a translation-invariant random walk distribution $p(x,y)=p(0,y-x), x,y\in S$.
The formal generator of our process is given by
\begin{equation}
\index{generator}
\label{eq:generator}
  (Lf)(\eta) = \sum_x \lambda \delta_1( \eta^B(x) ) [ f(\eta^{(x)}) - f(\eta) ]
             + \sum_{(x,y)} \eta^B(x) p(x,y) [ f(\eta^{(xy)}) - f(\eta) ]
\end{equation}
for $f:\Sigma\to\R$, where
$$
  \eta^{(x)}(z) = \begin{cases} \eta(x), & z\ne x \\ (1,0),& z=x \end{cases}
$$
and
$$
  \eta^{(xy)}(z) = \begin{cases}
    \eta(z), & z\ne x,y \\
    \bigl(0,\eta^B(x)-1\bigr), & z = x \\
    \bigl(0,\eta^A(y)+\eta^B(y)+1 \bigr), & z=y
  \end{cases}.
$$

Of course the evolution is clearly well defined for any initial configuration in $\Sigma'' = \{\eta\in\Sigma: \sum_x \eta^{AB}(x) <\infty \}$. 
That is, given a probability distribution $P$ of the initial configuration $\eta \in \Sigma''$, there is a unique Markov process defined for $t\in[0,\infty)$ with values in $\Sigma''$ and having~(\ref{eq:generator}) as generator.
Let $\PP$ denote the law of the process with initial distribution $P$.

Given $M\in\N$,
take $P^\mu_M$ as a product measure on $\Sigma$ having marginals
$P^\mu_M\bigl(\eta^A(x)=0\bigr) = 1\ \forall\ x\in S$ and
$$
  P^\mu_M\bigl(\eta^B(x)=n\bigr) =
      \begin{cases} \displaystyle\frac{e^{-\mu}\mu^n}{n!}, & \|x\| \leqslant M \\ 
      \delta_{0}(n),          & \|x\| > M
      \end{cases}
      .
$$
As above, $\PP^\mu_M$ will denote the law of the process with initial distribution $P^\mu_M$.

The dynamics we want to construct cannot be defined on the whole space $\Sigma$ and we ought to restrict the set of allowed configurations.
Define
$\alpha(x)=\sum_{n=0}^\infty 2^{-n}p^{(n)}(x,0)$
and
$\|\eta\| = \sum_{x\in S}\bigl(|\eta^A(x)|+|\eta^B(x)|\bigr)\alpha(x)$.
Let $\Sigma'=\{\eta\in\Sigma:\|\eta\|<\infty\}$.
Then $\Sigma''$ is dense in $\Sigma'$, $P^\mu(\Sigma')=1$ and straightforward adaptations of the Andjel's construction~\cite{andjel82} imply the existence of a unique measure $\PP^\mu$ on $\Omega=D\bigl([0,\infty),\Sigma')$ with the property that
\begin{equation}
\label{eq:finitebox}
  \PP^\mu(A) = \lim_{M\to\infty} \PP^\mu_{M}(A)
\end{equation}
for any event $A$ such that, for some finite $\Lambda\subseteq S$ and $0<t<\infty$,
$A$ is measurable with respect to $\bigl(\eta^{A,B}_s(x)\bigr)_{x\in\Lambda,s\in[0,t]}$ .

\section
[A construction that preserves monotonicity]
{A graphical construction that preserves monotonicity}

\label{sec:flag}

Since the model is not attractive, some type of comparison between different parameters (initial density, jump rate, sleep rate), as well as the obtainment of bounds, are highly desirable for its study.
In this section we shall give a graphical construction of the dynamics.
This construction recovers the monotonicity of some variables of interest, allowing the desired comparisons.

It is simpler, clearer and more convenient to state and prove such properties for initial states with only finitely many particles on the lattice. In particular we shall construct $\PP^\mu_{M}$ explicitly.
Finally in Section~\ref{sec:phasetransition} we shall use~(\ref{eq:finitebox}), together with the results presented here, to study the phenomenon of fixation.

The construction holds on a general graph $S$ and with a general random walk distribution.
For the sake of clarity we assume $S=\Z^d$, $d=1$, and consider nearest-neighbor random walks with probability $p$ of jumping to the right and $q=1-p$ of jumping to the left.

\index{envelopes}
The dynamics is described in a informal way as follows. Every particle waits an exponential time of rate $1+\lambda$ and when the clock rings this particle will perform an action, which may be a jump or an attempt to sleep.
But the particle will not toss a coin to decide what to do. Instead, it is the \emph{site} where the particle is at that moment that will toss the coin.
It seems at first that it makes no difference but it is with that intuition that we conceive the explicit construction. In this setting we place envelopes at the sites.

Suppose there are finitely many particles on the system.
Let there be a universal clock that will ring with the appropriate rate and let there be an i.i.d. sequence of labels, independent of the clock.
\index{envelopes}%
Also, at each site, let there be an i.i.d. sequence of envelopes, each one containing some instruction to be performed.
When the clock rings for the first time, it will ring for the particle indicated by the first element in the label sequence and at that moment this particle will perform some action.
If the particle is passive, nothing happens.
If the particle is active, it will open the first envelope at that site, burn the envelope and perform the action written inside. The instruction may be `jump to the left', `jump to the right', or `try to sleep'. So we say there are two types of envelopes: jump envelops and sleep envelopes.
The envelope is burned regardless if there are other particles on the same site and the particle unsuccessfully attempts to sleep.

This graphical construction, that is in particular a natural coupling, has two nice properties that we describe below.

\begin{figure}[!ht]
{
\small
\psfrag{x}{$x$}
\begin{center}
 \includegraphics[width=13cm]{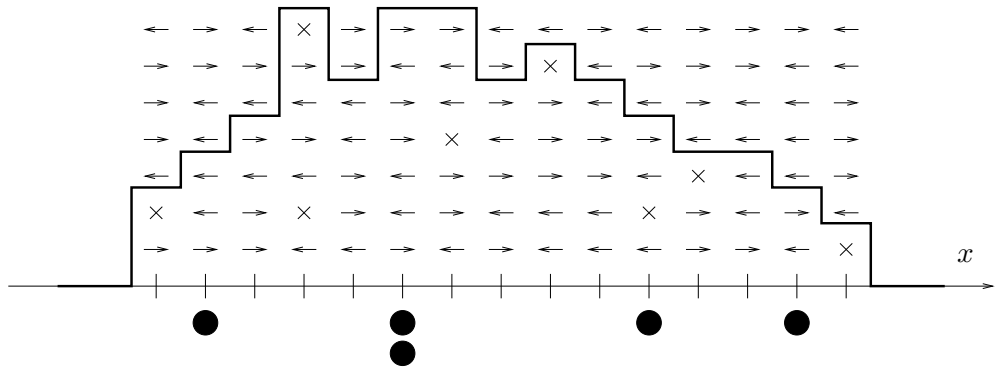}
\end{center}
\caption{
\small
\emph{commutativity}. The one-dimensional lattice with the particles' starting positions and the envelope sequences.
Particles are represented by black balls below the axis and envelopes are represented by arrows (jump envelopes) or crosses (sleep envelopes) just above the axis.
No matter in which order the evolution is performed, it will end up burning the envelopes below the bold line.}
\label{fig:commute}
}
\end{figure}

\index{commutativity}
The first property is \emph{commutativity} and says the following.
Suppose for a given realization of the process (universal clock, label sequence, envelopes) the system stabilizes, that is, all particles are passive for large enough time
(of course starting with finitely many particles this happens a.s.).
Then by changing the label sequence and the universal clock the system will stabilize at exactly the same state, except that some particles may be permuted.
Furthermore, the amount of envelopes burned at each site is also preserved.
In Figure~\ref{fig:commute} we show an example of a setup given by initial position of the particles and a stack of envelopes at each site.
The reader may evolve the system as she likes, and the evolution will terminate after burning exactly the envelopes indicated.

\begin{figure}[!ht]
{
\small
\psfrag{x}{$x$}
\begin{center}
 \includegraphics[width=10cm]{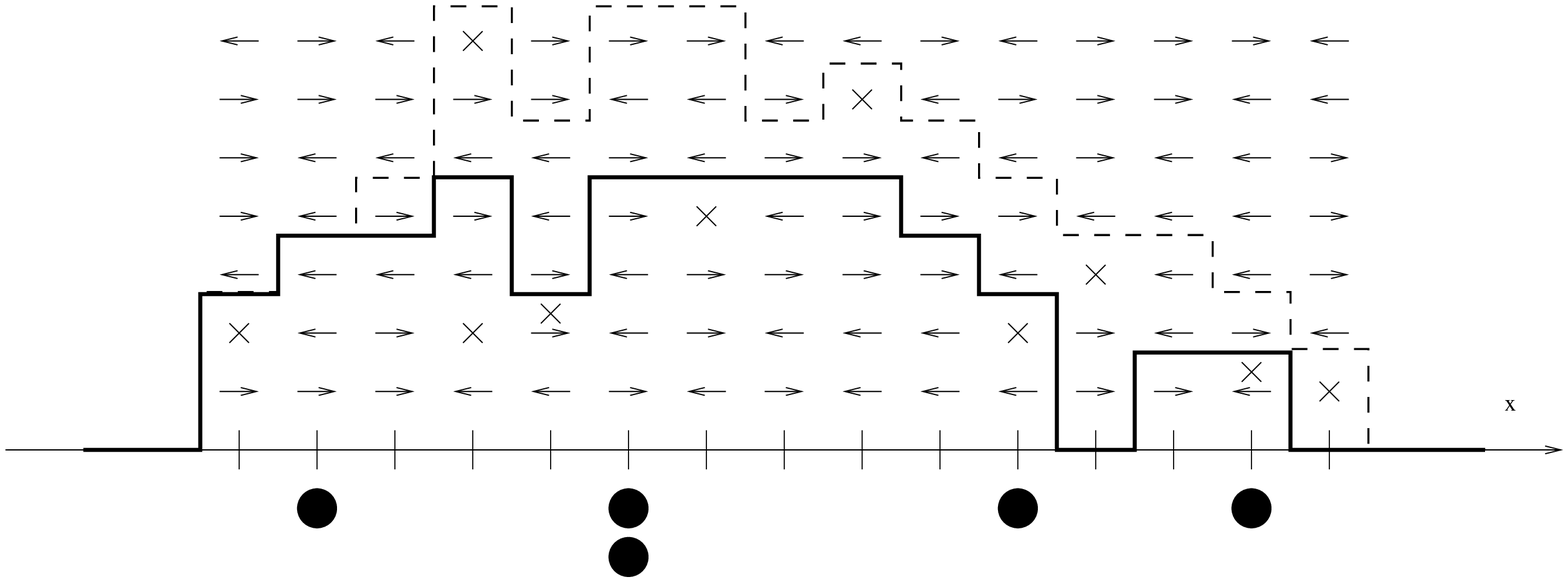}
 \vspace{7mm} \ 
 \includegraphics[width=10cm]{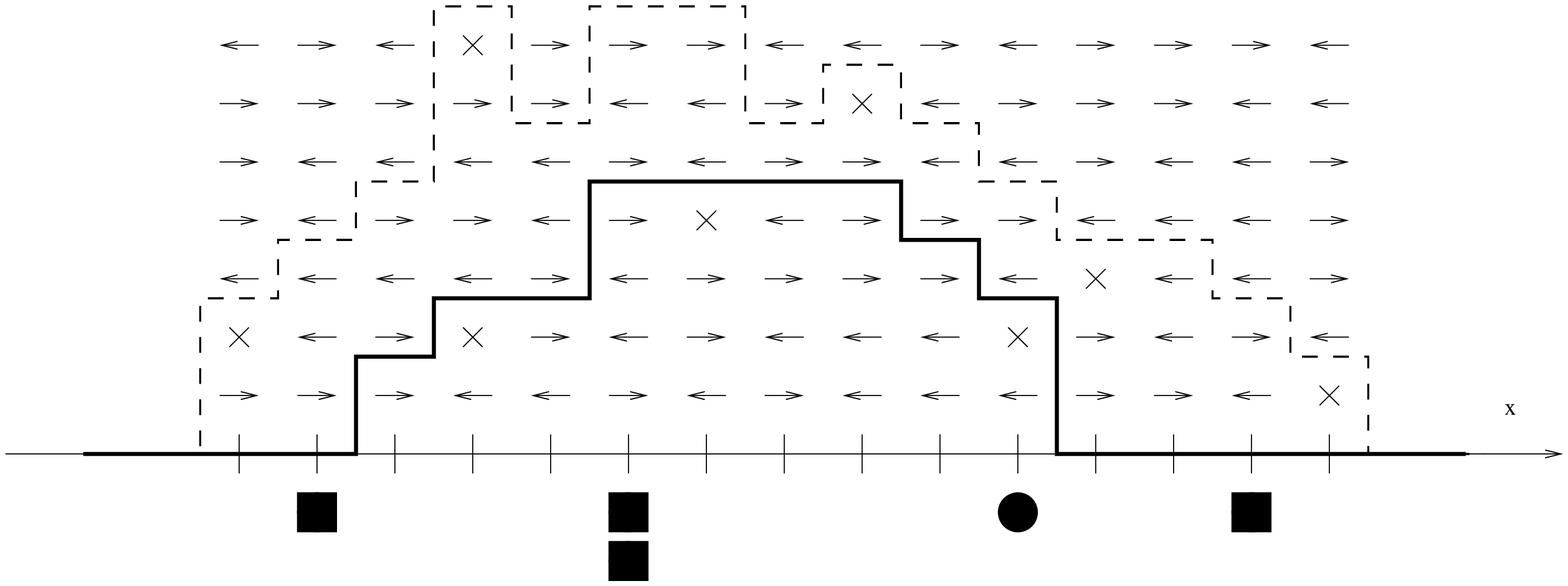}
 \vspace{7mm} \ 
 \includegraphics[width=10cm]{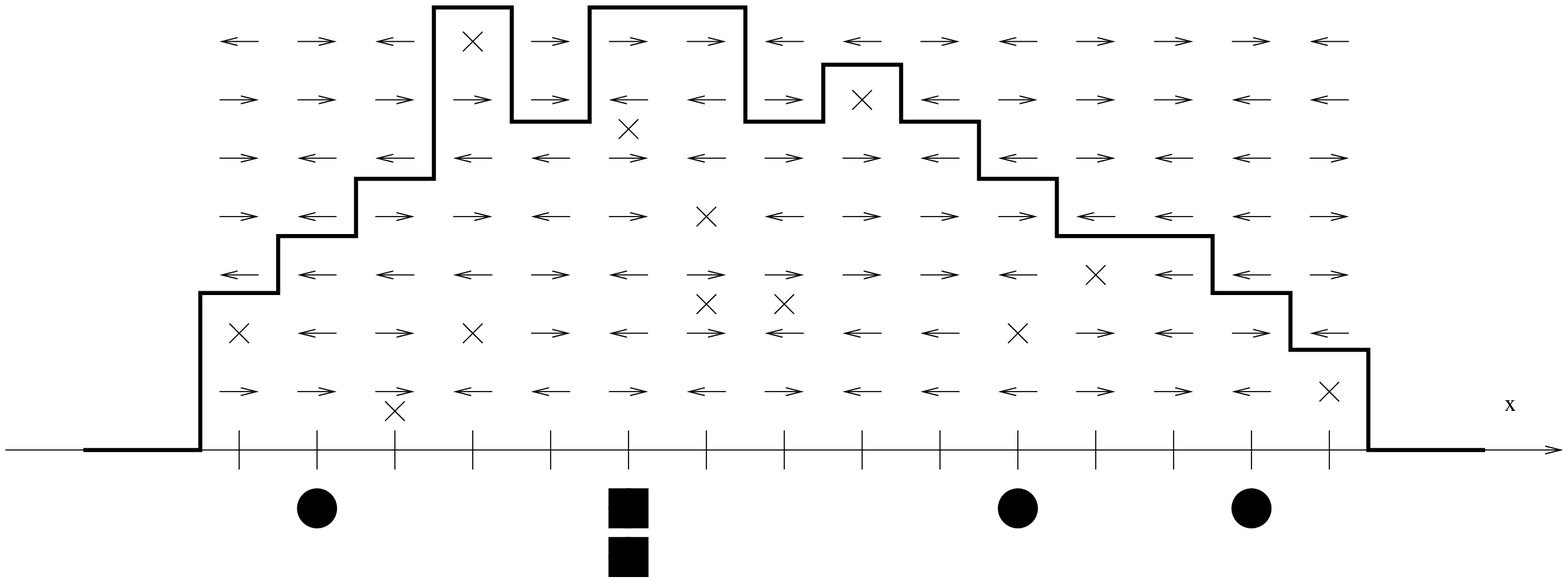}
\end{center}
\caption{
\small
\emph{monotonicity}.
First, inserting a couple of sleep envelopes has decreased the final profile of burned envelopes.
Second, starting with some particles passive instead of active (squares instead of balls) has caused an analogous effect.
Finally, even starting with quite a few extra sleep envelopes and some passive particles, the profile does not necessarily decrease.
}
\label{fig:monotone}
}
\end{figure}

So the final state of the system is determined by the initial conditions (positions and types of the particles) and the by sequences of envelopes.
\index{monotonicity}
The the second property, \emph{monotonicity}, states the following.
Suppose for some realization of the envelopes and initial conditions there is fixation, and take a new configuration by deleting some particles on the initial condition, changing the type of some particles from active to passive or inserting some sleep envelopes at some sites' sequences. 
Then for this new configuration the system also stabilizes and the final number of envelopes that are burned at each site (not counting the ones inserted) does not increase.
This is illustrated in Figure~\ref{fig:monotone}.

Graphical constructions in general have become standard tools in Probability and are usually attributed to Harris.
Actually it was known before but we have not determined exactly when it first appeared.
See Harris~\cite{harris78} and references therein to backtrack this history.
In the construction described above only the instructions are left at each site; the instant when each step should be performed is determined by some external factor.
What we call `envelope' here is called `router' or `rotor-router' by Levine~\cite{levine04}, Holroyd~et~al.~\cite{holroyd08}, and Levine and Peres~\cite{levine07}, or `stacks' by Propp and Wilson~\cite{propp98}, or finally `cards' by Diaconis and Fulton~\cite{diaconis91}.
The `commutativity property' is usually called `abelian property'.
\index{commutativity}
\index{abelian property|see{commutativity}}
As far as we have been able to find out, the earliest reference describing such a type of construction is \cite{diaconis91}.
The general framework considered by Eriksson~\cite{eriksson96} includes both the abelian sandpile model and the card game of Diaconis and Fulton.

The novelty introduced in this work is twofold.
The first thing is that even with the sleep envelopes this evolution is still commutative.
Finally, and perhaps the core of this chapter, the random variable `amount of visits to each site' increases with the initial amount of particles and decreases with the amount of sleep envelopes (and therefore it will be probably wise to try and approach the model by looking at this variable).

Before proceeding let us keep at the informal level for one more paragraph and explain how these properties may be used.

First, because of the commutativity property, one has the right to assume the label sequence is whatever looks convenient.
By doing so one can freely manipulate the scenario and dictate the destiny of each particle.
The evolution so obtained may have nothing to do with real life, but as long as it concerns the number of envelopes burned at each site (which is essentially the local time), things are fine.
Second, because of the monotonicity, the more particles are added to the initial state (which happens when the initial box $[-M,M]^d$ or the parameter $\mu$ increase) and the less sleep envelopes are present (which happens when $\lambda$ decreases), the more envelopes will be burned.
Therefore one can for instance ignore the existence of some sleep envelopes when manipulating the scenario, obtaining as a final result an upper bound for the local times.
All this will be made more clear in the next section, when these properties will be used in practice.

A formal description of our construction is given now.

Let $M$ be fixed, we shall construct the process having law $\PP^\mu_{M}$.
For each site $\|x\|\leqslant M$, let $\eta^B_x(0)$ be i.i.d. with distribution Poisson($\mu$).
Let $m=\sum_x\eta^B_x(0)$ be the total number of particles.
Label these particles $\rho_1,\dots,\rho_m$, let $x_1,\dots,x_m$ be their corresponding positions and let $\tau_1=\dots=\tau_m=B$ denote their type, that is, `active' or `$B$'.
Take a Poisson process $(0<t_1<t_2<\cdots)$ on $\R_+$ with intensity $m(1+\lambda)$.
Now take a sequence of independent, uniformly distributed indexes $N=(n_1,n_2,n_3,\dots)$ where $n_i\in\{1,\dots,m\}$.
For each site $x\in\Z^d$, take an i.i.d. sequence $(F_{x,j})_{j\in\N}$, where $F_{x,j}$ is $-1$ with probability $\frac q{1+\lambda}$, $+1$ with probability $\frac p{1+\lambda}$ and $0$ with probability $\frac \lambda{1+\lambda}$.
\index{Omega@$\tilde\Omega$}%
\index{omega@$\tilde\omega$}%
Let $\tilde\Omega$ be the space consisting on elements $\tilde\omega=\bigl( (x,\tau), F, T, N \bigr)$ such that $x\in S^m$, $\tau\in\{A,B\}^m$, $F \in \{-1,0,1\}^{S\times \N}$, $T=(0=t_0<t_1<t_2<\cdots)$ with $t_i\uparrow\infty$, and $N\in\NN_m$ for some $m\geqslant 0$, where $\NN_m \subseteq \{1,\dots,m\}^\N$ denotes the set of sequences where each symbol appears infinitely many times.
\index{N@$\NN_m$}
Let $\PPP^\mu_M$ denote the probability distribution on $\tilde\Omega$ given by the above
description.

We now show how to construct the process $\eta^{A,B}_t(x)$, $x\in S$, $t\geqslant0$
from $\tilde\omega$.
The system state is represented by $X = \bigl((x, \tau ), J\bigr)$, where $(x,\tau)=(x_n, \tau_n)_{n=1}^m$, $J = (j_x)_{x\in S}$, $x_n$ gives the position of the particle $\rho_n$, $\tau_n$ gives its state ($A = $ passive, $B =$ active) and $j_x$ tells how many envelopes have been burned at site $x$. We
denote the coordinates of $X$ by $x_n(X)$, $\tau_n(X)$ and $j_x(X)$. The set of possible states
with $m$ particles is denoted $\X_m = S^m \times \{A,B\}^m \times \Z_+^S$.

To recover the original representation, where the particles are undistinguishable, we define $\eta[X]$ by
$$
  \eta^{A}(x) = \sum_{n=1}^m \delta_x \bigl(x_n(X)\bigr)
                               \delta_{A} \bigl(\tau_n(X)\bigr)
$$
and analogously for $\eta^{B}(x)$.

With $F$ fixed, we take $X_0 = \bigl((x,\tau),J_0\bigr)$, where $J_0 \equiv 0$ and $\tilde\omega = \bigl((x,\tau),F,T,N\bigr)$.
We shall obtain $X_{i}$ by updating $X_{i-1}$ according to the action performed when the clock rings for particle $\rho_{n_{i}}$, and we denote this operation by $X_i = n_i\cdot X_{i-1}$.
This yields a sequence of states $(X_0,X_1,X_2,\dots)$.

The process $\bigl(\eta_t(\tilde\omega)\bigr)_{t\geqslant0}$ will then be given by
$
  \eta_t(\tilde\omega) = \eta[X_i]  \mbox{ for }  t\in[t_i,t_{i+1})
$.
It is easy to see that that the process $(\eta_t)_{t\geqslant0}$ so obtained is a \cadlag\ strong Markov process with countably many possible states and its infinitesimal generator indeed corresponds to the dynamics we are considering, that is, $\PPP^\mu_M(\eta(\tilde\omega)\in A) = \PP^\mu_M(A)$.

It remains to define $n\cdot X$.
When the clock rings for some active particle it will look for the next envelope in the sequence of the corresponding site $x$, that is, $F_{x,j+1}$ and perform its instruction: if $F_{x,j+1}=\pm1$ the particle jumps to site $x\pm1$, if $F_{x,j+1}=0$ the particle tries to fall asleep (and succeeds iff it is alone at site $x$).

Let $\tilde n \in\{1,\dots ,m\}$ and $X \in \X_m$. We are going to define $\tilde n \cdot X$.
If $\tau_{\tilde n} = A$, take $X' = X$.
Otherwise, let $\tilde x = x_{\tilde n}(X)$ and $\tilde j = j_{\tilde x}+1$.
If $F_{\tilde x,\tilde j} = 0$ take $x_n' = x_n$ for all $n$ and
$$
  \tau_n' = \begin{cases}
               \tau_n, & n\ne\tilde n \\
               A,      & n=\tilde n \mbox{ and } x_l\ne x_{\tilde n}\ \forall\ l\ne\tilde n \\
               B,      & n=\tilde n \mbox{ and } x_l =  x_{\tilde n} \mbox{ for some } l\ne\tilde n 
            \end{cases}.
$$
If $F_{\tilde x,\tilde j}=\pm1$, take
$$
  x_n' = \begin{cases}
              x_n,     & n\ne\tilde n \\
              x_n\pm1, & n = \tilde n
         \end{cases}
$$
and
$$
  \tau_n' = \begin{cases}
              \tau_n,     & x_n \ne x_{\tilde n}\pm 1 \\
              B,          & x_n  =  x_{\tilde n}\pm 1
         \end{cases}.
$$
Take
$$
  j_x' = \begin{cases}
              j_x,     & x\ne\tilde x \\
              j_x+1,   & x = \tilde x
         \end{cases}.
$$
Finally take $X'=\bigl((x_n',\tau_n')_n,(j_x')_x\bigr)$ and put $\tilde n\cdot X = X'$.

We now state and prove the properties already mentioned.

It is only assumed that the evolution stabilizes.
By that we mean that there exist $t_0$ and $\eta$ such that $\eta_t(\tilde\omega)=\eta\ \forall\ t\geqslant t_0$, so after some time all particles will be passive.
\index{stabilize}
We denote by $\tilde\Omega_s\subseteq\tilde\Omega$ the set of $\tilde\omega$ that stabilize.
Though one will immediately notice that $\PPP^\mu_M(\tilde\Omega_s)=1$, the argument hold in great generality and it involves no probability.

For $\tilde\omega\in\tilde\Omega_s$,
let $\eta_\infty(\tilde\omega)$ denote its final state
and
let $R_x(\tilde\omega)$ denote the total number of jump envelopes burned at the site $x$ during the evolution.
\index{Rx@$R_x$}

\index{commutativity}
\begin{theo}[commutativity]
\label{theo:commute}
Let $\tilde\omega = \bigl((x, \tau), F, T,N \bigr) \in\tilde\Omega$ be a realization of the
graphical construction that stabilizes.
Suppose
$\tilde\omega' = \bigl((x', \tau'), F', T',N' \bigr) \in\tilde\Omega$
has the same envelopes and initial state of $\tilde\omega$,
that is,
$(x', \tau')=(x, \tau)$ and $F'=F$.

Then $\tilde\omega'$ also stabilizes at the same state $\eta_\infty(\tilde\omega') = \eta_\infty(\tilde\omega)$.
Moreover, the total number of particles that visit each site is preserved:
$R_x(\tilde\omega') = R_x(\tilde\omega)$.
\end{theo}

The above property, besides being very useful in practice, says that only $(x,\tau)$ and $F$
matter to decide whether $\tilde\omega$ is in $\tilde\Omega_s$ and to determine $\eta_\infty$ and $R_x$.
The theorem below states that $R_x$ depends on $(x,\tau)$ and $F$ in a monotone way.
To be more precise, we write $\tilde\omega'\preceq\tilde\omega$ when:
all particles present in $\tilde\omega'$ are present in $\tilde\omega$ at the same position;
all particles that are active in $\tilde\omega'$ are also active in $\tilde\omega$;
the envelopes $F'$ can be obtained from $F$ by inserting some sleep envelopes in the sequence $(F_{x,j})_{j\in\N}$ for some sites $x$.

\index{monotonicity}
\begin{theo}[monotonicity]
\label{theo:monotone}
Let $\tilde\omega = \bigl((x, \tau), F, T,N \bigr) \in\tilde\Omega$ be a realization of the
graphical construction that stabilizes.
Suppose
$\tilde\omega'\preceq\tilde\omega$,
that is,
$\tilde\omega'$ 
has less particles than $\tilde\omega$ (and less active particles),
as well as more sleep envelopes.

Then $\tilde\omega'$ also stabilizes and the total number of particles that visit each site does not increase:
$R_x(\tilde\omega') \leqslant R_x(\tilde\omega)$.
\end{theo}

The remainder of this section is devoted to the proof of Theorems~\ref{theo:commute}~and~\ref{theo:monotone} and for another lemma that will be used in the sequel.

\textbf{Proof of Theorem~\ref{theo:commute}.}
Of course we should only consider the sequence of states of the system, regardless of the exact time they happen in the continuous-time setting.

Suppose we start at state $X_0$ and the first executed action is particle $\rho_n$, standing at site $x$ and performing the action given by $F=F_{x,j}$. It means $\tau_n=B$, $x=x_n$ and $j=j_x+1$. If $F$ is a jump envelope the particle performs that jump, 
and if it is a sleep envelope the particle becomes passive iff $x_{\tilde n}\ne x_n\ \forall\tilde n\ne n$, but in any case the envelope $F=F_{x,j}$ is burned.
This gives a new state $X_1 = C_0 \cdot X_0$, $C_0 = C(n,x,j)$, where $X_1$ differs from $X_0$ by the possible change in the position or state of particle $\rho_n$ (and the state of other particles that may be activated by that change) and by the change $j_x \mapsto j_x+1$.

So, for $C(n,x,j)\cdot X$ to make sense we must have $x_n(X)=x,j=j_x(X)+1,\tau_n(X)=B$. In this case we say that $C=C(n,x,j)$ and $X$ are \emph{compatible} and operator $C$ is called an \emph{elementary operator}.

By $\C$ we denote a (possibly finite) sequence of operators $\C=(C_1,C_2,C_3,\dots)$, which are either elementary operators or the identity operator. We say that $\C$ and $X$ are compatible if $C_k$ and $\left[(C_{k-1}\circ \cdots \circ C_1 ) \cdot X\right]$ are compatible for all $k\in\N$. We denote by $X_{k-1}$ the term in brackets and define $\C\cdot X=(X_0,X_1,X_2,\dots)$. (When $X_k=X'$ for all large $k$ we also denote $X'$ by $\C\cdot X$ and when $\C$ is finite we also denote by $\C\cdot X$ the last state in this sequence.)

For $\C=(C_1,C_2,\dots)$ and $X$ to be compatible it is necessary and sufficient that the following conditions are fulfilled.

For a fixed $\tilde n$, write $\C_{\tilde n}$ for the subsequence of $\C$ made of the elementary operators of the form $C(\tilde n,x,j)$, that is, the operators that correspond to an action performed by particle $\rho_{\tilde n}$.
Write $\C_{\tilde n}=(C_{k(1,\tilde n)},C_{k(2,\tilde n)},\dots)$, where $C_{k(l,n)}=C(n,x_{l,n},j_{l,n})$.
The sequences $\C_n$, $n=1,\dots,m$ must satisfy the two properties below.
For $l>1$, write $x=x_{k(l-1,n)}$, $j=j_{k(l-1,n)}$, $x'=x_{k(l,n)}$, $j'=j_{k(l,n)}$.
First, $\C_n$ describes a trajectory, that is, $x_{k(1,n)}=x_n(X)$ and, for $l>1$, one has
$x'=x$ if $F_{x,j}=0$ or $x'=x\pm1$ if $F_{x,j}=\pm1$.
Second, if $l>1$ and $F_{x,j}=0$, that is, particle $n$ has previously burned a sleep envelope, or if $l=1$ and $\tau_n(X)=A$, then there must be some other particle that will activate it before its next action.
More precisely, if $F_{x,j}=0$, there must be either an $n'$ such that $x_{n'}(X_{k(l-1,n)})=x$ or $k'$, $n'$, $j''$ and $x''=x\mp1$ such that $k(l-1,n)<k'<k(l,n)$, $C_{k'}=C(n',x'',j'')$ and $F_{x'',j''}=\pm1$.

For a fixed $\tilde x$ write $\C_{\tilde x}$ for the subsequence of $\C$ made of operators of the form $C(n,\tilde x,j)$, that is, the operators that correspond to an action performed at site $\tilde x$.
Write $\C_{\tilde x}=(C_{k(1,\tilde x)},C_{k(2,\tilde x)},\dots)$, where $C_{k(l,x)}=C(n_{l,x},x,j_{l,x})$.
Then the sequences must satisfy $j(1,x) = j_x(X)+1$ and, for $l>1$, $j(l,x)=j(l-1,x)+1$.

So, $\C$ and $X$ are compatible if and only if all the conditions above are satisfied for all $n$ and all $x$.
In this case the final position $x_n$ of each particle $\rho_n$ is determined by the sequence $\C_n$ and its state $\tau_n$ is $A$ if and only if the last operator in $\C_n$ corresponds to a sleep envelope and after this last operator no other particle is found in site $x=x_n$.
The final value of $j_x$ is determined by the sequence $\C_x$ and it is given by $j+1$ where the last operator in the sequence $\C_x$ is $C(n,x,j)$ for some $n$.

For $n\in\{1,\dots,m\}$, let $C=C(X,n)$ denote the action performed by particle $n$ if the system state is $X$. Namely, $C$ is the identity if $\tau_n(X)=A$ and $C=C(n,x,j+1)$ if $\tau_n(X)=B$, $x_n(X)=x$, $j_x(X)=j$.
Notice that the operation $n\cdot X$ introduced just before Theorem~\ref{theo:commute} is given by $C(X,n)\cdot X$.
For $N=(n_1,n_2,\dots)\in\{1,\dots,m\}^\N$ and $X=X_0$ a state, define $\C=\C(X,N)=(C_1,C_2,\dots)$ by $C_1=C(X_0,n_1)$, $X_1=C_1\cdot X_0$, $C_2=C(X_1,n_2)$ and so on.
Write $N\cdot X$ for the sequence $(X_0,X_1,\dots)$.
Of course in this case $\C(X,N)$ and $X$ are compatible.

Let $\pi$ be a permutation of $\{1,\dots,m\}$.
Denote $\pi(X)=\bigl( (x_{\pi^{-1}(n)}, \tau_{\pi^{-1}(n)})_{n=1}^m, (j_x)_{x\in S}\bigr)$.
Then $\pi(n\cdot X)=\pi(n)\cdot\pi(X)$.
Notice that $\eta[X]=\eta[\tilde X]$ iff $\tilde X=\pi(X)$ for some $\pi$.

For $N=(n_1,n_2,n_3,\dots)$ define $\pi(N) = \bigl(\pi(n_1),\pi(n_2),\pi(n_3),\dots\bigr)$, $\pi(C)=C(\pi(n),x,j)$ for $C=C(n,x,j)$ and $\pi(C)=C$ when $C$ is the identity, $\pi(\C)=\bigl(\pi(C_1),\pi(C_2),\pi(C_3),\dots\bigr)$ when
$\C=\bigl(C_1,C_2,C_3,\dots\bigr)$,
finally
$\pi(\textbf{X})=(\pi(X_0),\pi(X_1),\dots)$ for $\textbf{X}=(X_0,X_1,\dots)$.
It is immediate that $\pi(\C)$ and $\pi(X)$ are compatible whenever $\C$ and $X$ are so and $\pi(\C)\cdot\pi(X) = \pi(\C\cdot X)$. Also $\C(\pi(X),\pi(N)) = \pi(\C(X,N))$ and $\pi(N\cdot X) = \pi(N)\cdot\pi(X)$.

Let $\theta$ denote a shift, that is, $\theta N = (n_2,n_3,n_4,\dots)$, $\theta \textbf{X} = (X_1,X_2,X_3,\dots)$, and $\theta\C = (C_2,C_3,C_4,\dots)$; for $N = (n_1,n_2,n_3,\dots)$, $\textbf{X} = (X_0,X_1,X_2,\dots)$, $\C = (C_1,C_2,C_3,\dots)$.
Then $(\theta\C)\cdot(C_1\cdot X) = \theta(\C\cdot X)$, $(\theta N)\cdot(n_1\cdot X) = \theta(N \cdot X)$.

We say that the sequence $\textbf{X}=(X_0,X_1,X_2,\dots)$ stabilizes if there is $\tilde X$ and $k_0$ such that $X_k=\tilde X$ for all $k>k_0$ and $\tau_n(\tilde X)=A$ for all $n=1,\dots,m$. Of course if $\textbf{X}$ stabilizes at $\tilde X$ then $\theta\textbf{X}$ stabilizes at $\tilde X$ and $\pi(\textbf{X})$ stabilizes at $\pi(\tilde X)$.

The theorem will follow from
\begin{claim}
\label{claim:commute}
Let $X$ and $(F_{x,j})$ be given.
Suppose the sequence $N\cdot X$ stabilizes at $\tilde X$ for some $N\in\{1,\dots,m\}^\N$.
Then, for any $N'\in\NN_m$ there is a permutation $\pi$ such that the sequence $N'\cdot X$ will stabilize at $\pi(\tilde X)$.
\end{claim}

To prove the claim, let us write $d(X',X'') = \sum_x j_x(X'')-j_x(X')$ when $X''$ is obtained from $X'$ by applying elementary operators. In this case $d(X',X'')$ will be the number of such operators. We shall prove the statement of the claim by induction on $d(X,\tilde X)$.
For $d(X,\tilde X)=0$ we have $X=\tilde X$, so at state $X$ all particles are passive and for any sequence $N'$ we have $N'\cdot X = (X,X,X,\dots)$.

Now suppose $d(X,\tilde X)>0$. We can assume $\C = \C(X,N) = (C_1,C_2,\dots)$ contains no identity operators by suppressing the elements in the sequence $N$ that would correspond to such operators.
Take $N'=(n_1',n_2',\dots)\in\NN_m$. Since $d(X,\tilde X)>0$, there is at least one $n\in\{1,\dots,m\}$ such that $\tau_n(X)=B$.
Take $i_1=\min\{i:\tau_{n_i'}(X)=B\}\in\N$, such $i_1$ exists by definition of $\NN_m$ and the last observation ($i_1$ is the first time in the sequence $N'$ when the clock will ring for a particle that is active in $X$).
Writing $\C'=\C(X,N')=(C_1',C_2',\dots)$, we get $C_1'=C_2'=\cdots=C_{i_1-1}'=I\!d$ and $C_{i_1}' = C(\bar n,\bar x, \bar j + 1)$, where $\bar n = n_{i_1}'$, $\bar x = x_{\bar n}(X)$, and $\bar j=j_{\bar x}(X)$. We can assume $i_1=1$ by considering $\theta^{i_1-1} N'$ instead of $N'$.
Take $i^*=\min\{i:C_i=C(n,\bar x, j+1) \mbox{ for some } n,j\}$, such $i^*\in\N$ exists since $\C\cdot X$ stabilizes, $\tau_{\bar n}=B$, and $x_{\bar n}(X)=\bar x$.
Write $C_{i^*}=C(n^*, x^*, j^*+1)$.
Since $\C$ and $X$ are compatible, we must have $j^* = j_{\bar x}(X)=\bar j$ because of the definition of $i^*$.

There are two possibilities.

First case: $\bar n = n^*$. As we shall see below, if one takes
$N''=(n_{i^*},n_1,n_2,\dots, n_{i^*-1}, n_{i^*+1},n_{i^*+2},\dots)$ then it is the case that $N''\cdot X$ stabilizes and its final state is given by $\tilde X$.
Now $n_1'' = n_1'$, so $N''\cdot X = (X,\theta N''\cdot X_1'')$ and $N'\cdot X = (X,\theta N'\cdot X_1')$, where $X_1' = n_1'\cdot X = n_1''\cdot X = X_1''$.
But $d(X_1,\tilde X)=d(X,\tilde X)-1$, so we apply the induction hypothesis for $X_1'$, $\theta N''$ and $\theta N'$ to conclude that $N'\cdot X$ stabilizes at state $\pi(\tilde X)$ for some $\pi$, finishing the proof for this first case.

Second case: $\bar n \ne n^*$. Take $\pi = \pi^{\bar n,n^*}$ given by
$$
  \pi^{\bar n,n^*}(n) = 
  \begin{cases} {\bar n}, & $if $ n=n^* \cr n^*, & $if $ n=\bar n \cr n, & $otherwise$,
  \end{cases}
$$
and take $N''' = \bigl(n_1,n_2,\dots,n_{i^*-1},\pi(n_{i^*}),\pi(n_{i^*+1}),\pi(n_{i^*+2}),\dots\bigr)$. Write $N'''\cdot X = (X_0''',X_1''',X_2''',\dots)$.
Of course $X_k''' = X_k$ for $0\leqslant k < i^*$.
By definition of $i^*$, since $N$ and $X$ are compatible, it follows that $x_{n^*}(X_{i^*-1})=x_{\bar n}(X_{i^*-1})=\bar x$ and $\tau_{n^*}(X_{i^*-1})=\tau_{\bar n}(X_{i^*-1})=B$, so $\pi(X_{i^*-1})=X_{i^*-1}$. Since $\theta^{i^*-1} N''' = \pi(\theta^{i^*-1} N)$, $\theta^{i^*-1}\textbf X = (\theta^{i^*-1}N)\cdot X_{i^*-1}$ and $\theta^{i^*-1}\textbf X''' = (\theta^{i^*-1}N''')\cdot X_{i^*-1}=\pi\bigl((\theta^{i^*-1}N)\cdot X_{i^*-1}\bigr)$, we have that $\textbf{X}'''$ stabilizes at $\pi(\tilde X)$.
Also, $d\bigl(X,\pi(\tilde X)\bigr) = d\bigl(X,X_{i^*-1}\bigr)+ d\bigl(X_{i^*-1},\pi(\tilde X)\bigr)= d\bigl(X,X_{i^*-1}\bigr)+ d\bigl(X_{i^*-1},\tilde X\bigr)= d\bigl(X,\tilde X\bigr)$.
So it suffices to show that there is some $\pi'$ such that $N'\cdot X$ stabilizes at $\pi'(\tilde X)$. When one compares $N'\cdot X$ and $N'''\cdot X$, and repeats the current proof with $N'''$ instead of $N$, it will happen that the first case will hold and, as it has already been shown, there will be some $\pi''$ such that $N'\cdot X$ stabilizes at $\pi''(\pi(\tilde X))$. Take $\pi'=\pi''\circ\pi$ to complete the proof of the second case.

To finish the proof we must show that $N''\cdot X$ stabilizes at $\tilde X$ in the first case. Since $\theta^{i^*}N'' = \theta^{i^*}N$
we have $\theta^{i^*}(N\cdot X)$ = $(\theta^{i^*}N)\cdot X_{i^*}$ and $\theta^{i^*}(N''\cdot X)$ = $(\theta^{i^*}N'')\cdot X_{i^*}''$ = $(\theta^{i^*}N)\cdot X_{i^*}''$, so it suffices to show that $X_{i^*}'' = X_{i^*}$.
By definition of $n^*$ and $i^*$, we have that $x_k \ne x^*$ and $n_k \ne n^*$ for $k=1,\dots,i^*-1$, where $C_k=C(n_k,x_k,j_k)$.
Take $\bar N=(n_1,n_2,\dots,n_{i^*})$ and $\bar N''=(n_1'',n_2'',\dots,n_{i^*}'') = (n_{i^*},n_1,n_2,\dots,n_{i^*-1})$.
Let us show that $\C(X,\bar N'')=(C_{i^*},C_1,C_2,\dots,C_{i^*-1})$.
First, $n_1''=n^*$ and by the definition of $n^*$ we have $x_{n^*}(X)=x^*$, $j_{x^*}(X)=j^*$ so $C(X,n_1'')=C(n^*,x^*,j^*+1)=C_{i^*}=C_1''$. Now $X_0$ and $X_1''$ coincide except that the position and state of particle $\rho_{n^*}$ may be different, $j_{x^*}$ has increased by one unit and some other particles in $X_1''$ might be in state $B$ instead of $A$ because of a jump of $\rho_{n^*}$. So $C_1=C(X,n_1)=C(n_1,x_1=x_{n_1}(X), j_1 = j_{x_1}(X))=C(X_1'',n_1)$, since $\tau_{n_1}(X)=\tau_{n_1}(X_1'')=B$, $n_1\ne n^*$ and $x_1\ne x^*$.
Again, $X_1=C_1\cdot X$ and $X_2''=C_1\cdot X_1''$ coincide except that $j_{x^*}(X_1'')=j_{x^*}(X)+1$, $\tau_n(X_1)=A$, and $\tau_n(X_2'')=B$ for the particles $\rho_n$ that were activated by the first jump of $\rho_{n^*}$.
By repeating this procedure, we get $C(X_k'',n_k)=C(X_{k-1},n_k)$ for $k=1,\dots,i^*-1$; $X_{i^*}''$ and $X_{i^*-1}$ coincide except that $j_{x^*}(X_{i^*}'')=j_{x^*}(X_{i^*-1})+1$, $\tau_n(X_{i^*-1})=A$, and possibly $\tau_n(X_{i^*}'')=B$ for the particles $\rho_n$ that were activated by the first jump of $\rho_{n^*}$.
Now $C(X_{i^*-1},n^*)=C(n^*,x^*,j^*+1)=C_{i^*}$, since $x_{n^*}(X_{i^*-1})=x^*$ and $j_{x^*}(X_{i^*-1})=j^*$.
So $X_{i^*}=C_{i^*}\cdot X_{i^*-1}$ and $X_{i^*}$ differs from $X_{i^*-1}$ on the following features.
First, an increment at $j_{n^*}$ and the change in the position or state of $\rho_n^*$ according to $J_{x^*,j^*+1}$.
Finally, they may differ on the state of particles that would be activated after this last jump, which are all those activated at $X\mapsto X_1''$ except those that have been activated and left to some other site during the course of $X\mapsto X_1 \mapsto\dots\mapsto X_{i^*-1}$.
This implies $X_{i^*}=X_{i^*}''$.
This finishes the proof of Claim~\ref{claim:commute} and of the theorem.
$\hfill\square$
\\

\textbf{Proof of Theorem~\ref{theo:monotone}.}
We have to prove the statement of the theorem when $\tilde\omega'\in\tilde\Omega$ is obtained from $\tilde\omega\in\tilde\Omega_S$ by modifying the following features: some particles may be removed at the initial state, particles that are kept may start passive instead of active, and sleep envelopes may be inserted.
We can perform such modifications one by one, so it suffices to prove the statement when either only one particle is removed, or one particle starts passive instead of active, or else one sleep envelope is inserted at a given site.

We continue to use the objects, notations and results of the proof of Theorem~\ref{theo:commute}.

Given a sequence $\C$ compatible with a state $X$ there is (at least) one sequence $N$ such that $\C=\C(X,N)$. We shall denote by $N(\C)$ the smallest one. For $\C=(I\!d,\dots,I\!d,C(n_1,x_1,j_1), \dots, I\!d,\dots,C(n_2,x_2,j_2),\dots)$, remove the identity operators, obtaining $\C'=(C(n_1,x_1,j_1),C(n_2,x_2,j_2),C(n_3,x_3,j_3),\dots)$. Of course $\C$ and $\C'$ are compatible with the same states $X$ and furthermore $\C\cdot X$ and $\C'\cdot X$ are the same except that some positions of the later sequence may appear repeated several times on the former one, in particular $\C\cdot X$ stabilizes at $\tilde X$ if and only if $\C'\cdot X$ does. We define $N(\C)=(n_1,n_2,n_3,\dots)$ in this case and it follows that $N(\C)\cdot X = \C'\cdot X$.

Let us consider the first case where the envelope configuration is changed so that $(F'_{x,j})$ is given by inserting one sleep envelope at one given site in $(F_{x,j})$. Take $N\in\NN_m$ and consider the sequence $N\cdot_{F'} X$. If the sequence stabilizes before even burning this extra envelope, the result follows trivially. So suppose otherwise there is some $k,\tilde n$ such that $C(X_{k-1},n;F')=C(n,\tilde x,\tilde j+1)$, where $F'_{\tilde x,\tilde j+1}$ is the inserted envelope. Then $X_k$ is given by $X_{k-1}$ except that the counter $j_{\tilde x}$ increases by one unit and that maybe $\tau_n(X_k)=A$. The first difference does not violate the state of the theorem since, for its purpose, burned sleep envelopes are not counted. The other difference, that the particle is passive in the evolution ruled by $F'$ when it should be active in the evolution ruled by $F$, is exactly the second case discussed below.

Second case: $x_n(X)=x_n(X')$ for all $n$, $\tau_n(X)=\tau_n(X')$ for all $n=1,\dots,m$ except for the last particle $\rho_m$, for which $\tau_m(X)=B,\tau_m(X')=A$ (we have assumed for simplicity that the difference is for particle $\rho_m$). We have to prove that for some $N\in\NN_m$, $N\cdot X'$ stabilizes at $\tilde X'$ and $j_x(\tilde X')\leqslant j_x(\tilde X)\ \forall\ x$, under the hypothesis that $N\cdot X$ stabilizes at $\tilde X$.
So fix $N\in\NN_m$.
Take $\bar\C=\C(X,N)$ and $\bar N=N(\bar\C)$.
Then $\bar N=(n_1,\dots,n_{\bar k})$, where $\bar k=d(X,\tilde X)$.

Take $\bar N'=(1,2,\dots,m-1)$ and $X_{(1)}'=\bar N'\cdot X'$, where by abuse of notation we call $N'\cdot X'$ the last term in the $m$-tuple.
If $X_{(1)}'=X'$ it means $\tau_n(X')=A\ \forall\ n=1,\dots,m-1$ and, since $\tau_m(X')=A$, we have by Theorem~\ref{theo:commute} that $X'=\tilde X'$. Then of course $j_x(\tilde X) \geqslant j_x (X) = j_x(X') = j_x(\tilde X')$ and we are done. If $\tau_m(X_{(1)}')=B$, since $m\not\in\bar N'$ we have $X_{(1)}' = \bar N'\cdot X$, so take $N' = (\bar N',N)\in\NN_m$, it follows by Theorem~\ref{theo:commute} that $N'\cdot X$ and $N'\cdot X'$ coincide after the first $m$ coordinates, therefore $N'\cdot X'$ stabilizes at $\tilde X$ and we are done.
So suppose $X_{(1)}' \ne X'$ and $\tau_m(X_{(1)}')=A$.
Then $d_1=d(X',X_{(1)}') \geqslant 1$.
Take $X_{(2)}' = \bar N'\cdot X_{(1)}'$.
If $X_{(2)}'=X_{(1)}'$ it means $\tau_n(X_{(1)}')=A\ \forall\ n=1,\dots,m-1$ and, since $\tau_m(X_{(1)}')=A$, we have by Theorem~\ref{theo:commute} that $X_{(1)}'=\tilde X'$. Then of course $j_x(\tilde X) \geqslant j_x (\bar N'\cdot X) = j_x(\bar N'\cdot X') = j_x(\tilde X')$ and we are done.
If $\tau_m(X_{(2)})'=B$, since $m\not\in\bar N'$ we have $X_{(2)}' = \bar N'\cdot X$, so take $N' = (\bar N',\bar N',N)\in\NN_m$, it follows by Theorem~\ref{theo:commute} that $N'\cdot X$ and $N'\cdot X'$ coincide after the first $2m-1$ coordinates, therefore $N'\cdot X'$ stabilizes at $\tilde X$ and we are done.
Now if we suppose $X_{(2)}' \ne X_{(1)}'$ we have $d(X_{(1)}',X_{(2)}')\geqslant 1$, so $d_2 = d(X',X_{(2)}')\geqslant 2$.
Now either we will continue to define $X_{(3)}',X_{(4)}',\dots$ indefinitely, or this construction will halt at some $X_{(k)}$.
All we need to do is to rule out the first possibility.
For that purpose we are going to prove that $d_k \leqslant \bar k$; this suffices, since by the above construction we have $d_k\geqslant k$.
For simplicity we consider $k=2$.
Let $\C_2 = \C\bigl(X',(\bar N',\bar N')\bigr)$, $N_2 = N(C_2)$ and $\bar\C_2=\C(X',N_2)$.
Briefly, $\bar\C_2$ is just the $d_2$-tuple of elementary (non-trivial) operators that take $X'$ to $X_{(2)}''$.
Since the only difference between $X$ and $X'$ is that at state $X$ the particle $\rho_m$ is active, $\bar\C_2$ being compatible with $X'$ implies it is compatible with $X$ and $\bar\C_2\cdot X = \bar N_2\cdot X$. By Theorem~\ref{theo:commute}, $(\bar N_2,N)\cdot X$ stabilizes at $\tilde X$, so $d(X,\tilde X)\geqslant d_2$, finishing the proof of the second case.

Third case: $X=\bigl((x_n,\tau_n)_{n=1}^m,(j_x)_{x\in S}\bigr)$, and $X'=\bigl((x_n,\tau_n)_{n=1}^{m-1},(j_x)_{x\in S}\bigr)$ (we have assumed for simplicity that the missing particle is $\rho_m$).
As before, we have to prove that for some $N\in\NN_{m-1}$, $N\cdot X'$ stabilizes at $\tilde X'$ and $j_x(\tilde X')\leqslant j_x(\tilde X)\ \forall\ x$, under the hypothesis that $N\cdot X$ stabilizes at $\tilde X$.
We fix $N\in\NN_m$ and take $\bar\C=\C(X,N)$ and $\bar N=N(\bar\C)$.
Then $\bar N=(n_1,\dots,n_{\bar k})$, where $\bar k=d(X,\tilde X)$.
We consider the same $\bar N'$, and keep on defining $X_{(1)}',X_{(2)}',\dots$ as before, and now the condition for this construction to stop is when $X_{(k)}'=X_{(k-1)}'$, which means $\tau_n(X_{(k-1)}')=A$ for $n=1,\dots,m-1$.
Again, if $\bar\C_k$ is compatible with $X'$, it is also compatible with $X$, with some abuse of notation since now we are considering elementary operators acting on different state spaces.
So by the same argument as for the second case, $k \leqslant d_k = d(X',X_{(k)}') \leqslant \bar k$ and for some finite $k$ we have $\tilde X'=X_{(k)}' = (C_{d_k}'\circ\cdots\circ C_1')\cdot X'$ satisfies $\tau_n(X_{(k)})=A$ for $n=1,\dots,m-1$.
Again in this case $X_{(k)}=(C_{d_k}'\circ\cdots\circ C_1')\cdot X$ is well defined and it differs from $\tilde X'$ for the presence of an extra particle and because for some $n$ it might be that $\tau_n(X_{(k)})=B$ and $\tau_n(\tilde X')=A$ but $j_x(\tilde X') = j_x(X_{(k)})$.
Now $j_x(X_{(k)}) \leqslant j_x(\tilde X)$ since, by Theorem~\ref{theo:commute}, $N\cdot X_{(k)}$ stabilizes at $\tilde X$. This completes the proof of the third case.
$\hfill\square$

A different coupling allows the comparison between different initial conditions for a fixed time.
We shall need the fact below for a mild technical question of exchanging the order of a certain limit in Section~\ref{sec:phasetransition}.
For fixed time $t>0$ let $R_x^t$ denote the number of jump envelopes burned at site $x$ up to time $t$.
\index{Rxt@$R_x^t$}
\begin{lemma}
\label{lemma:fixedt}
Let $P,P'$ denote probability distributions on $\Sigma''$ such that $P\succeq P'$.
Take $\PPP$ and $\PPP'$ as the corresponding laws on $\tilde\Omega$ according to the construction described above.
For each $x\in S$, $t>0$, and $r\in\N$,
$\PPP\bigl(R_x^t\geqslant r\bigr)\geqslant\PPP'\bigl(R_x^t\geqslant r\bigr)$.
\index{monotonicity}
\end{lemma}

We shall just highlight the main steps of the proof:

Let $\bigl(\eta^{A,B}_0(x)\bigr)_{x}$ and $(F_{x,j})_{x,j}$ be given.
We shall construct $(\eta_t)_{t\geqslant0}$ in a slightly different way.
Take $(\tilde T_{x,j})_{x,j}$ as i.i.d. $\exp(1+\lambda)$ random variables and $T_{x,j} = \sum_1^j \tilde T_{x,i}$.
Let $J^0_x = 0$ for $x\in S$.
Write $D_x(t) = \int_0^t \eta_s^B(x) ds$.
Take $t_1 = \inf\{t: D_x(t) = T_{x,j=1}(t) \mbox{ for some } x\in S\}$ and take $x_1$ as the (a.s. unique) $x$
such that $D_x(t_1) = T_{x,j=1}(t_1)$.
Let $\eta_s = \eta_0$ for $s < t_1$ and $\eta_{t_1}$ be given by the new state obtained from $\eta_0$ by opening envelope $F_{x_1,1}$ and set $J^1_x = \delta_{x_1}(x)$.

Suppose $\eta_s$ has been defined for $s\in[0,t_k]$.
Take $t_{k+1} = \inf\{t: D_x(t) = T_{x,J^k_x+1}(t) \mbox{ for some } x\in S\}$ and take $x_{k+1}$ as the (a.s. unique) $x$ such that $D_x(t_{k+1}) = T_{x,J^k_x+1}(t_{k+1})$.
Let $\eta_s = \eta_{t_k}$ for $t_k \leqslant s < t_{k+1}$ and $\eta_{t_{k+1}}$ be given by the new state obtained from $\eta_{t_k}$ by opening the envelope $F_{x_{k+1},J^k_x+1}$. Set $J^{k+1}_x = J^k_x + \delta_{x_{k+1}}(x)$.

Now by a routine argument, similar to that used to prove super-additivity for first passage percolation, one proves the stochastic domination claimed on the lemma.
$\hfill\square$

\section{The phase transition}
\label{sec:phasetransition}

In this section we present the proof of Theorem~\ref{theo:fixation}.
The proof strongly relies on the commutativity and monotonicity properties shown in Section~\ref{sec:flag} and approximations by finite boxes~(\ref{eq:finitebox}) will be important.

We shall first give a simple characterization of phase transition in terms of burned envelopes.
From this we immediately prove the existence of $\mu_c\in[0,\infty]$ and that $\mu_c$ is increasing in $\lambda$.
Then we show that $\mu_c\leqslant 1$ for one-dimensional simple random walks.
Finally, in the remainder of the section we give the proof that $\mu_c\geqslant\frac\lambda{1+\lambda}$ for simple random walks, which is the main result.

By translation invariance, the phenomenon of fixation is equivalent to the fact that a.s. there is fixation at the origin, that is, for large enough time $\eta^B_t(0)=0$.
Consider the random variable $R_0^t$ given by the number of times a particle jumps out of the origin up to time $t$.
As $t$ increases one has $R_0^t\uparrow R_0$ and fixation is equivalent to $R_0<\infty$.
Define the event $A_r^t = [R_0^t\geqslant r]$.
Now $A_r^t \uparrow_{t\to\infty} A_r \downarrow_{r\to\infty} A_\infty$, where $A_r = [R_0 \geqslant r]$,
and fixation is equivalent to $\PP^\mu(A_\infty)=0$.

Note that by~(\ref{eq:finitebox}) one has
$$
  \PP^\mu(A_r^t) = \lim_{M\to\infty} \PP^\mu_{M}(A_r^t),
$$
so fixation means
\begin{equation}
\label{eq:unifconvergence}
  \PP^\mu_{M}(A_r) \mathop{
  \longrightarrow}_{r\to\infty} 0
  \qquad\mbox{uniformly in $M$.}
\end{equation}
Now, within the framework of Section~\ref{sec:flag}, the event $A_r$ means that at least $r$ jump envelopes are burned at the origin.

It is immediate that (\ref{eq:unifconvergence}) is sufficient for fixation. 
To see why it is necessary we remark that $\PP^\mu_M(A_r^t)$ is increasing in both $M$ and $t$, so switching limits we get
$\PP^\mu_{M}(A_r^t) \leqslant \PP^\mu(A_r)$.
Monotonicity in $t$ is obvious and monotonicity in $M$
follows from Lemma~\ref{lemma:fixedt}.

Now suppose (\ref{eq:unifconvergence}) holds for some $\mu$.
By Theorem~\ref{theo:monotone}, the same limit will hold for any $\mu'<\mu$.
So one can take $\mu_c=\sup\{\mu\geqslant0:\mbox{limit (\ref{eq:unifconvergence}) holds}\}$.
(Notice that it does not exclude the possibilities $\mu_c=0$ or $\mu_c=\infty$.)
It follows from Theorem~\ref{theo:monotone} that $\mu_c$ is nondecreasing in $\lambda$.
All this holds in any dimension $d\geqslant1$.

To prove $\mu_c < \infty$ in one dimension we argue as follows.
Suppose $\mu>1$.
With high probability there will be an amount of particles of order $\mu M$ in a given large region (say, in $[-M,0]$).
Since no more than $M$ particles can stay passive in this region, an excess of particles (at least about $(\mu-1)M$) will have to exit such region, and because of the topological constraints of the one-dimensional lattice, they must exit either by its left border or by its right border, that is, the origin.
This implies that whp at least $(\mu-1)M/2$ particles will eventually pass by the origin and therefore $R_0 = \infty$ almost surely, so the system does not fixate.

This actually proves that $\mu_c \leqslant 1$.
We have implicitly assumed that the random walks are of nearest-neighbor type but with slight modifications in the argument one easily shows that $\mu_c \leqslant 1$ for any finite-range random walk.

Let us now show that $\mu_c \geqslant \frac \lambda{1+\lambda}$ for nearest-neighbor walks. With simple adaptations the same argument proves that $\mu_c>0$ for finite-range random walks.

First, we give a heuristic proof. The rigorous argument follows next.

\index{trap}
We shall build a global trap for the particles using disjoint sites to guarantee that they cannot collaborate (activate each other) to overcome the trap.
The density of individual traps one is going to find per site will be given by
$\frac \lambda {1+\lambda}$
because this is the probability that between two consecutive jump envelopes there is at least one sleep envelope.
If the density of particles is smaller, that is, $\mu<\frac \lambda {1+\lambda}$, then on a large scale the system will be trapped whp.
In order to be able take advantage of this large-scale behavior, we choose a large enough region $\Lambda$ around the origin and we sweep $\Lambda$.
Once this is done the particles cannot come close to the origin again because, due of the large-scale domination of the traps over the particles, the probability that the trap fails will be small.

\index{sweep}
To `sweep' the region $\Lambda$ we let the particles move inside this area until they reach its boundary.
The particles will be swept one by one, this is possible because of Theorem~\ref{theo:commute}.
None of them will become passive until this step is finished (here we use Theorem~\ref{theo:monotone}).
During this procedure the probability that more than $r$ particles cross the origin is small if $r$ is taken large enough.

We prove now that whp no particle will come back to the origin after $\Lambda$ is swept.
Place an imaginary barrier at the origin and pick the particle that is closest to this barrier.
Let this particle jump (Theorem~\ref{theo:commute} again) until it hits the origin (or some fixed site far away).
While we make the particle jump we are ignoring that each envelope this particle is opening could have been a sleep envelope
(Theorem~\ref{theo:monotone} again).
Then we start coming back along this path while looking for sleep envelopes that might have been missed, until we succeed to find the first one.

The set of envelopes we have discovered so far has the following property.
Suppose the label sequence were exactly as we had assumed, and the sleep envelopes we have ignored were indeed absent.
Then all the particles that started in $\Lambda$ would move until they reach the border of $\Lambda$.
Next, one of these particles would start jumping until it opens the sleep envelope we have managed to find.
Therefore this particle would become passive before reaching the barrier.

Now place a new imaginary barrier at the site where the trap for the first particle has been set.
Pick the particle that is closest to this new barrier.
Again let it jump until it hits the new barrier.
Then trace back its path looking for sleep envelopes that had been ignored, until the first one is found.
As in the previous paragraph, the evolution we are building makes the particles evacuate $\Lambda$, attracts the first particle to the first trap and now it attracts this second particle to this other trap without touching the origin and without activating the first particle.

Again place a new barrier where the second particle should be trapped and keep repeating this procedure until a trap has been set for each of the particles.
The global configuration so obtained is enough to guarantee that all the particles will be passive without coming back to the origin.

As mentioned in the previous section, this evolution may have nothing to do with that governed by $\PP^\mu_{M}$. Nevertheless the real evolution (with a random sequence of labels and a few more sleep envelopes) will have the same value for $R_0$ by Theorem~\ref{theo:commute}, or maybe even smaller by Theorem~\ref{theo:monotone}.

It is possible that this procedure fails. It will fail if, at some stage, say the $n$-th stage, we do not manage to find a sleep envelope between the position of the $n$-th particle (after sweeping) and the $(n-1)$-th barrier.
But the probability that this will happen is very small because the density of traps is smaller than that of particles and $\Lambda$ has been chosen large enough.

The remainder of this section is devoted to formalizing the algorithm just described.

For each $x\in S$ split the envelope sequence $(F_{x,j})_j$ into two sequences: one consisting on the jump envelopes and one consisting on the sleep envelopes.
The first sequence is given by $\tilde F_{x,k} = F_{x,j(x,k)}$ where $j(x,0)=0$, $j(x,k)=\min\{j>j(x,k-1):F_{x,j}=\pm 1\}$.
So $\tilde F_{x,k}$ is the $k$-th jump envelope at $x$ and between the $(k-1)$-th and $k$-th jump envelopes there are $D_{x,k}=[j(x,k)-j(x,k-1)-1]$ sleep envelopes in $(F_{x,j})_j$.
Of course $(F_{x,j})_j$ is in bijection with $(\tilde F_{x,k},D_{x,k})_{k}$.
Also, inserting a sleep envelope in $(F_{x,j})_j$ means increasing the value of $D_{x,k}$ for some $k$.

It is routine to check that $(\tilde F_{x,k})_{x,k}$ are i.i.d. having value $+1$ with probability $p$ and $-1$ with probability $q=1-p$.
Also $(D_{x,k})_{x,k}$ are i.i.d. geometric variables that are positive with probability $\frac\lambda{1+\lambda}$.
We finally remark that $(\tilde F_{x,k})_{x,k}$ and $(D_{x,k})_{x,k}$ are independent.

We are going to show that the system fixates if $\mu<\frac\lambda{1+\lambda}$.
We first present the proof for the symmetric case $p=q=\frac12$ and after we give the modifications in the case $p=1-q>\frac12$.

Let $\epsilon>0$. Take $K$ with the following property: if $S_n$ and $S_n'$ denote the sum of the $n$ first random variables of an i.i.d. sequence distributed respectively as Bernoulli$(\frac\lambda{1+\lambda})$ and Poisson$(\mu)$) then the probability that $S_{n+1}'>S_n$ for some $n\geqslant K$ is less than $\epsilon$.
The reason for this choice of $K$ will be clear later.

In the sequel we shall consider the evolution given by the envelope sequence $(F'_{x,j})$ instead of $(F_{x,j})$, where $(F'_{x,j})$ is obtained from the same $(\tilde F_{x,j})$ but $(D'_{x,k})\equiv 0$.
It means that we take the configuration $\tilde\omega'=\bigl((x_n,\tau_n)n,F',T,N'\bigr)\succeq\tilde\omega$.
Here $F'$ is obtained by removing all the sleep envelopes and $N'$ will be constructed in an algorithmic way.
Of course for this strange configuration $\tilde\omega'$ the particles will never become passive since all sleep envelopes have been removed.
So later on we shall put back a few sleep envelopes present in $\tilde\omega$ to construct an intermediate scenario $\tilde\omega''\succeq\tilde\omega$ such that the particles become appropriately trapped in the evolution given by $\tilde\omega''$.

\label{sweeping}
Let $m'$ be the number of particles in $\tilde\omega$ and relabel them $\rho_1,\dots,\rho_m,\dots,\rho_{m'}$ in a way that $|x_n|\leqslant K$ for $n=1,\dots,m$ and $|x_n|> K$ for $n=m+1,\dots,m'$.
Let $X_0=\bigl((x_n,\tau_n)_n,(j_x)_x\bigr)$ be the initial state, with $j_x\equiv0$.
Take $n_1 = 1$ and $X_1 = n_1\cdot X_0$.
If $x_1(X_1)=\pm K$ we say particle $\rho_1$ has been swept and we move to particle $\rho_2$; otherwise we take $n_2=1$, $X_2=n_2\cdot X_1$ and so on.
Once we have $x_1(X_{k_1})=\pm K$ we move to particle $\rho_2$.
Then we take $n_{k_1+1}=2$, $X_{k_1+1} = n_{k_1+1}\cdot X_{k_1}$ and keep doing this, thus obtaining $X_0,X_1,X_2,\dots$ until $x_2(X_{k_2})=\pm K$ for some $k_2$.
After this we take $n_{k_2+1}=3$, $X_{k_2+1} = n_{k_2+1}\cdot X_{k_2}$ and keep doing this until $x_3(X_{k_3})=\pm K$ for some $k_3$.
At this point particles $\rho_1$, $\rho_2$ and $\rho_3$ have been swept to the boundary of $[-K,K]$.
We repeat this procedure until we have swept all of $\rho_1,\dots,\rho_m$, thus getting $(n_1,n_2,n_3, \dots,n_{k_1},n_{k_1+1},\dots,n_{k_2},\dots,n_{k_m})$ and $(X_0,X_1,\dots,X_{k_m})$.

After this step the region $[-K,K]$ is clean and $|x_n(X_{k_m})| \geqslant K$ for all $n=1,\dots,m'$.
Moreover, the number of particles at each boundary $\eta^{B}(\pm K)$ with $\eta=\eta[X_{k_m}]$ is Poisson$(\mu(K+\frac{1}2))$ and is independent of $(\eta^B(x))_{|x|>K}$.
Also, $j_x(X_{k_m})$ is a.s. finite, so for there is $r_0$ such that the probability that
$j_x(X_{k_m})\geqslant r_0$ at the origin $x=0$ is less than $\epsilon$.

We are now going to build a trap to show that after sweeping no particle can come back to the origin.

To start this second step relabel the particles in $X_{k_m}$, that is, take
$\tilde X_{0}= \pi (X_{k_m})$
for some $\pi$ such that $K\leqslant x_1\leqslant x_2\leqslant \cdots\leqslant x_{m''}$
and $-K\geqslant x_{m''+1}\geqslant x_{m''+2}\geqslant\cdots\geqslant x_{m'}$ in $\tilde X_0$.

Let $a_0 = 0$, $b_0 =\lceil (3 + 1/\lambda) M\rceil$ and $W_0 = \emptyset$.
We shall build a sequence of intervals $[a_n,b_n]\supseteq[a_{n-1},b_{n-1}]\supseteq\cdots$ where the particles $\rho_n,\dots,\rho_{m''}$ will remain confined.
Finally $W_0\subseteq W_1\subseteq\cdots$ will be the sequence of traps where each particle will become passive in a modified configuration $\tilde\omega'' \preceq \tilde\omega$.

Take
$n_1' = 1$, $\tilde X_{1}=n_1'\cdot \tilde X_{0}$,
$n_2' = 1$, $\tilde X_{2}=n_2'\cdot \tilde X_{1}$,
$n_3' = 1$, $\tilde X_{3}=n_3'\cdot \tilde X_{2}$,
and so on until for some $l_1>0$ the position $x_1(\tilde X_{l_1})$ of particle $\rho_1$ is either $a_0$ or $b_0$. 

The idea is to set up a trap for $\rho_1$ strictly between the origin and $\rho_1$'s initial position in a way that (i) this particle does not reach the origin and (ii) this operation consumes as little space as possible.
The first requirement is exactly what we want to prove and the second requirement will be important because when setting up the trap for the other particles $\rho_2,\rho_3,\dots$ we want to guarantee that none of them will activate $\rho_1$ again.

Suppose $x_1(\tilde X_{l_1})=a_0$.
The set $T_1 =\{ a_0+1,\dots,x_1(\tilde X_0)-1 \}$ contains the possible sites where the trap for $\rho_1$ will be set up.
For $x\in T_1$ let $k_1(x)=j_x(\tilde X_{l_1})$.
Now take $\tilde x_1 = \min\{x\in T_1:D_{x,k_1(x)}>0\}$.
If this set is empty we have failed in finding a trap for $\rho_1$ and we stop the construction.
Otherwise we declare the first trap successful.
The construction of the first trap is complete on defining $a_1=\tilde x_1$, $b_1=b_0$
and $W_1 = W_0\cup \{(\tilde x_1,k_1(\tilde x_1))\}$.

In the case $x_1(\tilde X_{l_1})=b_0$ we take $T_1=\{M+1,\dots,b_0-1\}$, $\tilde x_1 = \max\{x\in T_1: D_{x,k_1(x)}>0\}$, $a_1=a_0$, $b_1=\tilde x_1$ and the same expression for $W_1$.

Now let us trap $\rho_2$.
Again $\rho_2$ should neither touch the origin nor activate $\rho_1$ and the trap for $\rho_2$ should consume as little space as possible.
Take
$n_{l_1+1}' = 2$, $\tilde X_{l_1+1}=n_{l_1+1}'\cdot \tilde X_{l_1}$,
$n_{l_1+2}' = 2$, $\tilde X_{l_1+2}=n_{l_1+2}'\cdot \tilde X_{l_1+1}$,
$n_{l_1+3}' = 2$, $\tilde X_{l_1+3}=n_{l_1+3}'\cdot \tilde X_{l_1+2}$,
and so on until for some
$l_2>0$
the position 
$x_2(\tilde X_{l_2})$
of particle $\rho_2$ is either $a_1$ or $b_1$. 

If $x_2(\tilde X_{l_2})=a_1$ we take $T_2 =\{ a_1+1,\dots,x_2(\tilde X_0)-1 \}$,
for $x\in T_2$ we set $k_2(x)=j_x(\tilde X_{l_2})$, then take $\tilde x_2 = \min\{x\in T_2:D_{x,k_2(x)}>0\}$.
Again if such set is empty we have failed in finding a trap for $\rho_2$ and we stop the construction.
Otherwise we declare the second trap successful.
Finally we define $a_2=\tilde x_2$, $b_2=b_1$ and $W_2 = W_1\cup \{(\tilde x_2,k_2(\tilde x_2))\}$.
If $x_2(\tilde X_{l_2})=b_1$ we take $T_2=\{M+1,\dots,b_1-1\}$, $\tilde x_2 = \max\{x\in T_2: D_{x,k_2(x)}>0\}$, $a_2=a_1$, $b_2\tilde x_2$ and the same expression for $W_2$.

In this way we construct $[a_0,b_0],W_0,[a_1,b_1],W_1,\dots,[a_{m''},b_{m''}],W_{m''}$.
Starting with $\tilde a_0 =-b_0$, $\tilde b_0=0$ and $\tilde W_0 = \emptyset$
we perform a similar construction for $[\tilde a_0,\tilde b_0],\tilde W_0,[\tilde a_1,\tilde b_1],\tilde W_1,\dots,[\tilde a_{m'-m''},\tilde b_{m'-m''}],\tilde W_{m'-m''}$.

Suppose all the traps have been successful.
Take $W = W_{m''} \cup W_{m'-m''}$, the total set of traps. $W$ has exactly $m'$ elements, one trap corresponding to each particle.
By construction $D_{x,k}>0$ for all $(x,k)\in W$.
Take $D_{x,k}'' = \I_{(x,k)\in W}$ so that $D_{x,k}'' \leqslant D_{x,k}$ and consider the envelope sequence $(F_{x,j}'')$ obtained from $(\tilde F_{x,j})$ and $(D_{x,k}'')$.
Write $n_{k_m+i}=\pi^{-1}(n_i')$ for $i=1,\dots,l_{m'}$ and take $N'' = N' = (n_1,n_2,\dots,n_{k_m+l_{m'}})$.
Let $(X_0'',X_1'',\dots)$ represent the evolution of $\tilde\omega'' = \bigl((x_n,\tau_n)_n,F'',T'',N''\bigr) \succeq \tilde\omega$.
We leave for the reader the task of checking 
that $x_n(X_i'') \ne 0$ for all $i\geqslant k_m$ and that the evolution $(X_1'',X_2'',X_3'',\dots)$ stabilizes at $\pi^{-1}(\bar X)$, where $x_n(\bar X) = \tilde x_n$.
In other words, the particles leave $(-K,+K)$ in the first step and during the second step they get stuck at the traps at $\tilde x_n$ without being at the origin again.

The procedure may stop at some point, namely if one fails to set a certain trap.
Given that the $n$-th trap was successful and it was placed at site $\tilde x_n$, the probability that the $(n+1)$-th trap will be at site $\tilde x_n+1 + k$ is given by $(\frac{1}{1+\lambda})^k \frac{\lambda}{1+\lambda}$. By the definition of $K$, the probability that one of the $m''$ first traps is not successful is less than $\epsilon$. Similarly, the probability that one of the $m'-m''$ other traps fails is also bounded by $\epsilon$. By the choice of $r_0$ we have that $\PP^\mu_M(A_r) < 3\epsilon$ for any $r\geqslant r_0$, where $r_0$ does not depend on $M$, finishing the proof for the simple symmetric random walk.

For the simple asymmetric random walk the proof has slight adaptations that we indicate in the following three paragraphs.

First we explain how to define a suitable $K$.
Suppose $p>q$, so the random walks are biased to the right.
For $k > 0$, consider the following system of random walks.
Start with Poisson$(2 k \mu)$ number particles at site $x=k$ and, independently of that, Poisson$(\mu)$ particles at every site $x>k$.
Then let all these particles perform independent random walks with probability $p$ of jumps to the right and $q=1-p$ of jumps to the left.
Let $x_0$ be the random variable given by the minimum position among all such paths.
It is easy to see that $x_0 \to \infty$ in probability as $k\to\infty$.
So choose $k_1$ such that $x_0$ is positive with probability at least $1-\epsilon$.
We are now going to choose $k_2$.
Let $S_n$ (resp. $S_n'$) denote the sum of the $n$ first random variables of an i.i.d. sequence distributed as Bernoulli$(\frac\lambda{1+\lambda})$ (resp. Poisson$(\mu)$).
Then there is $k_2>0$ such that the probability that $S_{n+k}'>S_n$ for some $n\geqslant 1$ is less than $\epsilon$ if $k\geqslant k_2$.
Take $K=k_1\vee k_2$.

There is also a change in the sweeping procedure defined on p.~\pageref{sweeping}.
Instead of stopping to move a particle until it reaches $\pm K$ we move it until it hits $+K$.
This will happen almost surely since the random walk is biased to the right.

The second step, namely that of finding traps, is the same.
Once $K$ is chosen, the choice of $r$ is also analogous.
So, the probability that this procedure will fail will again be bounded by $3\epsilon$.

This finishes the proof of Theorem~\ref{theo:fixation}.
$\hfill\square$

\section{Concluding remarks}
\label{sec:concluding}

The proof that $\mu_c \geqslant \frac \lambda{1+\lambda}$ for simple random walks presented in Section~\ref{sec:phasetransition} can be easily adapted for finite-range walks to show that $\mu_c > 0$, though the authors conjecture that $\mu_c \geqslant \frac \lambda{1+\lambda}$ should still hold in this case.

The initial distribution being Poisson with parameter $\mu$ plays no special role in the proof.
Indeed the proof of the lower bound for $\mu_c$ works with any family $\{P_\mu\}_{0\leqslant\mu<\infty}$ of distributions on $\{0,1,2,\dots\}$ parametrized by the first moment.
For the uniqueness of the phase transition we only need the natural assumption that $P_{\mu}\succeq P_{\mu'}$ when $\mu\geqslant\mu'$, that is, the bigger $\mu$, the more particles.

Yet we conjecture that this phase transition is actually universal with respect to the distribution.
There should be a $\mu_c$ such that for any spatially-ergodic initial distribution the system fixates if $E[\eta^{AB}(0)]<\mu_c$ and does not fixate if $E[\eta^{AB}(0)]\geqslant\mu_c$, $E[\eta^{A}(0)]>0$.
In the later case the system should converge to an ergodic distribution that depends only on $\mu$.

Some questions remain open.
We do not know how to answer such questions, but nevertheless we strongly believe that the same scheme developed here could be extremely useful for this purpose.

To complete the qualitative description of this model, it is still open to rule out two possibilities, which are conjectured not to hold.
These possibilities are physically absurd, but so far there is no mathematical proof that they do not hold.
We have shown that $\mu_c$ increases as $\lambda$ increases, that $0<\mu_c\leqslant 1$ for all $\lambda$ and that $\mu_c\to1$ as $\lambda\to\infty$.
What remains to be proven is that $\mu_c<1$ for all $\lambda$ and $\mu_c\to0$ as $\lambda\to0$.
Strikingly enough, even the fact that for some $\lambda>0$ and $\mu<1$ the system does not fixate (which would be implied by either of these) remains unproven.

One can also consider $\mu$ fixed and study the phase transition at $\lambda_c$, as was exposed at the introductory chapter.
Theorem~\ref{theo:fixation} then says that $\lambda_c<\infty$ for $\mu<1$, that $\lambda_c=\infty$ for $\mu\geqslant1$ and that $\lambda_c\to0$ as $\mu\to0$.
What remains open then is that $\lambda_c\to\infty$ as $\mu\nearrow1$ and that $\lambda_c>0$ for all $\mu>0$.

In dimensions $d\geqslant2$, the fact that $\mu_c>0$ is a very interesting open problem.
Whereas for $d=1$ the trick of letting the particles evolve until hitting barriers and then setting traps worked well, for higher dimensions different ideas will be necessary.

The proof that $\mu_c\leqslant 1$ for dimensions $d\geqslant2$ was recently given by Shellef~\cite{shellef08}, his proof is quite simple within the framework we consider here.
The fact that $\mu_c < 1$ still remains as an open question in any dimension.

The open problem that is, in my opinion, the most interesting but the hardest to approach is that of whether there is fixation at $\mu=\mu_c$.
Unlike the most common models that exhibit phase transition, such as percolation or Ising model, we believe that the critical behavior of this model bears the most important features from the supercritical phase.
Namely, at $\mu=\mu_c$ the system should not fixate.
The only case considered so far~\cite{hoffman07} is that of one-dimensional, totally asymmetric, nearest-neighbor random walks.

\chapter{The Broken Line Process}
\label{chap:brokenline}

In this chapter we introduce the continuum broken line process and derive some of its properties.
This process is similar in spirit to the Arak-Surgailis process.
Its discrete version, introduced by Sidoravicius, Surgailis and Vares, is presented first.
A natural generalization to a continuous object living on the discrete lattice is then proposed and studied.
The broken lines are related to the Young diagram and the Hammersley process and are useful for computing last passage percolation values and finding maximal oriented paths.

For a class of passage time distributions there is a family of boundary conditions that make the process stationary and reversible.
For such distributions there is a law of large numbers and the process can be extended to the infinite lattice.

One application is a simple proof of the explicit law of large numbers for last passage percolation known to hold for exponential and geometric distributions.

This chapter is based on a research program carried out in collaboration with
V.~Sidoravicius, D.~Surgailis, and M.~E.~Vares.

\section{Introduction}

The main motivation of our approach is the discrete \emph{geometric
broken line process}.
\index{broken line process!geometric}
Informally it can
be described as the following particle system with creation and
annihilation.
Consider a sequence $\xi_{t,x}$
of independent and geometrically distributed
($P(\xi_{t,x} = k) = (1-\lambda^2) \lambda^{2k}, k=0,1,\dots $
with parameter $0< \lambda^2 < 1 $)
random variables, indexed by points
$(t,x) \in {\Z}^2 $ such that both $t$ and $x$ simultaneously
are either even, or odd (the set of all such $(t,x)$ is denoted by
$\tilde {\Z}^2$.)
Assume that at each point $(t,x) \in \tilde
{\Z}^2$, $\xi_{t,x} $ pairs of particles are born with
opposite velocities $\pm 1$.
\index{birth process}%
\index{creation|see{birth process}}%
The born particles then move with
constant velocities.
At the moment when a moving particle collides
with another moving particle having opposite velocity, they
annihilate each other.
\index{annihilation}
Because of the possibility of multiple
particles being at the same time-space point, the annihilation
rule says that at the moment of collision, first annihilated are
`older' particles.
\index{annihilation!rule}

\index{broken line}
Trajectories of particles in space-time give rise to a
process of (partially coinciding)
`discrete' random broken lines whose segments (corresponding to
individual particles) lie at the angle $\pm \pi/4 $ to the
time axis.
The broken line process is stationary in time and space
with respect to natural translations of
$\tilde {\Z}^2 $ and can be
defined on $\tilde {\Z}^2 $ by consistency, as the limit of similar
broken line processes on bounded domains.
Let $\phi_{t,x} $ be the number of broken lines which pass
through the point $(t,x) \in \tilde {\Z}^2 $.
The process $\{ \phi_{t,0} \} $ is called the \emph{intersection process},
it is stationary but not i.i.d.
\index{intersection process}
and it satisfies a law of large numbers.
This setup is described in Section~\ref{sec:discretebl}, which is based on the unpublished preprint~\cite{sidoraviciusXX}.

A generalization of this object is what we call the \emph{continuum broken line process},
introduced in Section~\ref{sec:continuumbl}.
\index{broken line process!continuum}
Within this framework, instead of having a certain number 
$\xi\in\Z_+$
of pairs of particles being born at each space-time point, we have a mass whose value is in principle
non-integer.
A few subtle differences will arise from this difference but one can overcome them by looking for broken lines of finite size, even when the process is defined on the whole lattice.
In view of Proposition~\ref{prop:llpp} and its explicit construction, it will hopefully be clear that the broken line process and its continuum version are indeed natural objects.
Also presented in Section~\ref{sec:continuumbl} is the discussion of necessary and sufficient conditions for reversibility, which turns out to be crucial for the applications that appear after.
Applications to last passage percolation are shown in Section~\ref{sec:appl}.

\section{The geometric broken line process}
\label{sec:discretebl}

In this section we shall describe the evolution of a particle system with discrete time step.
This system is first defined for finite space-time domains.
Later on it is described as a non-homogeneous Markov chain and we shall see that this evolution is reversible and consistent.
As a consequence of the above properties, the model is extended to the whole lattice.
We shall define the broken lines starting from the space-time trajectories of the particles.
It is instructive to consider the discrete geometric broken lines process first, before proceeding to full generality.
This presentation is based on~\cite{sidoraviciusXX}.

\subsection{The evolution of a particle system}
\label{sec:evolution}

\index{domain|see{hexagonal domain}}
\index{hexagonal domain}
Let
\begin{equation*}
S = \bigl\{(t,x)\in \tilde {\Z}^2: t_0 \leqslant t\leqslant t_1,
x_{t,-} \leqslant x\leqslant x_{t,+} \bigr\},
\end{equation*}
where $t_0, t_1 \in {\Z}, t_0 < t_1 $ are given
points, and $x_{t,-} \leqslant x_{t,+}, t_0 \leqslant t\leqslant t_1$
are paths in $\tilde {\Z}^2 $ such that for some
$t_0 \leqslant t_{0,1}^\pm \leqslant t_1$,
$ x_{t+1, \pm} - x_{t,\pm} = \pm 1 \,(t_0 \leqslant t < t^\pm_{0,1}) $
and $x_{t+1,\pm} - x_{t,\pm} = \mp 1 \,(t_{0,1}^\pm\leqslant t < t_1)$;
in other words, $x_{t,+} $ increases on the interval $[t_0, t^+_{0,1}]$
and decreases on $[t^+_{0,1}, t_1]$, while $x_{t,-} $ decreases on
$[t_0,t^-_{0,1}] $ and increases on $[t^-_{0,1}, t_1]$. Any such set
$S \subset \tilde {\Z}^2$ will be called  a \emph{hexagonal domain}.
\index{hexagonal domain}
\begin{figure}[!htb]
{
\small
\psfrag{x}{$x$}
\psfrag{t01-}{$t_{0,1}^-$}
\psfrag{t01+}{$t_{0,1}^+$}
\psfrag{t0}{$t_{0}$}
\psfrag{t1}{$t_{1}$}
\psfrag{t}{$t$}
\begin{center}
 \includegraphics[scale=.3]{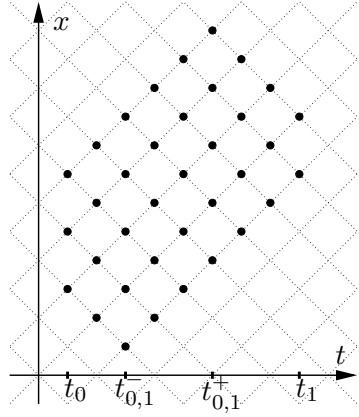}
\end{center}
\caption{
\small
a hexagonal domain}
\label{fig:hexdomain}
}
\end{figure}
Let
$\partial_0 S = \bigl\{(t,x)\in S: t=t_0, x_{t_0,-} \leqslant x\leqslant
x_{t_0,+}\bigr\} $ denote
the left vertical boundary of $S$, and
$\partial_+S =
\bigl\{(t,x) \in S: t_0\leqslant t\leqslant t^+_{0,1}, x= x_{t,+} \bigr\},
\partial_- S = \bigl\{(t,x) \in S: t_0 \leqslant t\leqslant t^-_{0,1},
x= x_{t,-}\bigr\} $ be the `northwest' and the `southwest' boundaries,
respectively. Put also $\partial_+^\circ S
= \partial_+ S \backslash \partial_0 S, \partial_-^\circ S =
\partial_-S\backslash \partial_0 S$.

For any point $(t,x) \in S$, let $\eta^+_{t,x}, \eta^-_{t,x} $
denote the numbers of ascending
(moving with velocity $+1$)  and descending
(moving with velocity $-1$) particles
which `leave' $(t,x)$, and $\zeta^+_{t,x}, \zeta^-_{t,x} $ the
respective numbers of ascending and descending particles which `come' to
$(t,x)$. The relation
\begin{equation}
\label{eq:2.2}
\zeta^+_{t,x} - \zeta^-_{t,x} =
\eta^+_{t,x} - \eta^-_{t,x} 
\end{equation}
says that at each point $(t,x)\in S$, particles can be created or killed
by pairs only. More precisely, we assume that each pair of incoming particles
with opposite velocities is annihilated, and a certain number
$\xi_{t,x}=0,1, \dots $ of \emph{pairs} of particles with opposite
velocities is born at $(t,x)$, so that
\begin{equation}
\label{eq:2.3}
\eta^+_{t,x} = \xi_{t,x} + \bigl[\zeta^+_{t,x}- \zeta^-_{t,x}\bigr]^+,
\quad
\eta^-_{t,x} = \xi_{t,x} + \bigl[\zeta^-_{t,x}- \zeta^+_{t,x}\bigr]^+.
\end{equation}
Moreover,
\begin{equation}
\label{eq:2.4}
\zeta^\pm_{t,x} = \eta^\pm_{t-1,x\mp 1}, \quad (t,x) \in
S\backslash
(\partial_\mp S \cup \partial_0 S),
\end{equation}
as all transformations
of particles may occur only at lattice points $(t,x) \in S$.
Put
\begin{eqnarray*}
& &
 \eta = (\eta^+_{t,x}, \eta^-_{t,x}: (t,x) \in S),
\\ & &
 \xi = (\xi_{t,x}: (t,x) \in S),
\\ & &
 \zeta_0 = (\zeta^+_{t,x},\zeta^-_{t,x}: (t,x) \in \partial_0 S),
\\ & &
 \zeta_+ = (\zeta^-_{t,x}: (t,x) \in \partial_+^\circ S),
\\ & &
 \zeta_- = (\zeta^+_{t,x}: (t,x) \in \partial_-^\circ S).
\end{eqnarray*}
\label{page:zetacirc}%
$\zeta^\circ = (\zeta_0, \zeta_+, \zeta_-)$,
and regard $\eta $ as configuration in $S$ and
$\zeta^\circ $ as boundary
data.
It is clear that there is a 1-1 correspondence between
$(\zeta^\circ, \xi ) $ and $(\zeta^\circ, \eta)$, where
$\zeta^\pm_{t,x}, \eta^\pm_{t,x}, \xi_{t,x} $ are related
by~(\ref{eq:2.2})-(\ref{eq:2.4}).
Note that~(\ref{eq:2.2})~and~(\ref{eq:2.4}) imply
\begin{equation}
\label{eq:2.5}
\eta^+_{t,x} - \eta^-_{t,x} = \eta^+_{t-1,x-1}
- \eta^-_{t,x+1}, \quad
(t,x) \in S\backslash (\partial_0S \cup \partial_+S \cup
\partial_- S).
\end{equation}
The probability measure $P_S(\eta)$ corresponding to the evolution of this particle
systems in the hexagonal domain $S$ can now be defined as follows.
Let $0<\lambda < 1$ be a parameter. Assume
that all $\zeta^\pm_{t,x}$'s in $\zeta^\circ $ are independent
and Geom($\lambda$)-distributed, i.e.
$$
P(\zeta^\pm_{t,x} =k) = (1-\lambda) \lambda^k, \quad k=0,1, \dots,
\quad (t,x) \in \partial_0 S \cup \partial_\mp S.
$$
Assume also that all $\xi_{t,x}$'s in $\xi $ are independent and
Geom($\lambda^2$)-distributed and, moreover, $\zeta^\circ $ and
$\xi $ are independent. Then $P_S $ is the distribution of $\eta =
(\eta^+_{t,x}, \eta^-_{t,x}: $ $ (t,x) \in S)$.
In other words,
$P_S$ is the distribution of the number of outgoing particles from
$(t,x) \in S$, assuming that $\zeta^\pm_{t,x}, (t,x) \in
\partial_0S \cup \partial_\mp S$ (the number of particles which
`immigrate' to $S$ through its left boundary $\partial_0S \cup
\partial_+S \cup \partial_-S $) are i.i.d.
Geom($\lambda$)-distributed, and $\xi_{t,x}, (t,x)\in S$ (the
number of pairs of particles with opposite velocities born inside
$S$) are i.i.d. Geom($\lambda^2$)-distributed, independent of the
$\zeta^\pm_{t,x}$'s.
In the evolution, the born particles move
with constant velocities $+1$ or $-1$ until they collide at some
lattice point $(t,x) \in S$, after which the colliding particles
die (annihilate).

\subsection
{An equivalent description for the evolution of the particles}

Below we provide another
description of the evolution of particle systems introduced
above.
In this setting various properties are more easily seen to hold true.
Let $\Sigma_S $ be the set of all configurations
$\eta $, where $\eta^\pm_{t,x}$ take values
$0,1,\dots $ and satisfy relations~(\ref{eq:2.5}).
Let $\Sigma_{\partial^\circ S} $ be the set of all
configurations $\zeta^\circ = (\zeta_0, \zeta^+, \zeta_-), \,
\zeta_0 = (\zeta^+_{t,x}, \zeta^-_{t,x}: (t,x) \in \partial_0S),
\zeta_+ = (\zeta^-_{t,x}: (t,x) \in \partial^\circ_+ S),
\zeta_- = (\zeta^+_{t,x}: (t,x) \in \partial^\circ_- S) $
on $\partial^\circ S = \partial_0S \cup \partial_+S \cup
\partial_- S$. Let
$\bar \Sigma_S $ be the set of all configurations
\index{Sigmas@$\Sigma_S$}
\index{eta@$\bar\eta$}
\label{Sigma}
$\bar \eta = (\zeta^\circ, \eta), $ where
$\zeta^\circ \in {\Sigma}_{\partial^\circ S},
\eta \in {\Sigma}_S, $ and, moreover,
\begin{equation}
\label{eq:2.6}
\eta^+_{t,x} - \eta^-_{t,x}
= \begin{cases}
 \zeta^+_{t,x}- \eta^-_{t-1,x+1}, & \mbox{ if } (t,x)\in \partial_-S \backslash \partial_0 S \\
 \eta^+_{t-1,x-1}- \zeta^-_{t,x}, & \mbox{ if } (t,x)\in \partial_+S \backslash \partial_0 S \\
 \zeta^+_{t,x}- \zeta^-_{t,x},    & \mbox{ if } (t,x)\in \partial_0 S
 \end{cases}.
\end{equation}
Finally, for a given configuration $\zeta^\circ \in
{\Sigma}_{\partial^\circ S}$,
let
${\Sigma}_{S|\zeta^\circ} $ denote the class of all 
$\eta \in {\Sigma}_S $ which satisfy~(\ref{eq:2.6}).

We shall define a probability measure  $Q_{S|\zeta^\circ}(\eta)$
on ${\Sigma}_{S|\zeta^\circ} $
as a (non homogeneous)
Markov chain $\eta_t, t_0 \leqslant t \leqslant t_1 $, whose
values at each time $t$ are restrictions of
a configuration $\eta = \bigl(\eta^+_{u,x}, \eta^-_{u,x}:
(u,x) \in S\bigr) \in {\Sigma}_{S|\zeta^\circ} $ on $t=u$;
in other words, $\eta_t = (\eta^+_{t,x}, \eta^-_{t,x}: x\in S_t)$, where
$S_t = \bigl\{x: (t,x) \in S\bigr\}$.

The transition probabilities of the Markov chain are defined as
follows (for clarity, we assume $t^+_{01} = t^-_{01} = t_{01}$
below).

\noindent (i) At $t=t_0$, the distribution of
$\eta_{t_0} = (\eta^+_{t_0,x}, \eta^-_{t,x}: x\in \partial_0 S)$
depends only on $\zeta_0 = (\zeta^+_{t,x}, \zeta^-_{t,x}:
(t,x) \in \partial_0 S)$ and is given by
\begin{equation}
\label{eq:2.7}
Q_{S|\zeta^\circ}(\eta_{t_0}|\zeta_0) = \prod_{x\in \partial_0 S}
q(\eta^+_{t_0,x}, \eta^-_{t,x}| \zeta^+_{t,x}, \zeta^-_{t,x}),
\end{equation}
where
\begin{equation}
\label{eq:2.8}
q(n^+, n^-| m^+, m^-) = Z^{-1}_{m^+,m^-} \lambda^{n^++ n^-}
\delta_{\{n^+ -n^- = m^+-m^-\}},
\end{equation}
$n^\pm, m^\pm =0,1,\dots,\,
0< \lambda <1 $ is a parameter, and
\begin{equation}
\label{eq:2.9}
Z_{m^+,m^-} = \sum_{n^+,n^-} \lambda^{n^++n^-}
\delta_{\{n^+-n^- = m^+-m^-\}}
= \frac{\lambda^{|m^+-m^-|}}{(1-\lambda^2)}.
\end{equation}

\noindent (ii) Let $t_0 < t \leqslant t_{01}$. Then
$\eta_t = (\eta^+_{t,x}, \eta^-_{t,x}, x\in S_t) $
depends only on $\eta_{t-1} =
(\eta^+_{t-1,x}, \eta^-_{t-1,x}, x\in
S_{t-1}) $ and $ \zeta_{t,+} = \zeta^-_{t, x_{t,+}},
\zeta_{t,-} = \zeta^+_{t,x_{t,-}} $, according to the transition
probability
\begin{equation}
\label{eq:2.10}
\begin{array}{rcl}
Q_{S|\zeta^\circ}(\eta_t|\eta_{t-1}, \zeta_{t,\pm})
&=& q(\eta^+_{t,x_{t,+}}, \eta^-_{t,x_{t,+}}|\eta^+_{t-1,x_{t-1,+}},
\zeta^-_{t,x_{t,+}})\\
&& \times\ \ q(\eta^+_{t,x_{t,-}}, \eta^-_{t,x_{t,-}}|\eta^+_{t-1,x_{t-1,-}},
\zeta^+_{t,x_{t,-}})\\
&& \times\ \ \prod_{x\in S_t, x\neq x_{t,\pm}} q(\eta^+_{t,x}, \eta^-_{t,x}|
\eta^+_{t-1,x-1}, \eta^-_{t-1,x+1})
\end{array}.
\end{equation}

\noindent (iii) Let $t_{01} < t \leqslant t_1$. Then
$\eta_t = (\eta^+_{t,x}, \eta^-_{t,x}, x\in S_t) $
depends only on $\eta_{t-1} =
(\eta^+_{t-1,x}, \eta^-_{t-1,x}, x\in
S_{t-1}) $ and
\begin{equation}
\label{eq:2.11}
Q_{S|\zeta^\circ}(\eta_t|\eta_{t-1}) =
\prod_{x\in S_t} q(\eta^+_{t,x}, \eta^-_{t,x}|
\eta^+_{t-1,x-1}, \eta^-_{t-1,x+1}).
\end{equation}

Let $\Pi_\lambda (\zeta^\circ) $ be the product geometric distribution, i.e. all random variables
$\zeta^+_{t,x}, (t,x)\in \partial_0 S \cup \partial_-S,
\zeta^+_{t,x}, (t,x) \in \partial_0 S \cup
\partial_+ S  $
are independent and Geom($\lambda $)-distributed.
Define the probability measures
$Q_S (\bar \eta) = Q_S (\zeta^\circ, \eta), Q_S (\eta) $ on
$\bar {\Sigma}_S, {\Sigma}_S $  by
$$
Q_S(\zeta^\circ, \eta ) =  \Pi_\lambda (\zeta^\circ)
Q_{S|\zeta^\circ}(\eta),  \qquad Q_S (\eta) = \sum_{\zeta^\circ \in
{\Sigma}_{\partial^\circ S}} Q_S(\zeta^\circ, \eta),
$$
respectively. We claim that
\begin{equation}
\label{eq:2.13}
Q_S = P_S,
\end{equation}
where $P_S $ was defined in Section~\ref{sec:evolution} above. Because of the
$1-1$ correspondence $(\zeta^\circ, \eta) \leftrightarrow
(\zeta^\circ, \xi)$, (2.13) follows from
\begin{equation}
\label{eq:2.14}
Q_{S|\zeta^\circ}(\xi) :=
Q_{S|\zeta^\circ}\bigl(\eta: \xi (\zeta^\circ, \eta) = \xi\bigr)
= \Pi_{\lambda^2}(\xi).
\end{equation}
Equation~(\ref{eq:2.14}) can be shown by induction in $t_1-t_0=0,1,\dots $.
Let $t_1=t_0$.
As $2\xi_{t,x} = \eta^+_{t,x} +
\eta^-_{t,x} - |\zeta^+_{t,x}-\zeta^-_{t,x}|,
(t,x) \in \partial_0 S$, see~(\ref{eq:2.3}), so by~(\ref{eq:2.8}),(\ref{eq:2.9})
\begin{eqnarray*}
&&
P(\xi_{t,x}= k|\zeta^+_{t,x} = m^+, \zeta^-_{t,x} = m_-)
=
\\
&&
= \sum_{n^++n^- = 2k-|m^+-m^-|} q(n^+,n^-|m^+,m^-)\\
&&
= Z^{-1}_{m^+,m^-}
\sum_{n^++n^- = 2k-|m^+-m^-|}
\lambda^{n^++n^-} \delta_{n^+-n^-=m^+-m^-} \\
&&
=
(1-\lambda^2) \lambda^{2k} = \pi_{\lambda^2} (k),
\end{eqnarray*}
where $\pi_\lambda $ stands for Geom($\lambda$)-distribution.
Hence and by~(\ref{eq:2.7}) it follows that~(\ref{eq:2.14}) holds for $t_0 =t_1$.
In the general case,~(\ref{eq:2.14}) and~(\ref{eq:2.13}) follow by induction
in $t_1-t_0$ and a similar computation, based on the
explicit form of transition probabilities
$Q_{S|\zeta^\circ}(\eta_t|\eta_{t-1})$ given by~(\ref{eq:2.10})-(\ref{eq:2.11}).

In the sequel we discuss some properties of $P_S (\bar \eta)$.

\subsection{Duality}
The transition probability~(\ref{eq:2.8})
satisfies the following relation
\begin{equation}
\label{eq:2.15}
\pi_\lambda (m^+) \pi_\lambda (m^-) q(n^+,n^-|m^+,m^-) =
\pi_\lambda (n^+) \pi_\lambda (n^-)
q(m^-,m^+|n^-,n^+),
\end{equation}
$m^\pm, n^\pm = 0,1,\dots, $ where $\pi_\lambda $
is Geom($\lambda$)-distribution.  A similar
duality relation holds for transition probabilities of
the Markov chain $\eta_t, t_0 \leqslant t\leqslant t_1$. For example,
if $t_{01} < t\leqslant t_1$, then from~(\ref{eq:2.15}),(\ref{eq:2.11}) it follows that
\begin{eqnarray*}
&\prod_{x\in S_t} \pi_\lambda (\eta^+_{t-1,x+1}) \pi_\lambda(\eta^-_{t-1,x+1})
q(\eta^+_{t,x},\eta^-_{t,x}| \eta^+_{t-1,x-1},
\eta^-_{t-1,x+1}) \\
&= \prod_{x\in S_t} \pi_\lambda (\eta^+_{t,x}) \pi_\lambda (\eta^-_{t,x})
q(\eta^-_{t-1,x+1},\eta^+_{t-1,x+1}| \eta^-_{t,x},
\eta^+_{t,x}).
\end{eqnarray*}
The above
duality implies that the construction of
$P_S (\bar \eta) = Q_S(\bar \eta) $ can be reversed in time. Namely, in the
`dual picture', $\eta^+_{t,x} $ is the number of
`descending' particles which `come' to $(t,x)$ from the right,
and $\zeta^-_{t,x}$ is the number of `ascending' particles which
`leave' $(t,x)$ in the same direction. A similar `dual interpretation'
can be given to $\eta^-_{t,x}, \zeta^+_{t,x}$. In the dual
construction, $\eta^\circ = \{(\eta^+_{t,x}, \eta^-_{t,x},
(t,x) \in \partial_1 S), (\eta^+_{t,x}, (t,x) \in \partial^+ S
\backslash \partial_1S), (\eta^-_{t,x}, (t,x) \in
\partial^- S \backslash \partial_1 S)\} $ is the boundary
condition; here, $\partial_1 S$ is the right
vertical boundary, and $\partial^+ S, \partial^- S $ are the
`northeast' and `southeast' boundaries of $S$,
respectively.
The probability measure $\hat Q_{S|\eta^\circ} (\zeta) $ is defined
on the set $\hat {\Sigma}_{S|\eta^\circ} $ of all configurations
$\zeta = (\zeta^+_{t,x}, \zeta^-_{t,x}, (t,x) \in
S) $ which satisfy conditions analogous to~(\ref{eq:2.5}),(\ref{eq:2.6}); namely,
$$
\zeta^-_{t,x} - \zeta^+_{t,x} =
\zeta^-_{t+1,x-1} - \zeta^+_{t+1,x+1}, \quad
(t,x) \not\in \partial_1 S \cup \partial^+ S \cup \partial^- S,
$$
and
$$
\zeta^-_{t,x} - \zeta^+_{t,x} = 
\begin{cases}
\zeta^-_{t+1,x-1} - \eta^+_{t,x}, & \mbox{ if } (t,x) \in \partial^+S \backslash \partial_1 S \\
\eta^-_{t,x} - \zeta^+_{t+1,x+1}, & \mbox{ if } (t,x) \in \partial^- S\backslash \partial_1 S \\
\eta^-_{t,x} - \eta^+_{t,x},      & \mbox{ if } (t,x) \in \partial_1 S
\end{cases}.
$$
The definition of $\hat Q_{S|\eta^\circ} (\zeta) $ is
completely analogous to that of $Q_{S|\zeta^\circ}(\eta) $ and
uses a Markov chain $\hat \eta_t, t_0 \leqslant t\leqslant t_1$ run
in the time reversed direction, whose transition
probabilities are analogous to~(\ref{eq:2.7}),(\ref{eq:2.10}),(\ref{eq:2.11}). Then
if $\Pi_\lambda (\eta^\circ) $ is the product geometric distribution on
configurations $\zeta^\circ$, we set
$$
\hat Q_S (\zeta, \eta^\circ) = \Pi_\lambda (\eta^\circ)
\hat Q_{S|\eta^\circ}(\zeta)
$$
and obtain the equality  for the probability measures on $\bar {\Sigma}_S$
\begin{equation}
\label{eq:2.16}
\hat P_S = \hat Q_S
= Q_S = P_S.
\end{equation}

\subsection{Consistency}

\index{restricted process}
Let $S' \subset S''$ be two bounded hexagonal domains in
$\tilde {\Z}^2 $, and let $P_{S'}, P_{S''} $
be the probability distributions on the
configuration spaces $\bar {\Sigma}_{S'}, \bar {\Sigma}_{S''}$,
respectively, as defined above. The probability measure
$P_{S''}$ on $\bar {\Sigma}_{S''}$ induces a probability
measure $P_{S''|S'}$ on $\bar {\Sigma}_{S'}$ which is the
$P_{S''}-$
distribution of the
\emph{restricted process}
\index{restricted process}
$\bar \eta_{t,x}, (t,x) \in S$. Then the following
\emph{consistency property} is true:
\index{consistency property}
\begin{equation}
\label{eq:2.17}
P_{S''|S'} = P_{S'}.
\end{equation}
The proof of~(\ref{eq:2.17}) uses~(\ref{eq:2.16}) and the argument in
Arak and Surgailis~\cite[Theorem~4.1]{arak89}.
Let $\chi_\uparrow, \chi_\nearrow, \chi_\searrow $ denote a vertical line, an ascending
line, a descending line, respectively, in $\tilde {\Z}^2 $, with
the last two having slopes $\pi/4, -\pi/4$, respectively. Any such
line $\chi $ partitions $\tilde {\Z}^2 $ into the left part
$\tilde {\Z}^2_{\chi,-} =
\bigl\{(t,x)\in \tilde {\Z}^2: t\leqslant t', x\leqslant x' \mbox{ for some }
(t',x') \in \chi \bigr\} $ (which contains the line
$\chi $ itself)
and the right part $\tilde {\Z}^2_{\chi,+} =
\tilde {\Z}^2 \backslash \tilde {\Z}^2_{\chi,-}$.
Note it
suffices to show~(\ref{eq:2.17}) for
$$
S' = S'' \cap \tilde {\Z}^2_{\chi,\pm},
$$
for any line $\chi = \chi_\uparrow, \chi_\nearrow, \chi_\searrow $
of the above type.

In the simplest case $S'= S''\cap \tilde {\Z}^2_{\chi_\uparrow,-} $ (\ref{eq:2.17}) trivially follows from the
construction: in this case, $P_{S''|S'} $ is nothing else but the
evolution $\eta_t, t_0\leqslant t\leqslant t'$ observed up to the moment
$t'\leqslant t_1, (x,t') \in \chi_\uparrow $, and therefore coincides
with  $P_{S'}$. The case $S' = S''\cap \tilde {\Z}^2_{\chi_\nearrow,-} $ follows by the observation that $\bar
\eta_{t,x}, (t,x) \in \tilde S\cap \tilde {\Z}^2_{\chi_\nearrow,+} $ do not participate in the definition of
the probability of $\bar \eta_{t,x}, (t,x) \in S\cap \tilde {\Z}^2_{\chi_\nearrow, -}$: the evolution of the particles after
they exit through $\chi_\nearrow $ has no effect on the evolution
before they exit this line. The case $S' = S''\cap \tilde {\Z}^2_{\chi_\searrow, -} $ is completely analogous.

Non-trivial cases of~(\ref{eq:2.17}) are $S' = S''\cap
\tilde {\Z}^2_{\chi,+},\,\chi = \chi_\uparrow, \chi_\nearrow,
\chi_\searrow $. However, the `reversibility'~(\ref{eq:2.16}) allows to exchange the right and left directions
by replacing $P_S $ by the `reversed' process $\hat P_S$ and thus
reducing the problem to the above considered cases
of $S' = S'' \cap \tilde {\Z}^2_{\chi,-} $.

By consistency~(\ref{eq:2.17}), the evolution of particles defined in finite
hexagonal domains, can be extended to the evolution on the
whole lattice $\tilde {\Z}^2 $. Let
${\Sigma}= {\Sigma}_{\tilde {\Z}^2} $ be the set of all
configurations $\eta = (\eta^+_{t,x}, \eta^-_{t,x},
(t,x) \in \tilde {\Z}^2 \}$ satisfying~(\ref{eq:2.5}) for each
$(t,x) \in \tilde {\Z}^2 $. Then there exists a probability
measure $P= P_{\tilde {\Z}^2} $ on
${\Sigma}$ whose restriction $P_{|S}$ to an arbitrary hexagonal
domain $S $ coincides with $P_S$:
$$
P_{|S} = P_S.
$$
Furthermore, $P$ is invariant with respect to translations of
$\tilde {\Z}^2 $.

\subsection{The discrete broken line process}

In this section we shall
describe the construction of the broken line process
in a finite hexagonal domain $S \subset \tilde {\Z}^2 $.
This construction is more similar to the construction in the continuous Poisson model, and will yield a more streamlined proof of the the law of large numbers of the intersection process defined below.

Given a hexagonal domain
$S\subset \tilde {\Z}^2 $,
we denote by \label{page:sbar}%
$\bar S = \bigl\{(t,x) \in \tilde {\Z}^2: |(t,x) -(t',x')| \leqslant \sqrt{2}
\mbox{ for some }(t',x') \in S \bigr\}$.
In other words, $\bar S$
consists of all points of $S$ plus the neighboring points which lie
at the distance $\sqrt{2} $ to the boundary of $S$.
A
\emph{broken line}
\index{broken line}
in $\bar S$
is a finite path
$\gamma = \{ (t_i,i): (t_i,i)  \in \bar S \}$ such that
$|t_{i+1}- t_i| = 1$ for all  $i \in {\Z}$
such that $(t_i,i) \in \bar S, (t_{i+1},i+1) \in \bar S$.
A broken line in $\bar S$ can be
identified with the graph of a continuous function
$x\mapsto \gamma (x): [i_-(\gamma), i_+(\gamma)]
\to {\R}$ which is
linear on each interval $[i,i+1]$ and satisfies
$\gamma (i) = t_i$ for each  $i\in {\Z}$;
here $i_-(\gamma) < i_+(\gamma) $ are integers
which satisfy $(t_{i_\pm(\gamma)}, i_\pm (\gamma)) \in \bar S\backslash S$.
Clearly, the graph of the function
$ \gamma (\cdot)$ belongs to the union of all
edges of $\bar S$ which are identified with
closed intervals of the length $\sqrt{2}$ between nearest
neighbors of $\bar S $. For each $(t,x) \in
\bar S$, we shall
denote by
$e^\pm_{t,x}$ the edge (interval)
between the points $(t,x)$ and $(t+1,x\pm1)$.

Given two broken lines $\gamma', \gamma''$, we write
$\gamma' \preceq \gamma'' $ if $\gamma'(x) \leqslant \gamma''(x)
\,\forall x\in [i_-(\gamma'),$ $i_+(\gamma')]\cap
[i_-(\gamma''),i_+(\gamma'')]$, and $\gamma' \prec \gamma'' $ if
$\gamma' \preceq \gamma'', \gamma' \neq \gamma''$ holds.
For any $(t,x) \in S$, we denote by
$e^\nearrow_{t,x}, e^\nwarrow_{t,x}, e^\searrow_{t,x},
e^\swarrow_{t,x} $ the edge of
$\tilde {\Z}^2 $ incident with $(t,x)$ and
lying in northeast, northwest, southeast, southwest direction,
respectively, so that $e^\nearrow_{t,x} = e^\swarrow_{t+1,x+1},
e^\nwarrow_{t,x} = e^\searrow_{t-1,x+1} $. 
Clearly, each edge incident to a point of the form $(t,x) \in S$, connects this point with some point of $\bar S$.

With each configuration
$\bar \eta = (\zeta^\circ, \eta) \in \bar {\Sigma}_S$
(see Section~\ref{Sigma}), one can associate
a finite partially ordered system
$\{\gamma_j, j=1, \dots, M \} $ of broken lines in $\bar S$,
such that for any point
$(t,x) \in \tilde {\Z}^2 $, the relations
\begin{equation}
\label{eq:2.19}
\nu (e^\nearrow_{t,x}) =
\eta^+_{t,x},
\quad \nu (e^\searrow_{t,x}) = \eta^-_{t,x},
\quad \nu (e^\nwarrow_{t,x}) = \zeta^-_{t,x},
\quad \nu (e^\swarrow_{t,x}) = \zeta^+_{t,x},
\end{equation}
hold, where
$\nu (e) $ is the number of broken lines
$\gamma_i $ which pass through a given edge $e$
of $\tilde {\Z}^2 $, and where
$\eta^+_{t,x}, \eta^-_{t,x}, \zeta^+_{t,x}, \zeta^-_{t,x} $
are related by~(\ref{eq:2.2}),(\ref{eq:2.4}), and
denote the respective numbers of outgoing and incoming
particles to a given site $(t,x) \in S$, see Section~\ref{sec:evolution}.
Because of these relations and~(\ref{eq:2.19}), it is convenient to denote
configurations $\bar \eta \in \bar {\Sigma}_S $ as
\begin{equation}
\label{eq:2.20}
\bar \eta = \bigl(\eta (e), e\in \E (\bar S)\bigr),
\end{equation}
where $\E (\bar S)$ is the set of all edges of $\bar S$, and
\index{ES@$\E(\bar S)$}
where $\eta (e) $ is the number of particles which pass through
along the edge $e$; $\eta (e) = \eta^+_{t,x}, \eta^-_{t,x},
\zeta^+_{t,x}, \zeta^-_{t,x} $ for
$e = e^\nearrow_{t,x}, e^\swarrow_{t,x}, e^\swarrow_{t,x},
e^\nwarrow_{t,x}$, respectively.
\index{eta@$\bar\eta$}
\index{flow field}
This field $\bar\eta$ is what we shall call \emph{flow field}.

According to this notation (\ref{eq:2.19}) becomes
$$
\nu (e) = \eta (e), \quad e\in \E(\bar S). 
$$

The broken lines for a given configuration $\bar \eta \in
\bar {\Sigma}_S $ can be constructed as follows.
For each particle passing through a given edge
$e\in E(\bar S)$, we define a label $p\in\{1,\dots, \eta (e)\}$, which is
interpreted as the
\emph{relative age}
\index{relative age}
of that particle among
all particles which pass through the same edge.
The particle whose relative age is $p=1$ is the oldest and that with $p=\eta(e)$ is the youngest.
This way, we put a coherent order on moving particles; more formally, any particle is characterized by a pair
$$
(e, p),  \quad e\in E(\bar S), \,\,p=1,2, \dots, \eta (e). 
$$

\index{association rules!discrete case}
\index{generation}
We shall define now a relation between labeled particles
$(e_1,p_1), (e_2, p_2)$ on
two adjacent edges $e_1, e_2 \in E(\bar S)$, which we denote by
\begin{equation*}
(e_1,p_1) \sim (e_2,p_2).
\end{equation*}
If relation $(e_1,p_1)\sim(e_2,p_2)$ holds, we say that
they are \emph{associated}, or \emph{belong to the same generation}.
The adjacent edges $e_1, e_2$
by definition are any two edges $e^\nearrow_{t,x}, e^\nwarrow_{t,x}, e^\searrow_{t,x},
e^\swarrow_{t,x} $ incident with some $(t,x) \in S$.
The relation $(e_1,p_1) \sim (e_2,p_2)$ is defined in the following cases:
\index{association rules!discrete case}

Case 1: 
$e_1 = e_{t,x}^\swarrow$, $e_2 = e_{t,x}^\nwarrow$;
$p_1,p_2 \leqslant \eta(e_1)\wedge\eta(e_2)$

Case 2:
$e_1 = e_{t,x}^\searrow$, $e_2 = e_{t,x}^\nwarrow$;
$\eta(e_2) > \eta (e^\swarrow_{t,x})$, $p_2 > \eta (e^\swarrow_{t,x})$.

Case 3:
$e_1 = e_{t,x}^\swarrow$, $e_2 = e_{t,x}^\nearrow$;
$\eta(e_1) > \eta (e^\nwarrow_{t,x})$, $p_1 > \eta (e^\nwarrow_{t,x})$.

Case 4:
$e_1 = e_{t,x}^\searrow$, $e_2 = e_{t,x}^\nearrow$;
$\eta(e_i)-\xi_{t,x} < p_i \leqslant \eta(e_i)$.

Namely, $(e_1,p_1) \sim (e_2,p_2)$ holds if and only if

Case 1:
$p_2 = p_1$;

Case 2:
$p_1 = p_2 - \eta(e_{t,x}^\swarrow)$;

Case 3: 
$ p_2 = p_1 - \eta(e_{t,x}^\nwarrow)$;

Case 4:
$ \eta(e_1) - p_1 = \eta(e_2) - p_2$.

It is clear that associated particles are either the particles
which annihilate each other at $(t,x)$ (Case 1), or the
particles born at $(t,x)$ and moving in opposite directions (Case 4), or, as
in Cases 2 and 3, we associate a
younger incoming particle which is not killed at $(t,x)$ with the
corresponding older outgoing particle (both moving in the
same direction).

Now, we say that a finite path
$\gamma = \{(t_i,i): (t_i,i)\in \bar S \} $ is a broken line
in a given configuration $\bar \eta $ (\ref{eq:2.20})
if any two adjacent edges of this path are associated, i.e.,
belong to the same generation. Let $\{\gamma_j, j=1, \dots, M \} $
be the family of all broken lines in a given configuration
$\bar \eta $.
The family $\{\gamma_j,j=1,\dots, M\} $
is partially order by the relation $\preceq $ defined above
(with the possible exception of broken lines that exit
$S$ at its left or right boundaries, in which case two lines
may not be ordered).
This follows from the definition of a broken line and the fact that different pairs of associated particles cannot `cross' each other: if $ (e'_1, p'_1) \sim (e'_2, p'_2), \,(e''_1, p''_1) \sim (e''_2, p''_2) $ and $e'_1,e'_2, e''_1, e''_2 $ are all incident with the same vertex $(t,x) $, then the paths $(e'_1,e'_2) $ and $(e''_1, e''_2)$ cannot cross each other, in other words, they cannot lie on different lines intersecting at $(t,x)$.

\subsection{The intersection process}

Let $\chi $ be an arbitrary
line in $\tilde {\Z}^2 $ of the type
$\chi = \chi_\nearrow$,
$\chi_\searrow$, $\chi_\uparrow$, $\chi_\rightarrow$, where $\chi_\rightarrow $ stands for a
horizontal line: $\chi_\rightarrow =
\{(t,x) \in \tilde {\Z}^2: x= x'\} $ for some
$x'\in \Z $. For any such line $\chi $ and
a finite hexagonal domain $S \subset \tilde {\Z}^2 $
we shall define the
\emph{intersection process}
\index{intersection process}
as the number of broken lines
$\gamma_i$'s which intersect $\chi $ at some point $(t,x) \in \chi
\cap S$. However because of the peculiarity of the discrete situation, the definition
of the intersection process requires some care. We shall consider
separately the four above types of $\chi $.

1. Case $\chi = \chi_\nearrow $. Let
$$
\tau_{t,x}(\chi) := \nu (e^\nwarrow_{t,x}),  \quad (t,x) \in
\chi_\nearrow \cap S
$$
be the number of broken lines which intersect $(t,x) \in \chi_\nearrow$ from
northwest.

2. Case $\chi = \chi_\searrow $. Let
$$
\tau_{t,x}(\chi) := \nu (e^\swarrow_{t,x}),  \quad (t,x) \in
\chi_\searrow \cap S
$$
be the number of broken lines which intersect $(t,x) \in \chi_\searrow$ from
southwest; we recall that
$e^\swarrow_{t,x} = e^\nearrow_{t-1,x-1}$ is the
edge between $(t,x)$ and $(t-1,x-1)$.

3. Case $\chi = \chi_\uparrow $. Let
$$
\tau^-_{t,x}(\chi) := \nu (e^\nwarrow_{t,x}), \quad
\tau^+_{t,x}(\chi) := \nu (e^\swarrow_{t,x}),  \quad (t,x) \in
\chi_\uparrow \cap S
$$
be the number of broken lines which intersect $\chi_\uparrow$
at $(t,x) $ from northwest and southwest, respectively.

4. Case $\chi = \chi_\to $. Let
$$
\tau^+_{t,x}(\chi) := \nu (e^\nearrow_{t,x}), \quad
\tau^-_{t,x}(\chi) := \nu (e^\nwarrow_{t,x}),  \quad (t,x) \in
\chi_\to \cap S\leqno
$$
be the number of broken lines which intersect $\chi_\to$
at $(t,x) $ from northeast and northwest, respectively.

The intersection processes are correspondingly defined as
\begin{eqnarray*}
&&
\tau_S (\chi_\nearrow) =
\bigl(\tau_{t,x}(\chi_\nearrow)\bigr)_{(t,x) \in \chi_\nearrow \cap S}
\\
&&
\tau_S( \chi_\uparrow) =
\bigl(\tau^+_{t,x}(\chi_\uparrow), \tau^-_{t,x}(\chi_\uparrow)\bigr)_
{(t,x) \in \chi_\uparrow\cap S}
\\
&&
\tau_S (\chi_\searrow) = 
\bigl(\tau_{t,x} (\chi_\searrow) \bigr)_{ (t,x) \in \chi_\searrow \cap S}
\\
&&
\tau_S (\chi_\to ) =
\bigl(\tau^+_{t,x}(\chi_\to), \tau^-_{t,x}(\chi_\to)\bigr)_
{(t,x) \in \chi_\to \cap S}
\end{eqnarray*}

The properties below follow from definition~(\ref{eq:2.19}) %
and the consistency~(\ref{eq:2.17}).

The intersection processes
$ \tau_S (\chi_\nearrow), \tau_S (\chi_\searrow), \tau_S (\chi_\uparrow) $
are i.i.d.~sequences of Geom$(\lambda)$-distributed random
variables. The process $ \tau_S (\chi_\to )$ is stationary, in
the sense that there exists a stationary process
$$
\tau (\chi_\to) = \bigl(\tau^+_{t,x}(\chi_\to), \tau^-_{t,x}(\chi_\to)\bigr)_
{ (t,x) \in \chi_\to },
$$
such that, for each finite
hexagonal domain $S\subset \tilde {\Z}^2 $,
the distribution of the process $\tau (\chi_\to)$ restricted to
$(t,x) \in S$ coincides with the distribution of
$\tau_S (\chi_\to) $. 

The marginal probability laws of all four processes
coincide, i.e.
$$
\tau_{t,x}(\chi_\nearrow)
\stackrel d=
\tau_{t,x}(\chi_\searrow)
\stackrel d=
\tau^\pm_{t,x} (\chi_\uparrow)
\stackrel d=
\tau^\pm_{t,x}(\chi_\to)
\stackrel d=
{\rm Geom}(\lambda). 
$$

Note that the process $\tau (\chi_\to)$ is not i.i.d.
This fact follows by an explicit
computation of the bivariate distribution
$$
P(\tau^-_{t,x}(\chi_\to) = m, \tau^+_{t,x}(\chi_\to) = n)
= (1-\lambda)(1-\lambda^2) \lambda^{2(m+n)}
(\lambda^{-m} + \lambda^{-n} - 1-\lambda).
$$
However,
$$
E \tau^+_{t,x}(\chi_\to) = E \tau^-_{t,x}(\chi_\to) = \lambda/(1-\lambda).
$$

Let
$$
T_n := \sum_{t=1}^n (\tau^+_{t,0}(\chi_\to) +
\tau^-_{t,0}(\chi_\to))
$$
be the number of broken lines which intersect interval
$(0,n)$ of the time axis (the sum above is taken over all
points $(t,0) \in \tilde {\Z}^2, 0\le t\le n$).
Then a.s.
$$
\lim_{n\to \infty} n^{-1}T_n = \lambda/(1-\lambda),
$$
that is, the intersection process satisfies a law of large numbers.

\sectionmark{The continuum broken line process}

\section
{The continuum broken line process on a discrete domain}
\label{sec:continuumbl}

\sectionmark{The continuum broken line process}

In this section we present a fairly natural generalization of the broken line process, namely the continuum broken line process. Within this framework, instead of having a certain number $\xi\in\Z_+$ of pairs of particles being born at each site, we have a mass $\xi\in\R_+$. In this case a broken line is described not only by the sites it occupies, this object has some thickness as well.
The continuum broken lines are not determined by sequentially associating pairs $(e,p)$, we rather need to associate mass. Accordingly, neither can they be understood by just following the mass associations over and over, as this mass may well branch due to the association rules. 
At first we consider only broken lines of finite length. But this restriction is in fact empty: all the important features of the process are well captured by describing finite broken lines.

Applications to last passage percolation are shown in Section~\ref{sec:appl}. In view of Proposition~\ref{prop:llpp} and its explicit construction, it will hopefully be clear that the broken line process and its continuum version are indeed natural objects.
The proofs of Theorem~\ref{theo:exp} and Theorem~\ref{theo:llngeo} illustrate how this process can be useful.

We have chosen to frequently abuse the notation to avoid obscuring the argument.
Such abuse should not give rise to any confusion.

\subsection{Evolution of a continuum-particle system on the discrete lattice, flow fields}

\index{rectangular domain}
\index{domain|see{rectangular domain}}
Let $S$ be a fixed hexagonal domain in $\tilde\Z^2$.
Consider $\bar S$ defined on p.~\pageref{page:sbar}, and take $\partial\bar S=\bar S\backslash S$,
$\E(\bar S)=\bigl\{e=\langle y,y'\rangle\in\E(\tilde\Z^2):y\in S\bigr\}$.
Usually we will denote the elements of $\tilde\Z^2$ by $y=(t,x)$ and the elements of $\E(\tilde\Z^2)$ by $e$.
We call $S$ a \emph{rectangular domain} if it is a degenerate hexagonal domain, i.e., $\#\partial_0S=\#\partial_1S=1$.
It is convenient to define, for a rectangular domain $S$, the boundaries
$\partial_\pm\bar S = \bigl\{(t,x)\in\partial\bar S:(t+1,x\mp1)\in S\bigr\}$ and
$\partial^\pm\bar S = \bigl\{(t,x)\in\partial\bar S:(t-1,x\mp1)\in S\bigr\}$.

\index{boundary condition}
\index{birth process}
Consider $\{\xi_y\}_{y\in S}$, $\xi_y\geqslant0\ \forall {y\in S}$, the \emph{particle birth process} in $S$.
$\zeta^\circ$ is the \emph{boundary condition}, or the \emph{particle flow entering} $S$, as defined on p.~\pageref{page:zetacirc}.
One takes $\eta_y\geqslant0,\ y\in S$, the particle flow inside $S$ and exiting $S$, defined by~(\ref{eq:2.3}),(\ref{eq:2.4}).

\index{flow field}
\index{etazetaxi@$\bar\eta(\zeta^\circ,\xi)$}
Define $\bar\eta = \bigl\{\eta(e)\geqslant 0,\ e\in{\cal E}(\bar S)\bigr\}$, called the \emph{flow field} in $\E(\bar S)$ associated with $(\zeta^\circ,\xi)$.
As in the discrete case, there is a 1-1 correspondence between $(\zeta^\circ,\xi)$, $(\zeta^\circ,\eta)$, and $\bar\eta$.
We shall write $\bar\eta(\zeta^\circ,\xi)$ to denote the flow field $\bar\eta$ corresponding to the birth process $\xi$ and the flow $\zeta^\circ$ entering the boundary of $S$, analogous for $\bar\eta(\zeta^\circ,\eta)$.
\index{flow field}

In general, any nonnegative field $\bar\eta$ that satisfies
\begin{equation}
 \label{eq:lcons}
 n^+-n^-=m^+-m^-
\end{equation}
for each vertex $y$, where $n^+=\eta(e_y^\swarrow)$, $n^-=\eta(e_y^\nwarrow)$, $m^+=\eta(e_y^\nearrow)$ and $m^-=\eta(e_y^\searrow)$, is a flow field. Such a flow field may be defined either on a hexagonal domain $\E(\bar S)$ or on all $\E(\tilde\Z^2)$. In the former case there always exists a unique pair $(\zeta^\circ,\xi)$ such that $\bar\eta=\bar\eta(\zeta^\circ,\xi)$. In the later case such an expression does not make sense, although it is possible to determine $\xi(\bar\eta)$ for a given flow field $\bar\eta$. In Section~\ref{sec:lext} we shall discuss the extension of $\bar\eta$ to all of $\E(\tilde\Z^2)$ in the same spirit as for the geometric broken line process.

\subsection{The continuum broken line process}

\index{atoms}
Let $S$ be a fixed rectangular domain in $\tilde\Z^2$.
For a given flow field $\bar \eta$ in $\E(\bar S)$, we shall define a process of lines whose elementary constituents, called \emph{atoms}, are given by a point $p$ standing at a edge $e$, i.e., they are pairs of the form $(e,p),e\in \E(\bar S), p\in\bigl(0,\eta(e)\bigr]$.

We associate two atoms and write $(e_1,p_1)\sim(e_2,p_2)$ according to the following rules:

\index{association rules!continuous case}
Case 1: $e_1 = e_y^\swarrow$, $e_2 = e_y^\nwarrow$; $p_1\in\bigl(0,\eta(e_1)\wedge\eta(e_2)\bigr],\ p_2 = p_1$;

Case 2: $e_1 = e_y^\searrow$, $e_2 = e_y^\nwarrow$; $p_2\in\bigl(\eta(e_y^\swarrow),\eta(e_2)\bigr],\ p_1 = p_2 - \eta(e_y^\swarrow)$;

Case 3: $e_1 = e_y^\swarrow$, $e_2 = e_y^\nearrow$; $p_1\in\bigl(\eta(e_y^\nwarrow),\eta(e_1)\bigr],\ p_2 = p_1 - \eta(e_y^\nwarrow)$;

Case 4: $e_1 = e_y^\searrow$, $e_2 = e_y^\nearrow$; $p_i\in\bigl(\eta(e_i)-\xi_y,\eta(e_i)\bigr],\ \eta(e_1) - p_1 = \eta(e_2) - p_2$;

Notice that, given a vertex $y$, each atom standing at an edge incident to $y$ from above is associated with exactly one atom standing at another edge incident to $y$ from below and vice versa.
Notice also that `$\sim$' is not transitive.

By $J$ we will always denote an interval of the form $(p,p']$ and $|J|=p'-p$. When $p\geqslant p'$, $J=\emptyset$ and $|J|=0$. We associate two intervals of atoms standing at adjacent edges and write $(e_1,J_1)\sim(e_2,J_2)$ if for all $p_1\in J_1$ there is $p_2\in J_2$ such that $(e_1,p_1)\sim(e_2,p_2)$ and vice versa.
Notice that in this case $|J_1|=|J_2|$.

\index{broken line}
We define a \emph{broken line} as a tuple $\gamma = (y_0,e_1,J_1,y_1,e_2,J_2,y_2,\dots,e_n,J_n,y_n)$, $n\geqslant1$, such that 
\begin{equation}
 \label{eq:lbroktr}
\begin{array}{l}
e_i=\langle y_{i-1},y_i\rangle \\
x_i = x_{i-1}+1\\
t_i = t_{i-1}\pm 1
\end{array},
\end{equation}
and $y_i=(t_i,x_i)$ for $i=1,\dots,n$, and such that $|J_1|=|J_2|=\cdots=|J_n|$. If $J_i=\emptyset$ we identify $\gamma=\emptyset$.

\index{broken line!weight of}
We define next the \emph{weight} of a broken line $w(\gamma)$. If $\gamma\ne\emptyset$, we put $w(\gamma) = |J_1|>0$, otherwise we let $w(\gamma)=0$.

\index{broken line!domain of}
\index{domain|see{broken line domain}}
The \emph{domain} of $\gamma$, $D(\gamma)$, is given by $D(\gamma) = \{x_0,\dots,x_n\}$. Since for each $x\in D(\gamma)$ there is a unique $t$ such that $(t,x)\in\gamma$ we can denote such $t$ by $t(x)$.
(Notice that we abuse the symbol `$\in$' here, since $\gamma$ is a tuple and not a set.)

\index{broken trace}
\index{broken line!trace of}
For $\gamma = (y_0,e_1,J_1,y_1,e_2,J_2,y_2,\dots,e_n,J_n,y_n)$ we define what is called its \emph{trace} by $\ell(\gamma) = (y_0,e_1,y_1,e_2,y_2,\dots,e_n,y_n)$.
In general, any tuple $\ell$ satisfying (\ref{eq:lbroktr}) will be called a \emph{broken trace}.
Since either $(y_0,\dots,y_n)$ or $(e_1,\dots,e_n)$ are sufficient to determine $\ell$, we shall refer to any of these representations without distinction.

\index{broken trace!domain of}
\index{domain|see{broken trace domain}}
The \emph{domain} of a broken trace $\ell$ is defined as $D(\ell) = \{x:(t,x) \in \ell\}$. (Again we abuse the symbol `$\in$'.)
We also define $I(\ell)=\{t(x):x\in D(\ell)\}$.
It follows from~(\ref{eq:lbroktr}) that $D(\ell)$ is convex, i.e., it is of the form $\{a,a+1,a+1,\dots,b-1,b\}$ and that the same holds for $I(\ell)$.

\index{ls@$\ell\subseteq\bar S$}
\index{CS@$C(S)$}
We write $\ell\subseteq \bar S$ 
if $y_0\in S\cup\partial_- \bar S\cup\partial^- \bar S$, $y_n\in S\cup\partial_+ \bar S\cup\partial^+ \bar S$ and $y_1,\dots,y_{n-1}\in S$, 
or, equivalently, if $e_1,\dots,e_n\in\E(\bar S)$.
We say that $\ell$ \emph{crosses} $S$ if $\ell\subseteq \bar S$ and $y_0,y_n\in\partial\bar S$. Let $C(S)=\{\ell\subseteq\bar S:\ell \mbox{ crosses } S\}$. We write $\gamma\subseteq\bar S$ if $\ell(\gamma)\subseteq\bar S$ and we say that $\gamma$ \emph{crosses} $S$ if $\ell(\gamma)$ crosses $S$.

\index{Beta@$B(\bar\eta)$|see{broken line associated to $\bar\eta$}}
\index{broken line!associated to $\bar\eta$}
We say that the broken line $\gamma\subseteq\bar S$ is associated with a flow field $\bar\eta$ defined on $\E(\bar S)$ if $(e_{i-1},J_{i-1})\sim(e_i,J_i)$ for $i=1,\dots,n$. $B(\bar\eta)$ will denote the set of all broken lines associated with $\bar\eta$.

For $V\subseteq\tilde\Z^2$, and $\ell=(y_0,\dots,y_n)$, define $\ell\cap V = \{y_i:y_i\in V,i=0,\dots,n\}$.

\index{Lgamma@$L(\gamma)$}
\index{Ll@$L(\ell)$}
For a given broken line $\gamma$ associated with the field $\bar\eta(\zeta^\circ,\xi)$, its left corners $(t,x)$ correspond to part of the particle birth $\xi_{t,x}$ at $(t,x)$. Given a broken trace $\ell$, we denote by $L(\ell)$ the set of left corners of $\ell$, i.e., the points $(t,x)\in \ell$ such that $(t+1,x\pm1)\in\ell$. Also let $L(\gamma)=L(\ell(\gamma))$.

\index{xiell@$\xi(\ell)$}
\index{xigamma@$\xi(\gamma)$}
We define the fields $\xi(\ell)$ and $\xi(\gamma)$ in $\tilde\Z^2$ by  
\begin{equation}
 \label{eq:lxidef}
 [\xi(\ell)]_y =\I_{L(\ell)}(y)
\end{equation}
and $\xi(\gamma) = w(\gamma) \I_{L(\gamma)}$.

\index{etal@$\bar\eta(\ell)$}
\index{etagamma@$\bar\eta(\gamma)$}
In the same fashion, the extremal points of a broken line that crosses $S$ correspond to a particle flow entering or exiting $S$. So we also define
$$
  [\zeta^+(\ell)]_{t,x} = 
  \begin{cases}
    1, & (x,t)\in \ell\cap S \mbox{ and } (t-1,x-1)\in\ell\cap\partial\bar S  \\ 
    0 & \mbox{otherwise}
  \end{cases},
$$
$$
  [\zeta^-(\ell)]_{t,x} = 
  \begin{cases}
  1, & (x,t)\in \ell\cap S \mbox{ and }(t-1,x+1)\in\ell\cap\partial\bar S  \\ 
   0 & \mbox{otherwise}
  \end{cases},
$$
take the corresponding $\zeta^\circ$ and denote by $\zeta^\circ(\ell)$.
Define $\zeta^\circ(\gamma) = w(\gamma)\zeta^\circ(\ell(\gamma))$.
Also define $\bar\eta(\ell)$ by $\bar\eta(\ell)(e)=\I_{e\in\ell}$ and $\bar\eta(\gamma)=w(\gamma)\bar\eta(\ell(\gamma))$.
It is easy to see that $\bar\eta(\ell)=\bar\eta\bigl(\zeta^\circ(\ell),\xi(\ell)\bigr)$ for $\ell\in C(S)$, analogous for $\gamma$.

\index{broken line!maximal}
\index{gammal@$\gamma(\ell)$|see{broken line, maximal}}
Let $\bar\eta$ be given and fix some $\ell=(y_0,e_1,y_1,e_2,y_2,\dots,e_n,y_n)\subseteq\bar S$.
There is one maximal broken line that has trace $\ell$ and is associated with the field $\bar\eta$, which will be denoted $\gamma(\ell)$, the dependence on $\bar\eta$ is omitted.
By this we mean that there exist unique $J_1, J_2, \dots,J_n$ such that $\gamma(\ell)=(y_0,e_1,J_1,y_1,\dots,e_n,J_n,y_n) \in B(\bar\eta)$ and such that any $p_1,p_2,\dots,p_n$ with property $(e_1,p_1)\sim(e_2,p_2)\sim\cdots\sim(e_n,p_n)$ must satisfy $p_i\in J_i,\ i=1,\dots,n$.
The proof of this fact is shown in Appendix~\ref{sec:maximal}.

\index{wl@$w(\ell)$}
Let $w(\ell)$ denote the maximum weight of a broken line in $B(\bar\eta)$ that has trace $\ell$, which is given by $w\bigl(\gamma(\ell)\bigr)$.
The dependence on the field $\bar\eta$ is omitted when it is clear which field is being considered, otherwise we shall write $w_{\eta}(\ell)$.
Notice that $D\bigl(\gamma(\ell)\bigr) = D(\ell)$ if $w_\eta(\ell)>0$, and $D\bigl(\gamma(\ell)\bigr) =\emptyset$ otherwise.

\index{ll@$\ell\subseteq\ell'$}
We write $\ell\subseteq\ell'$ if $\ell'=(y_0,e_1,y_1,e_2,y_2,\dots,e_n,y_n)$ and,
for some $0\leqslant a<b\leqslant n$, $\ell=(y_a,e_{a+1},y_{a+1},e_{a+2},y_{a+2},\dots,e_b,y_b)$;
or, equivalently, if $D(\ell)\subseteq D(\ell')$ and $t(x)=t'(x)$ for all $x\in D(\ell)$.

Notice that $w(\ell_1) \geqslant w(\ell_2)$ when $\ell_1\subseteq\ell_2$. This is due to the successive branching caused by birth/collision that makes longer lines become thinner. It is even possible that for $\ell_1\subseteq\ell_2\subseteq\cdots\subseteq\ell_n\subseteq\cdots$ we have $w(\ell_n)\downarrow0$ and the limiting object could have one atom but null weight.

\index{ll@$\ell \succeq \ell'$}
For broken traces $\ell,\ell'$, we write $\ell \succeq \ell'$ if $\ell$ is to the right of $\ell'$. This means that $t(x)\geqslant t'(x)$ for all $x\in D(\ell)\cap D(\ell')$ and that $t(x)\geqslant t'(x')$ for some $x\in D(\ell),x'\in D(\ell')$.
(The last condition makes sense in case $D(\ell)\cap D(\ell')=\emptyset$.)
In general this relation is neither antisymmetric nor transitive.
Write $\ell \succ \ell'$ if $\ell \succeq \ell'$ and $\ell \ne \ell'$.

\begin{lemma}
\label{lemma:lorder}
 If $S$ is a rectangular domain, the following assertions hold:
\begin{enumerate}
 \item \label{item:partorder}
Relation $\succeq$, restricted to $\bigl\{\ell \in C(S)\bigr\}$, is a partial order.
 \item \label{item:ellextreme}
The elements of $C(S)$ are extremal in the following sense. If $\ell\in C(S)$, $\ell'\subseteq\bar S$ and $\ell\subseteq\ell'$, then $\ell=\ell'$.
 \item \label{item:useless}
Relations $\subseteq$ / $\succeq$ have a certain concavity in the following sense. Suppose that $\ell\subseteq\ell_1$ and $\ell\subseteq\ell_3$ for some $\ell\subseteq\bar S$ and $\ell_1\preceq\ell_2\preceq\ell_3\in C(S)$; then $\ell\subseteq\ell_2$.
\end{enumerate}
\end{lemma}
It should be clear that the items above hold true. The rigorous proof of this lemma consists of straightforward but tedious verifications and is postponed until Appendix~\ref{sec:lorder}.

\begin{lemma}
 \label{lemma:comparable}
 Let $S$ be a rectangular domain. If $w_\eta(\ell)>0$ and $w_\eta(\ell')>0$ for some flow field $\bar\eta$, then $\ell$ and $\ell'$ are comparable, that is, $\ell\succeq\ell'$ or $\ell'\succeq\ell$.
\end{lemma}
Lemma~\ref{lemma:comparable} is a consequence of the way the association rules have been defined.
In order to prove it, one keeps applying such association rules and the result follows by induction -- see Appendix~\ref{sec:comparable}.

\index{flow field!decomposition}
The theorem below is fundamental. It says we can decompose a given flow field in many other smaller fields, each one corresponding to one of the maximal broken lines that cross $S$ and is associated to the field. In this case the original flow field and all of its features, namely the birth process, the boundary conditions and the weight it attributes to broken traces, are additive in the sense that each of these is obtained by summing over the smaller fields.
On the other hand it tells that, given an ordered set of broken lines, it is possible to combine the corresponding fields and thereby determine the flow field that is associated to them.

\begin{theo}
\label{theo:etaell}
Let $S$ be a rectangular domain.

Given $\ell_1\prec\cdots\prec\ell_n\in C(S)$ and $w_1,\dots,w_n>0$, there is a unique flow field $\bar\eta$ in $\E(\bar S)$ such that
\begin{equation}
  \label{eq:lbtcarac}
  w_\eta(\ell) =
  \begin{cases} 
     w_j, & \mbox{if } \ell=\ell_j \mbox{ for some } j \\ 0, & \mbox{otherwise}
  \end{cases}
  \qquad\quad
  \mbox{for any $\ell\in C(S)$.}
\end{equation}
Moreover, for $\zeta^\circ,\xi$ such that $\bar\eta=\bar\eta(\zeta^\circ,\xi)$, the following decompositions hold:
\begin{equation}
 \label{eq:ldecomp}
\bar\eta = \sum_jw_j\bar\eta_j, \quad \zeta^\circ = \sum_jw_j\zeta^\circ_j, \quad \xi = \sum_jw_j\xi_j,
\end{equation}
where $\bar\eta_j=\bar\eta(\ell_j)$, $\zeta^\circ_j = \zeta^\circ(\ell_j)$ and $\xi_j = \xi(\ell_j)$. Furthermore,
\begin{equation}
  \label{eq:lweightcalc}
  w_\eta(\ell) = \sum_{\ell'\in C(S)} w_\eta(\ell') \I_{\ell\subseteq\ell'}
 \qquad\quad
 \mbox{for any } \ell\subseteq\bar S.
\end{equation}

Conversely, given a flow field $\bar\eta$, there are unique sets $\{\ell_j\}$ and $\{w_j>0\}$ that satisfy~(\ref{eq:lbtcarac}). The set $\{\ell_j\}$ is totally ordered by the relation `$\succeq$' and (\ref{eq:ldecomp})-(\ref{eq:lweightcalc}) hold in this case.
\end{theo}

\textbf{Proof.}
We start by the converse part, first proving that (\ref{eq:lweightcalc}) holds for any flow field $\bar\eta$, which is the most laborious work. The proof of (\ref{eq:lweightcalc}) is a mere formalization of the construction described below,
whereas (\ref{eq:lbtcarac})~and~(\ref{eq:ldecomp}) are immediate consequences, as discussed afterwards.

Then it will suffice to show that, given any pair of sets $\{\ell_1\prec\ell_2\prec\cdots\prec\ell_M\}$ and $\{w_1,\dots,w_M>0\}$,
there is some flow field $\bar\eta$ satisfying~(\ref{eq:lbtcarac}).
Uniqueness of such flow field follows from the converse part.
Again, by the converse part,~(\ref{eq:ldecomp})~and~(\ref{eq:lweightcalc}) will hold as well,
completing the proof of the theorem.
In order to prove that~(\ref{eq:lbtcarac}) holds we basically have to see that the same construction can be reversed.

Though the construction looks simple, several equivalent representations of a flow field may be seen on the same picture and the theorem will be deduced from these representations.
We remind that the sole condition for a field $\eta(e)\geqslant0$, $e\in\E(\bar S)$ to be a flow field is the conservation law below:
\begin{equation}
 \label{eq:lcons2}
  \eta(e_y^\nwarrow) + \eta(e_y^\nearrow) = \eta(e_y^\swarrow) +
 \eta(e_y^\searrow)\qquad\forall\ y\in S.
\end{equation}
The construction we start describing now strongly relies on this fact.

\begin{figure}[!ht]
\begin{center}
 \includegraphics[scale=.45]{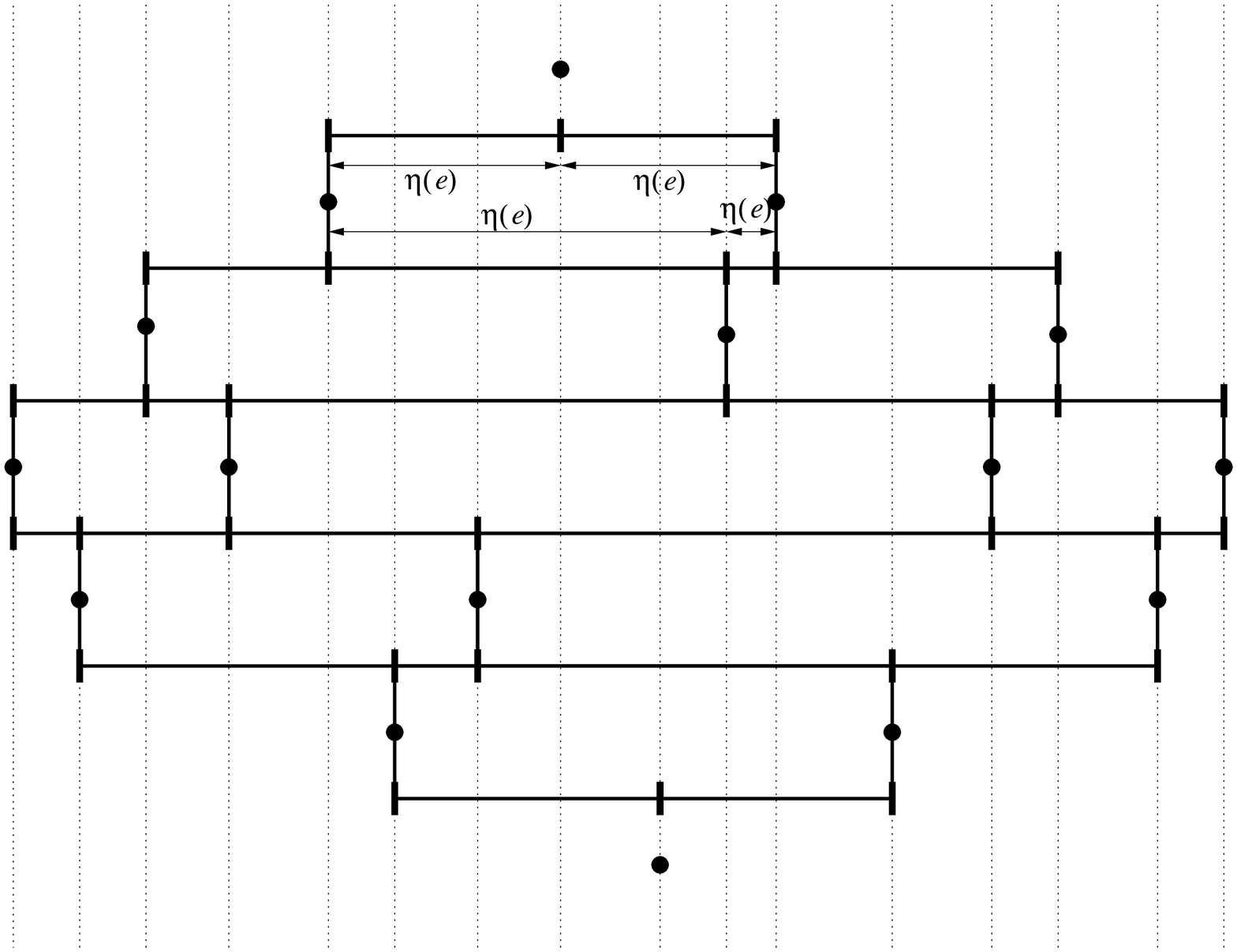}
\end{center}
\caption{
\small
construction of the brick diagram. In this example $S$ is a $3\times3$ rectangular domain. There are 6 horizontal lines with a total of 24 intervals forming 9 bricks. The diagram is divided into 15 strips.}
\label{fig:brick}
\end{figure}

First we plot on a horizontal line two adjacent intervals whose lengths correspond to the flow $\eta$ on the two topmost edges of $\E(\bar S)$, i.e., the two edges incident from above to the topmost site $y\in S$. Since the intervals are adjacent we may do that by marking three points on this line. On the next line (parallel to and below the first one), we plot four adjacent intervals having lengths corresponding to the flow on the four edges incident from above to the two sites on the 2nd row of $S$. It follows from~(\ref{eq:lcons2}) that one can position these intervals in a way that its 2nd and 4th points stay exactly below the 1st and 3rd points of the first line. By linking these two pairs of points we get a ``brick'' that corresponds to the topmost site $y$ of $S$.
The width of this brick is equal to the quantity expressed in~(\ref{eq:lcons2}), its top face is divided into two subintervals having lengths $\eta(e_y^\nwarrow)$ and $\eta(e_y^\nearrow)$ and its bottom face is divided into intervals of lengths $\eta(e_y^\swarrow)$ and $\eta(e_y^\searrow)$ -- see Figure~\ref{fig:brick}.

\index{brick diagram}
We carry on this procedure for the 3rd horizontal line, thereby getting two bricks that correspond to the two sites on the second row of $S$. We keep doing this construction until we mark intervals corresponding to the two bottommost edges of $\E(\bar S)$. In the final picture we have one brick corresponding to each site of $S$. Each pair of bricks on consecutive levels that have a common (perhaps degenerate) interval on their boundaries corresponds to adjacent sites $y$ and $y'$ in $S$; the length of this common interval equals $\eta\bigl(\langle y,y'\rangle\bigr)$. Once all the intervals have been plotted at the appropriate position, forming all the bricks, we draw a dotted vertical line passing by each point that was delimiting these intervals. By doing this we divide the whole diagram into strips, completing the construction of the \emph{brick diagram} as shown in Figure~\ref{fig:brick}.

Each strip corresponds to a maximal broken line that crosses $S$ and the weight of this broken line equals the width of the strip.
The sites/edges that compound the broken line correspond to the bricks/intervals the strip intersects.
Given a broken trace $\ell\subseteq\bar S$, we can determine the maximal broken line $\gamma(\ell)$ that passes through $\ell$ by looking which strips pass by all the sites/edges of $\ell$ (i.e., their corresponding bricks/intervals); the weight of this broken line is obtained by summing the width of such strips.
\begin{figure}[!ht]
\begin{center}
 \includegraphics[width=55mm]{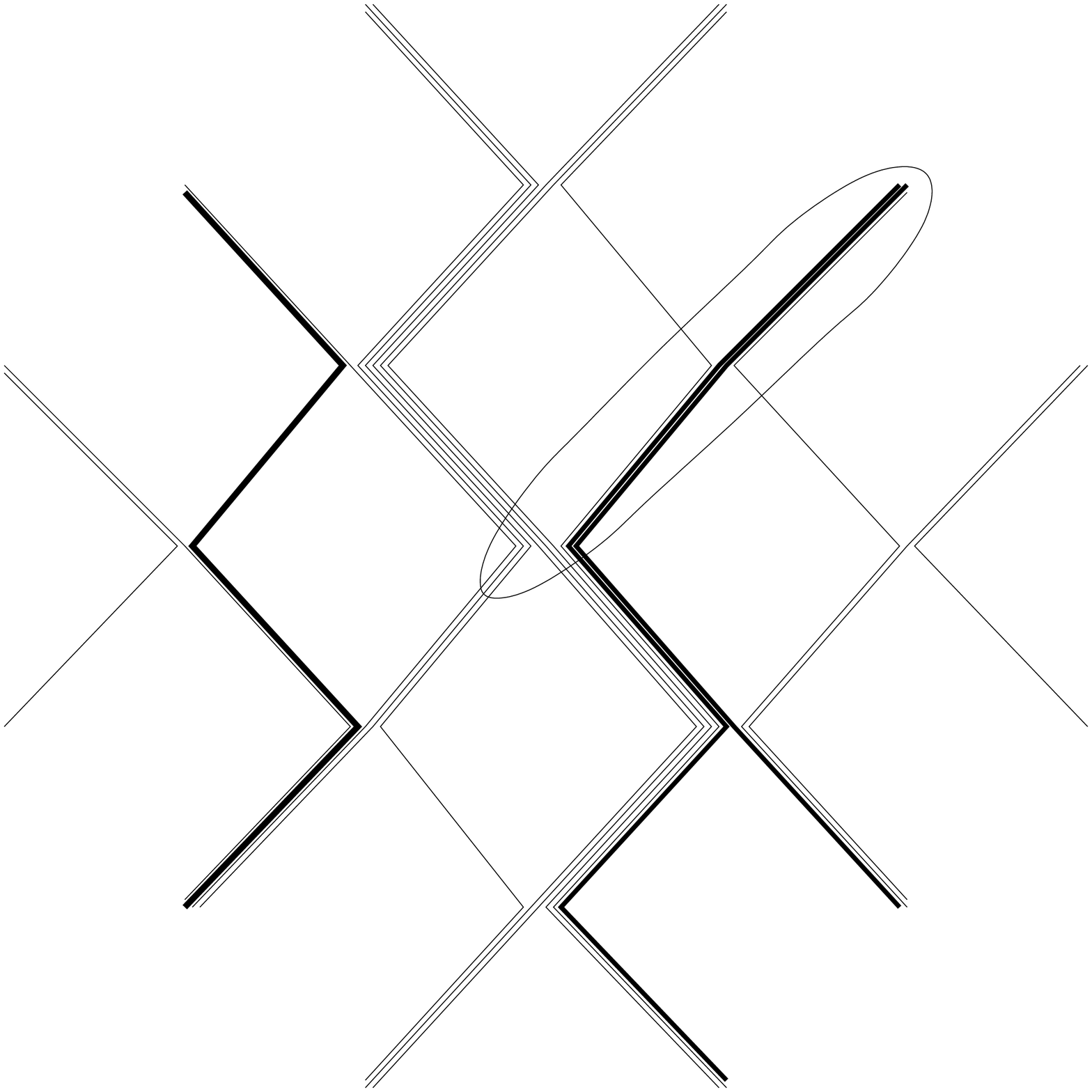}
 \hspace{0.4cm} 
 \includegraphics[width=55mm]{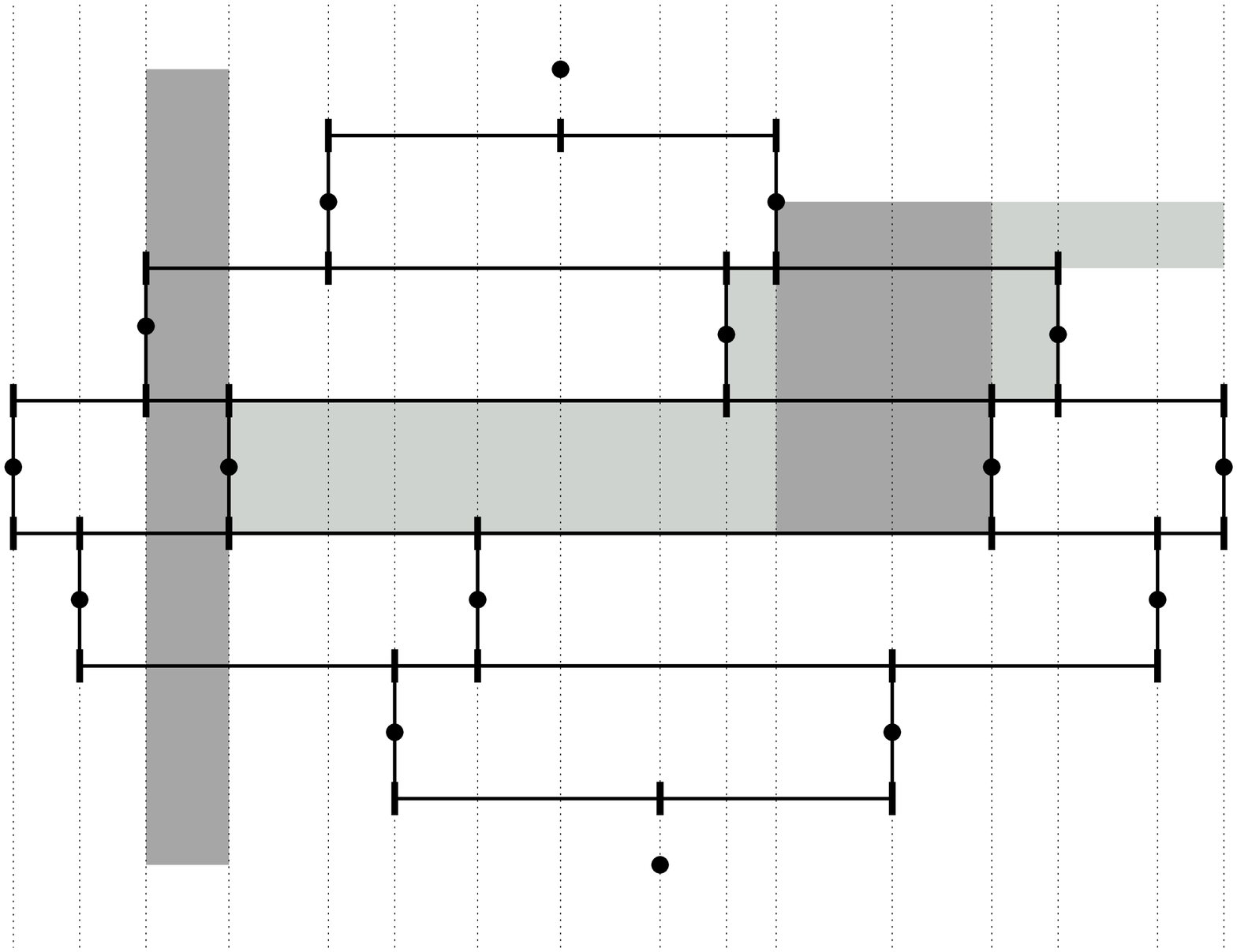}
\end{center}
\caption{
\small
broken lines that cross $S$ and calculation of the weight of a given broken trace. On the left we have a random configuration of broken lines on the $3\times3$ domain $S$ and on the right we have the brick diagram from which they arose.
The 3rd strip in the brick diagram is highlighted and the corresponding broken line appears in bold on the left. Also, for a broken trace $\ell$ that starts at the center of $S$, goes northeast, and then exits $S$ going northeast again, we have determined which broken lines in $C(S)$ contain $\ell$; they are the 11th and 12th ones and they also appear in bold. In this case $w(\ell)$ equals the sum of the widths of the 11th and 12th strips.}
\label{fig:domainlines}
\end{figure}
See Figure~\ref{fig:domainlines}.

So, given a flow field $\bar\eta$ we construct the brick diagram from which the broken line configuration is deduced, with the desired property that such broken lines satisfy~(\ref{eq:lbtcarac}) and~(\ref{eq:lweightcalc}). 

On the other hand, let a set of well ordered broken lines, that is, $\{\ell_1\prec\ell_2\prec\cdots\prec\ell_M\}$ and $\{w_1,\dots,w_M>0\}$, be given.
One can consider the corresponding broken line diagram, from which one constructs the brick diagram and the later gives a flow field satisfying~(\ref{eq:lbtcarac}).

In a first reading, one is encouraged to understand the above description with the corresponding figures. The more interested reader will find a detailed proof in Appendix~\ref{sec:etaell}.

Since~(\ref{eq:lweightcalc}) has been proven for any flow filed, notice that uniqueness of $\{\ell_j\}$ and $\{w_j\}$ is trivial by definition and well ordering follows from Lemma~\ref{lemma:comparable}.
Finally, (\ref{eq:ldecomp}) holds because of~(\ref{eq:lweightcalc}) when we write each process in terms of weights of certain broken lines.
For $\zeta^\circ,\xi,\eta$ such that $\bar\eta=\bar\eta(\zeta^\circ,\xi)=\bar\eta(\zeta^\circ ,\eta)$ we have $\zeta^-_y = \eta(e_y^\nwarrow)$, $\zeta^+_y = \eta(e_y^\swarrow)$, $\eta^-_y = \eta(e_y^\searrow)$, $\eta^+_y = \eta(e_y^\nearrow)$, and $\xi_y=\eta_y^+\wedge\eta_y^-$. Given $e=\langle y,y'\rangle\in\E(\bar S)$, take $\ell(e)=(y,e,y')$. Then of course $\eta(e)=w_\eta\bigl(\ell(e)\bigr)$. Also, given $y\in S$, take $\ell^<(y) =(y_-,e_-,y,e_+,y_+)$, where $y=(t,x)$, $y_{\pm}=(t+1,x\pm1)$ and $e_\pm=\langle y,y_\pm\rangle$. It is also easy to see that $w_\eta\bigl(\ell^<(y)\bigr)=\xi_y$.
Also notice that $\zeta^-_y(\ell)=\I_{e_y^\nwarrow\in\ell}= \I_{\ell(e_y^\nwarrow)\subseteq\ell}$.
Now we put all the pieces together to get $\zeta_y^-=\eta(e_y^\nwarrow) = w_\eta\bigl(\ell(e_y^\nwarrow)\bigr)=\sum_jw_j\I_{\ell(e_y^\nwarrow)\subseteq\ell_j}= \sum_jw_j(\zeta_j^-)_y$, i.e., $\zeta^-=\sum_jw_j\zeta_j^-$. The other equalities are deduced similarly.

This completes the proof of the converse part and, as discussed above, of the theorem.
$\hfill \square$\\

\begin{cor}
 Let $S$ be a rectangular domain and let $\bar\eta(\zeta^\circ,\xi)$ be given. Then
\begin{equation}
 \label{eq:lhzx}
 \left[ \sum_{y\in\partial_-S} \zeta_y^+  + \sum_{y\in\partial^-S} \eta_y^- \right] = \left[
 \sum_{y\in\partial_+S} \zeta_y^- + \sum_{y\in\partial^+S} \eta_y^+ \right] = \sum_{\ell\in C(S)}w(\ell).
\end{equation}
\end{cor}

\textbf{Proof:}
For each $\ell\in C(S)$ there is exactly one $e\in\ell$ such that $y\in\partial^+S,y'\in \partial^+\bar S$ or $y\in\partial_+S,y'\in \partial_+\bar S$, where $e=\langle y,y'\rangle$.
So
$1 =
\sum_{(t,x)\in\partial_+S} \I_{\langle(t,x),(t-1,x+1)\rangle\in\ell}
+
\sum_{(t,x)\in\partial^+S} \I_{\langle(t,x),(t+1,x+1)\rangle\in\ell}
=
\sum_{y\in\partial_+S}\zeta_y^-(\ell)
+
\sum_{y\in\partial^+S}\eta_y^+(\ell)$.
Multiplying by $w(\ell)$ and summing over all $\ell\in C(S)$ we get
$\sum_{y\in\partial_+S} \zeta_y^- + \sum_{y\in\partial^+S} \eta_y^+ = \sum_{\ell\in C(S)}w(\ell)$.
The proof of 
$\sum_{y\in\partial_-S} \zeta_y^+  + \sum_{y\in\partial^-S} \eta_y^- = \sum_{\ell\in C(S)}w(\ell)$ is similar.
$\hfill \square$\\

Denote by $H_S(\zeta^\circ,\xi)$ the quantity expressed in (\ref{eq:lhzx}), if $\zeta^\circ=0$ we just write $H_S(\xi)$.

\index{last passage percolation}
\index{percolation|see{last passage percolation}}
For a birth field $\xi\geqslant0$ on $S$ we define the \emph{last passage percolation value} $G_S(\xi)$ as the maximum of $\sum_{y\in\pi} \xi_y$ over all $\pi\in\Pi_S$, where
$$
 \Pi_S=\{\pi=(y_0,\dots,y_m):y_0\in\partial_0 S, y_m\in\partial_1 S, t_{i+1}=t_i+1, x_{i+1}=x_i\pm1\}.
$$
When it is clear which rectangular domain we are referring to we shall drop the subscript $S$ of $H_S$, $G_S$ and $\Pi_S$.

The next proposition illustrates the connection between broken lines and passage time. Furthermore, its proof gives an explicit algorithm for determining the optimal path (which is a.s. unique when $\xi$ has a continuous distribution).
\begin{prop}
 \label{prop:llpp}
 Let $S$ be a rectangular domain. The last passage percolation value $G_S(\xi)$ is given by the sum of the weight of the broken lines associated to the corresponding birth field: $G_S(\xi)=H_S(\xi)$.
\end{prop}

\begin{figure}[!htb]
\psfrag{t}{$t$}
\psfrag{x}{$x$}
\begin{center}
 \includegraphics[width=5.0cm]{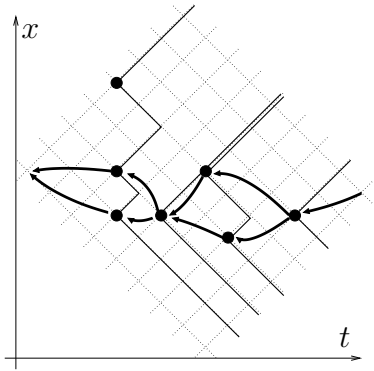}
\end{center}
\caption{
\small
construction of the maximal path. From $\xi$ one constructs the flow field and by following algorithm~(\ref{eq:algorithm}) one gets an optimal path. The theory developed for broken lines and Theorem~\ref{theo:etaell} guarantee this is indeed optimal.
The algorithm forbids the path to cross any broken line; for that purpose it suffices to require that it does not `cross' the flow field.}
\label{fig:lpp}
\end{figure}

\textbf{Proof:}
The proof consists on formalizing the following argument.
On one hand, an oriented path $\pi$ connecting $\partial_0S\to\partial_1S$ can cross at most one left corner of each broken line.
On the other hand, it is possible to assemble the path backwards, following a local rule that does not miss any broken line.
This is possible because the broken lines never cross each other.
See Figure~\ref{fig:lpp}.

Take $\bar\eta = \bar\eta(\zeta^\circ=0,\xi)$ and write $\{\ell\in C(S):w_\eta(\ell)>0\}=\{\ell_1\prec\ell_2\prec\cdots\prec\ell_M\}$. Now by (\ref{eq:lxidef}) and (\ref{eq:ldecomp}) it follows that $G(\xi)$ is the maximum over all $\pi\in\Pi$ of
$$
  \sum_{y\in\pi} \xi_y = \sum_{y\in\pi} \sum_j w(\ell_j)[\xi(\ell_j)]_y = 
  \sum_j w(\ell_j) \sum_{y\in\pi} \I_{L(\ell_j)}(y).
$$
Since the paths $\pi$ are oriented, they cannot intersect more than one left corner of each broken line, hence $\sum_{y\in\pi} \I_{L(\ell_j)}(y) \leqslant 1$ for each $j$.
We shall exhibit an algorithm for constructing a path $\pi^*$ that satisfies
\begin{equation}
 \label{eq:intersect}
 \sum_{y\in\pi^*} \I_{L(\gamma_j)}(y) \geqslant 1,
\end{equation}
which completes the proof.

The path $\pi^*$ is constructed by the following rule.
Let $y_m \in \partial_1 S$. For $i = m,m-1,m-2\dots,2,1$, let $t_{i-1}=t_i-1$ and
\begin{equation}
 \label{eq:algorithm}
  x_{i-1} = 
  \begin{cases}
     x_i-1,& \eta(e_{y_i}^\swarrow)\geqslant\eta(e_{y_i}^\nwarrow), \\ 
     x_i+1, & \mbox{otherwise}.
  \end{cases}
\end{equation}

It remains to show~(\ref{eq:intersect}),
i.e., that $\pi^*$ intersects $L(\ell_j)$ for each $j=1,\dots,M$.
Fix $j$ and write $\ell$ for $\ell_j$.
Assume without proof that $\pi^*$ intersects $\ell$ at some point
(a complete proof is shown in Appendix~\ref{sec:intersect}).
Take $n=\min\{i:y_i\in\ell\}$. If $y_0\in\ell$ we have $y_0\in L(\ell)$ due to the fact that $\zeta^\pm_{y_0}=0$.
So suppose $n>0$. By construction $x_{n-1}=x_n\pm1$; assume for simplicity $x_{n-1}=x_n+1$. Now $t_{n-1}=t_n-1$ and $n$ is minimal, so $(t_n-1,x_n+1)\not\in\ell$.
Since $y_n\not\in\partial\bar S$, $x_n+1\in D(\ell)$ and thus $(t_n+1,x_n+1)\in\ell$.
Now if $(t_n-1,x_n-1)$ were in $\ell$ there would be $p_1\in(0,\eta(e_{y_n}^\swarrow)]$, $p_2\in(0,\eta(e_{y_n}^\nearrow)]$ associated by Case~3, which implies $\eta(e_{y_n}^\swarrow)>\eta(e_{y_n}^\nwarrow)$, contradicting the choice of $x_{n-1}=x_n+1$, so $(t_n-1,x_n-1)\not\in\ell$. But also $x_n-1\in D(\ell)$, thus $(t_n+1,x_n-1)\in\ell$ and therefore $y_n\in L(\ell)$.
$\hfill \square$\\

\begin{cor}
 \label{lemma:lmonot}
 $H(\zeta^\circ,\xi)$ in nondecreasing in $(\zeta^\circ,\xi)	$.
\end{cor}

There is a nice direct proof for the above corollary that uses Theorem~\ref{theo:etaell}, and a much shorter one that relies on Proposition~\ref{prop:llpp}. The argument for the short proof will appear in the proof of Theorem~\ref{theo:exp}, just before~(\ref{eq:somenopt}), so it is not repeat it here. The next corollary is indeed an immediate consequence of Proposition~\ref{prop:llpp}, although not at all obvious when one thinks of the definition of $H_S(\xi)$ by itself.

\begin{cor}
 Let $S$ be a rectangular domain and $\xi\geqslant0$ a particle birth process. Consider the map $M:(t,x)\mapsto(-t,x)$, $S'=M(S)$, $\xi'(t,x)=\xi(-t,x)$. Then $H_S(\xi)=H_{S'}(\xi')$.
\end{cor}

\subsection{Reversibility}
\label{sec:lext}

In analogy with the discrete case, given the distributions of $\zeta^+$, $\zeta^-$ and $\xi$ satisfying (\ref{eq:2.15}), the evolution of the particles described as the evolution of the associated Markov chain is reversible.
More precisely, reversibility is determined by the following relation
\begin{equation}
\label{eq:lrev}
 \begin{array}{c}
  \pi_1(\dd m^+)\pi_2(\dd m^-)q(\dd n^+,\dd n^-|m^+,m^-) = \\ = \pi_1(\dd n^+)\pi_2(\dd n^-)q(\dd m^+,\dd m^-|n^+,n^-).
 \end{array}
\end{equation}

In general, all we have to check is that the relation (\ref{eq:lrev}) holds for the distribution $q(\dd n^+,\dd n^-|m^+,m^-)$, where the $q$ is the kernel that makes (\ref{eq:lcons})
happen a.s. and $(n^+ \wedge n^-)$ be independent of $m^+$ and $m^-$. (It is given by $\xi_y$, having law $\pi_3$.)

Thus any triple $\pi_1,\pi_2,\pi_3$ that satisfies (\ref{eq:lrev}) will define a family of measures $\{P_S:S\mbox{ hexagonal domain}\}$ which is consistent, i.e., satisfies~(\ref{eq:2.17}). In this case $\bar\eta$ can be consistently extended to ${\cal E}(\tilde \Z^2)$.
As a consequence, if one takes $(\zeta^+_y)_{y\in\partial_-S}$ i.i.d. distributed as $\pi_1$,
$(\zeta^-_y)_{y\in\partial_+S}$ i.i.d. having law $\pi_2$, and $(\xi_y)_{y\in S}$ i.i.d. with law $\pi_3$, then $(\eta_y^+)_{y\in\partial^+S}$ will be i.i.d. with law $\pi_1$ and 
$(\eta_y^-)_{y\in\partial^-S}$ will be i.i.d. with law $\pi_2$.

The geometric broken line process described in Section~\ref{sec:discretebl} is a particular case of $\eta^+\sim{\rm Geom}(\lambda_1)$, $\eta^-\sim{\rm Geom}(\lambda_2)$, $\xi\sim{\rm Geom}(\lambda_1\lambda_2)$, taking $\lambda_1=\lambda_2=\lambda$.
Below we characterize the distributions that satisfy~(\ref{eq:lrev}).

Since the particle system always evolves by keeping the relation (\ref{eq:lcons}), we can parametrize this hyperplane in $\R_+^4$ by $T:\R_+^3\to\R_+^4$, $T(r,s,t)=(r,s,t+[r-s]^+,t+[r-s]^-)$.
The joint distribution of $(n^+,n^-,m^+,m^-)$ is given by the left hand side of (\ref{eq:lrev}). It can be obtained by taking the projection $T_*\mu$, where $\mu = \pi_1 \times \pi_2 \times \pi_3$.

Now consider the operator in $\R_+^4$ given by $L(x,y,z,w)=(z,w,x,y)$. Reversibility (\ref{eq:lrev}) in our case just means that $L$ preserves $T_*\mu$. This is equivalent to the fact that
\begin{equation}
 \label{eq:lmpres}
  R_*\mu = \mu,
\end{equation}
where $R=T^{-1}LT$. Writing $R$ explicitly gives $R(r,s,t)=(t+[r-s]^+,t+[r-s]^-,r\wedge s)$ and $R^2 = I$, the identity operator. Note that $B_1$, $\partial B_1$, $B_2$, given by $B_1=\{(r,s,t):r>s\}$, $B_2=\{(r,s,t):r<s\}$, remain invariant under $R$ and they form a partition of $\R^3_+$. Also note that $\{r\wedge s = 0\}$ is in bijection with $\{t=0\}$.

Another example that satisfies~(\ref{eq:lmpres}) and thus~(\ref{eq:lrev}) is when the $\pi_i$'s are respectively the exponential distributions $\exp(\alpha_i)$ with $\alpha_3=\alpha_1+\alpha_2$.
In this case it is a simple exercise to check (\ref{eq:lmpres}). Assume $A\subseteq B_1$, the expression for $R$ becomes much simpler and it is straightforward that $\mu(R^{-1}(A))=\mu(A)$.
For $A\subseteq B_2$ the situation is analogous and $\mu(\partial B_1)=0$, completing the proof of (\ref{eq:lmpres}).

The question of finding triples of measures on $\R_+$ that satisfy~(\ref{eq:lmpres}) is an interesting question on its own.
In discussion with A.~Ram\'\i rez we have shown the following.
If the measures are supported on $\Z_+$ and assign positive weight to each $n\in\Z_+$, then they must be geometric distributions.
If the measures have continuous positive densities, then they must be exponential distributions.
Of course a lot remains between being supported on $\Z_+$ and having continuous positive density, but nonetheless we believe that the only examples are these two, modulo adding or multiplying by constants.

\sectionmark{Application: last passage percolation}

\section
{An application: geometric and exponential last passage percolation}
\label{sec:appl}

\sectionmark{Application: last passage percolation}

It is easy to see that the last passage percolation model satisfies a law of large numbers by super-additivity arguments.
However it is interesting that for the two-dimensional model and for the special case of i.i.d. geometric or exponential passage time distributions there is an explicit expression for the limiting constant.
For the exponential distribution it was found by Rost~\cite{rost81} and for geometric case by Jockusch, Propp, and Shor~\cite{jockusch98}.
Large deviations were studied by Johansson~\cite{johansson00} (lower tails) and by Sepp\"al\"ainen~\cite{seppalainen98} (upper tails).
Fluctuations were studied in~\cite{johansson00}.

With the aid of the broken line theory developed in the previous sections it is possible to re-obtain the explicit constants for the law of large numbers.
We also prove exponential decay for the probability of deviations.
This application is based on the law of large numbers for the intersection process, on Proposition~\ref{prop:llpp} and on Corollary~\ref{lemma:lmonot}.

In the same spirit, O'Connell~\cite{oconnell00} also devises such constants by simple probabilistic arguments.
We note that our approach is self-contained, except for using of Cram\'er's theorem for large deviations of i.i.d. sums.
A proof of Burke's theorem is implicitly contained in our considerations of reversibility.

The construction consists on first choosing the appropriate distributions of the boundary conditions that (i) make the broken line process reversible and (ii) provide the correct asymptotic behavior, and then dropping the boundary conditions afterwards.

What we present is an alternative proof that could give some geometric insight to the model.
Besides, the broken-line approach provides an explicit, linear algorithm for determining the maximal path (see the proof of Proposition~\ref{prop:llpp}). Our results also show that the boundary conditions give no asymptotic contribution to the total flow of broken lines that cross a given domain when they have this suitable distribution.

For $N,M\in\N$, we define the last passage percolation value on the square $\{1,\dots,N\}\times\{1,\dots,M\}$ as the random number $G(N,M)$ given by the maximum sum of $\xi_{y_j}$ over all oriented paths $(y_1,\dots,y_{N+M-1})$ from $(1,1)$ to $(N,M)$. $G(N,M)$ is random because so are the $\xi_y$'s.

\begin{theo}
 \label{theo:exp}
 Suppose $\xi_y$,
 $y\in\N^2$,
 are i.i.d and distributed as $\exp(\alpha)$, $\alpha>0$ and let $\beta>0$ be fixed.
 Then a.s.
$$
  \lim_{N\to\infty} \frac1NG(N,\lfloor\beta N\rfloor) = \frac{\left(1+\sqrt\beta\right)^2}\alpha.
$$
For each $\delta>0$, there exists $c=c(\delta)>0$ such that
\begin{equation}
 \label{eq:concentexp}
 P\left\{\left| \frac1NG(N,\lfloor\beta N\rfloor) - \frac{\left(1+\sqrt\beta\right)^2}\alpha \right| > \delta\right\} \leqslant e^{-cN}
\end{equation}
for all $N\in\N$.
\end{theo}

\begin{theo}
 \label{theo:llngeo}
 Suppose $\xi_y$, 
 $y\in\N^2$,
 are i.i.d and distributed as ${\rm Geom}(\lambda)$, $\lambda\in(0,1)$ and let $\beta>0$ be fixed. Then a.s.
$$
  \lim_{N\to\infty} \frac1NG(N,\lfloor\beta N\rfloor) = \frac{ \left(1+\sqrt{\beta \lambda}\right)^2} {1-\lambda} -1.
$$
For each $\delta>0$, there exists $c=c(\delta)>0$ such that
\[
 P\left\{\left| \frac1NG(N,\lfloor\beta N\rfloor) - \left( \frac{ \left(1+\sqrt{\beta \lambda}\right)^2} {1-\lambda} -1\right) \right| > \delta\right\} \leqslant e^{-cN}
\]
for all $N\in\N$.
\end{theo}

\textbf{Proof of Theorem~\ref{theo:exp}.}
The central idea of the proof may be hidden among all the calculations, basically it is built on the following argument.

By Proposition~\ref{prop:llpp}, Corollary~\ref{lemma:lmonot} and the laws of large numbers for i.i.d. exponential r.v.'s,
$$
   G(N,\beta N) =  H_S(\xi) 
  \leqslant   H_S(\zeta^\circ,\xi) =  \sum_{y\in\partial_-S} \zeta_y^+ + 
  \sum_{y\in\partial^-S} \eta_y^- \approx \frac N {\alpha_+} + \frac {\beta N} {\alpha_-},
$$
where $S$ is a rectangular domain with $N\times \beta N $ sites.
Here $\alpha_+$ and $\alpha_-$ can be any pair of positive numbers that make $\alpha=\alpha_++\alpha_-$ and therefore the broken line process reversible when the $\zeta^\pm$ are distributed as $\exp(\alpha_\pm)$ and the $\xi$ are distributed as $\exp(\alpha)$.
As a consequence we have
\[
 \lim \frac1NG(N,\beta N) 
 \leqslant \inf_{\alpha_+,\alpha_-} \left[ \frac 1 {\alpha_+} + \frac \beta {\alpha_-} \right ]
  = \frac{\left(1+\sqrt\beta\right)^2}\alpha,
\]
the infimum being attained for $\alpha_\pm = \alpha/(1+\beta^{\pm1/2})$.
Now we want to have the opposite inequality. We argue that $H_S(\zeta^\circ,\xi)$ cannot be much bigger than $H_S(\xi)$. To compare both, consider, instead of  boundary conditions $\zeta^\pm$ in $\partial S$, a slightly enlarged domain $S'$ without boundary conditions but where to each extra site we associate $\xi'$ corresponding to the previous $\zeta^\pm$, so that $H_{S'}(\xi')=H_S(\zeta^\circ,\xi)$. Now for $H_{S'}(\xi')$ to be considerably bigger than $H_S(\xi)$, it must be the case that the $\xi'$-optimal path in $S'$ occupies a positive fraction $\epsilon$ of the boundary $\partial_+S'$ and then takes a $\xi$-optimal path 
in the remaining $(N\times(1-\epsilon)\beta N)$-rectangle. But in fact it happens that an oriented walker, after visiting $\epsilon\beta N$ sites in $\partial_+ S$, looks ahead and realizes it is far too late to perform the $\xi'$-optimal path, as the following equation shows:
\begin{equation}
\label{eq:comparetoolate}
 \epsilon\beta \frac{1+\sqrt{1/\beta}}\alpha + \frac{(1+\sqrt{(1-\epsilon)\beta})^2}\alpha 
 = \frac{(1+\sqrt\beta)^2}{\alpha} - \frac{2\sqrt\beta}\alpha f(\epsilon).
\end{equation}
Here
\[
 f(\tau) = 1 - \tau/2 - \sqrt{1-\tau}
\]
is positive and increasing in $(0,\infty)$. Therefore the $\xi'$-optimal path in $S'$ cannot stay too long at the boundary of $S'$ and thus 
$H_S(\xi) \approx H_{S'}(\xi')=H_{S}(\zeta^\circ,\xi)$.

Now let us move to the proper mathematical proof.

We shall use the following basic fact. Given $\alpha_1,\alpha_2,\alpha_3,\rho,\delta,K>0$, there exist positive constants $c,C>0$ such that if $(X_j)$, $(Y_j)$ and $(Z_j)$ are sequences of i.i.d. r.v.'s distributed respectively as $\exp(\alpha_1)$, $\exp(\alpha_2)$ and $\exp(\alpha_3)$, then for all $N\in\N$,
\begin{equation}
 \label{eq:llemma}
 P\left\{
 \begin{array}{r}
 \exists l,m,n\in \{0,\dots,\lceil\rho N\rceil\}:\left|\sum_{j=1}^lX_j+\sum_{j=1}^mY_j + 
 \right.\\ \left. 
 + \sum_{j=1}^nZ_j-\frac l{\alpha_1} - \frac m{\alpha_2} - \frac n{\alpha_3} \right|
 \geqslant \delta N - K
 \end{array}
  \right\} \leqslant Ce^{-cN},
\end{equation}
regardless of the joint distribution of $(X,Y,Z)$.

We first map our problem in $\Z^2$ to the space $\tilde\Z^2$, where the theory of broken lines was developed. We do so by considering the rectangular domain $S(N,M)$ given by $S(N,M)=\{(t,x)\in\tilde\Z^2:0\leqslant t+x \leqslant 2(M-1), 0\leqslant t-x \leqslant 2(N-1)\}$ and the obvious mapping between $S(N,M)$ and $\{1,\dots,N\}\times\{1,\dots,M\}$. We write $S$ for $S(N,M)$.

In order to define a reversible broken line process in $S(N,M)$ with creation $\xi\sim\exp(\alpha)$ we can choose $\alpha_+$ and $\alpha_-$ such that $\alpha_++\alpha_-=\alpha$ and let $\zeta^+\sim\exp(\alpha_+)$, $\zeta^-\sim\exp(\alpha_-)$. Take $M= M(N) = \lfloor\beta N\rfloor$. Choosing $\alpha_+ = \alpha/(1+\beta^{1/2})$ and $\alpha_- = \alpha/(1+\beta^{-1/2})$ gives
\begin{equation}
  \label{eq:alphaopt}
  \frac 1 {\alpha_+} + \frac \beta {\alpha_-} = \frac{\left(1+\sqrt\beta\right)^2}\alpha.
\end{equation}

By Proposition~\ref{prop:llpp} and Corollary~\ref{lemma:lmonot},
\[
 G(N,M) = H_S(\xi) \leqslant H_S(\zeta^\circ,\xi) = \sum_{y\in\partial_-S} \zeta_y^+
  + \sum_{y\in\partial^-S} \eta_y^-,
\]
and it follows from~(\ref{eq:llemma},\ref{eq:alphaopt}) that
\begin{equation}
 \label{eq:concentg}
 P\left\{ \frac1NG(N,M) > \frac{\left(1+\sqrt\beta\right)^2}\alpha + \delta\right\} \leqslant C_0e^{-c_0N}.
\end{equation}

Now let us prove the lower bound to complete the concentration inequality above.
Consider 
$S'(N,M)=\{(t,x)\in\tilde\Z^2:-2\leqslant t+x \leqslant 2(M-1), -2\leqslant t-x \leqslant 2(N-1)\}$ and for $0\leqslant n\leqslant\beta N$, take
$\tilde S(N,M-n)=\{(t,x)\in\tilde\Z^2: 2n\leqslant t+x \leqslant 2(M-1), 0\leqslant t-x \leqslant 2(N-1)\}$,  
$\tilde S'(N,M-n)=\{(t,x)\in\tilde\Z^2: 2n-2\leqslant t+x \leqslant 2(M-1), -2\leqslant t-x \leqslant 2(N-1)\}$, 
$\tilde S(N-n,M)=\{(t,x)\in\tilde\Z^2: 0\leqslant t+x \leqslant 2(M-1), 2n\leqslant t-x \leqslant 2(N-1)\}$, 
$\tilde S'(N-n,M)=\{(t,x)\in\tilde\Z^2: -2\leqslant t+x \leqslant 2(M-1), 2n-2\leqslant t-x \leqslant 2(N-1)\}$.

For $(\zeta^\circ,\xi)$ defined on $S$ take $\xi'$ on $S'$ given by
\[
  \xi'_{t,x} = 
  \begin{cases}
     \xi_{t,x}, & y\in S \\
     \zeta^-_{t+1,x-1}, & t+1,x-1\in \partial_+S \\ 
     \zeta^+_{t+1,x+1}, & t+1,x+1\in \partial_-S \\ 
     0 & \mbox{otherwise}.
  \end{cases}
\]
For $(\tilde\zeta^\circ,\xi)$ defined on $\tilde S$ take $\tilde\xi'$ on $\tilde S'$ given by the analogous formulae.
It is easy to see that $H_S(\zeta^\circ,\xi)=H_{S'}(\xi')$ and $H_{\tilde S}(\tilde\zeta^\circ,\xi)=H_{\tilde S'}(\tilde\xi')$.

The two facts below will be important:
\begin{eqnarray}
 \label{eq:somenopt}
 \mbox{For some } n ,&&
 H_{S'}(\xi') = \sum_{j=0}^{n}\xi'_{j-1,j+1} + H_{\tilde S(N,M-n)}(\xi) \mbox{ or }\\&&\nonumber
 H_{S'}(\xi') = \sum_{j=0}^{n}\xi'_{j-1,-j-1} + H_{\tilde S(N-n,M)}(\xi).  \\
 \label{eq:allnopt}
 \mbox{For all } n ,&& H_{S}(\xi) \geqslant   H_{\tilde S(N,M-n)}(\xi) \mbox{ and } \\&&\nonumber
 H_{S}(\xi) \geqslant   H_{\tilde S(N-n,M)}(\xi).
\end{eqnarray}
Given any $0<\delta'<\delta$, by putting (\ref{eq:somenopt}) and (\ref{eq:allnopt}) together we see that for
\begin{equation}
 \label{eq:smalllpp}
 G(N,M) \leqslant N\left[ \frac{\left(1+\sqrt\beta\right)^2}\alpha - \delta \right]
\end{equation}
to hold, we must have either
\begin{equation}
 \label{eq:ldpreversiblebl}
 H_{S'}(\xi') \leqslant N\left[ \frac{\left(1+\sqrt\beta\right)^2}\alpha - \delta' \right],
\end{equation}
or, for some $n\in\{1,\dots,M-n\}$,
\begin{eqnarray}
 & \left( \sum_{j=0}^n \zeta_{j,j}^- \right ) + H_{\tilde S(N,M-n)}(\tilde\xi) \geqslant 
 N\left[ \frac{\left(1+\sqrt\beta\right)^2}\alpha - \delta' \right]
 \label{eq:largewayaround} \\
 & \sum_{j=0}^n \zeta_{j,j}^- \geqslant N[\delta-\delta']
 \label{eq:largerecovering},
\end{eqnarray}
or, for some $n\in\{1,\dots,N-n\}$,
\begin{eqnarray}
 & \left( \sum_{j=0}^n \zeta_{j,-j}^+ \right ) + H_{\tilde S(N-n,M)}(\tilde\xi) \geqslant 
 N\left[ \frac{\left(1+\sqrt\beta\right)^2}\alpha - \delta' \right]
 \label{eq:largewayaround2} \\
 & \sum_{j=0}^n \zeta_{j,-j}^+ \geqslant N[\delta-\delta']
 \label{eq:largerecovering2}.
\end{eqnarray}

The probability of (\ref{eq:ldpreversiblebl}) decays exponentially fast and this can be shown exactly as was done for (\ref{eq:concentg}).

We consider now the other possibility, (\ref{eq:largewayaround},\ref{eq:largerecovering}). The case (\ref{eq:largewayaround2},\ref{eq:largerecovering2}) is treated in a completely analogous way and is thus omitted. Let $\delta>0$ be fixed, take $\epsilon_0 = {\delta\alpha}/[{2\beta(1+\sqrt{1/\beta})}]$ and $\delta' = \frac\delta3 \wedge [\frac{\sqrt{\beta}}{\alpha} f(\epsilon_0)]$. With this choice of parameters 
\begin{equation}
 \label{eq:epssmallrecov}
 \epsilon\beta\frac{1+\sqrt{1/\beta}}\alpha \leqslant \frac
 \delta2 < \delta-\delta'
\end{equation} 
holds for $\epsilon\leqslant\epsilon_0$ and
\begin{eqnarray}
 \label{eq:longsum}
 \begin{array}{c}
 \displaystyle  \epsilon\beta\frac{1+\sqrt{1/\beta}}\alpha +
 (\beta-\epsilon\beta)\frac{1+(\beta-\epsilon_0\beta)^{-1/2}}\alpha + \frac{1+\sqrt{\beta-\epsilon_0\beta}}\alpha
  \leqslant \\ \qquad \qquad \qquad \qquad \qquad \qquad \qquad \qquad \qquad \qquad
  \displaystyle \leqslant \frac{(1+\sqrt\beta)^2}\alpha - 2 \delta'
 \end{array}
\end{eqnarray}
holds for $\epsilon\geqslant\epsilon_0$.

Consider the event that (\ref{eq:largerecovering}) happens for some $0\leqslant n\leqslant \lceil \epsilon_0\beta N \rceil = M_0$.
Since $\zeta_j^- \stackrel{\rm d}= \exp(\alpha/(1+\sqrt{1/\beta}))$ it follows from~(\ref{eq:epssmallrecov}) and~(\ref{eq:llemma}) that the probability of this event decays exponentially fast in $N$.

It remains to consider $n \geqslant M_0$ and show that in this case it is the probability of~(\ref{eq:largewayaround}) that decays exponentially fast. Now
\begin{eqnarray*}
&&
P\left\{
 \exists n \in\{M_0,\dots,M\}: 
 \sum_{j=0}^n \zeta_{j,j}^- + H_{\tilde S}(\xi) \geqslant 
 N \frac{\left(1+\sqrt\beta\right)^2}\alpha - N\delta'
\right\}
\leqslant \\&&
P\left\{
 \exists n \in\{M_0,\dots,M\}: 
 \sum_{j=0}^n \zeta_{j,j}^- + H_{\tilde S'}(\tilde\xi') \geqslant 
 N \frac{\left(1+\sqrt\beta\right)^2}\alpha - N\delta'
\right\}
= \\&&
P\left\{
 \exists n : 
 \sum_{j=0}^n \zeta_{j,j}^-  + \sum_{j=n}^{M-1} \tilde\zeta_{j,j}^- 
 + \sum_{j=0}^{N-1} \tilde\eta^+_{M-1+j,M-1-j} \geqslant 
 N \frac{\left(1+\sqrt\beta\right)^2}\alpha - N\delta'
\right\},
\end{eqnarray*}
where $\tilde\zeta^\pm$ are distributed as $\exp(\alpha/[1+(\beta-\epsilon_0\beta)^ {\pm1/2}])$ so that the broken line process on $\tilde S$ with $\xi$ distributed as $\exp(\alpha)$ is reversible, and therefore the $\tilde\eta^\pm$ are also distributed as the $\tilde\zeta^\pm$.

By~(\ref{eq:longsum}) the right-hand side of the inequality in the last line is greater than
\[
  n \frac{1+\sqrt{1/\beta}}\alpha + 
  (M-n)\frac{1+(\beta-\epsilon_0\beta)^{-1/2}}\alpha +
  N\frac{1+\sqrt{\beta-\epsilon_0\beta}}\alpha + N \delta'
\]
and by~(\ref{eq:llemma}) the last probability above decays exponentially fast. The proof is finished.
$\hfill \square$\\

\textbf{Proof of Theorem~\ref{theo:llngeo}.}
The proof is absolutely identical to that of the previous theorem, so we just highlight which equations should be replaced by their analogous counterparts.

In the heuristic part take $\lambda_+=[(\lambda+\sqrt{\beta \lambda})/(1+\sqrt{\beta \lambda})]\in(\lambda,1)$ and $\lambda_-=\lambda/\lambda_+\in(\lambda,1)$, so that $\lambda=\lambda_+\lambda_-$, the process is reversible for $\zeta^\pm\sim\exp(\lambda_\pm)$ and 
$$
  \frac{\lambda_+}{1-\lambda_+}+\beta\frac{\lambda_-}{1-\lambda_-} = 
  \frac{ \left(1+\sqrt{\beta \lambda}\right)^2} {1-\lambda} -1.
$$
Instead of (\ref{eq:comparetoolate}) consider
\[
 \epsilon\beta\frac{\lambda_-}{1-\lambda_-} + \frac{(1+\sqrt{(1-\epsilon)\beta\lambda})^2}{1-\lambda}-1
 = \frac{(1+\sqrt{\beta \lambda})^2}{1-\lambda}-1 - \frac{2\sqrt{\beta \lambda}}{1-\lambda}f(\epsilon).
\]

The proof of
\[
 P\left\{ \frac1NG(N,\lfloor\beta N\rfloor) > \frac{ \left(1+\sqrt{\beta \lambda}\right)^2} {1-\lambda} -1 + \delta\right\} \leqslant e^{-cN}
\]
is analogous to the proof of (\ref{eq:concentg}).

For the opposite inequality, we take $\epsilon_0=\delta(1-\lambda)/[2(\lambda\beta+\sqrt{\lambda\beta})]$, $\delta' = \frac\delta3 \wedge [\frac{\sqrt{\beta\lambda}}{1-\lambda} f(\epsilon_0)]$ so that instead of (\ref{eq:epssmallrecov}) and (\ref{eq:longsum}) the following estimates will hold, respectively for $\epsilon\leqslant\epsilon_0$ and $\epsilon\geqslant\epsilon_0$:
\[
 \epsilon\beta \frac{\beta\lambda+\sqrt{\beta\lambda}}{\beta(1-\lambda)} 
 \leqslant \frac \delta2 < \delta-\delta',
\]
\begin{eqnarray*}
 \epsilon\beta \frac{\beta\lambda+\sqrt{\beta\lambda}}{\beta(1-\lambda)} +
 (\beta-\epsilon\beta)\frac{(\beta-\epsilon_0\beta) \lambda+\sqrt{(\beta-\epsilon_0\beta)\lambda}}{(\beta-\epsilon_0\beta)(1-\lambda)}  +
 \frac{\lambda+\sqrt{(\beta-\epsilon_0\beta)\lambda}}{1-\lambda}
  \leqslant \\ \leqslant \frac{(1+\sqrt{\beta\lambda})^2}{1-\lambda} - 1 - 2 \delta'.
\end{eqnarray*}
The rest of the proof is the same.
$\hfill \square$\\

\clearpage
\thispagestyle{empty}

\appendix
\chapter{Technical Proofs}

\section{Existence of the maximal broken line}
\label{sec:maximal}

Here we prove that given a flow field and a broken trace there is a maximal broken line associated to that field and having that trace.

We claim that there exist unique $J_1=(a_1,b_1], J_2=(a_2,b_2], \dots,J_n=(a_n,b_n]$ such that 
\begin{equation}
 \label{eq:gammaell}
 (e_1,J_1)\sim(e_2,J_2)\sim\cdots\sim(e_n,J_n),
\end{equation}
and such that any $p_1,p_2,\dots,p_n$ with property
\begin{equation}
 \label{eq:pgammaell}
 (e_1,p_1)\sim(e_2,p_2)\sim\cdots\sim(e_n,p_n)
\end{equation} must satisfy $p_i\in J_i,\ i=1,\dots,n$.
If $J_1\ne\emptyset$ we define $\gamma(\ell) = (y_0,e_1,J_1,y_1,\dots,e_n,J_n,y_n)$ and by~(\ref{eq:gammaell}) we have $\gamma(\ell) \in B(\bar\eta)$; otherwise $\gamma(\ell)=\emptyset$.

To prove the claim start by observing some consequences of the association rules.
\\ 1. If $(e_1,p_1)\sim(e_2,p_2)$ then $(e_1,p_1-\delta)\sim(e_2,p_2-\delta)$ for some $\delta>0$.
\\ 2. If $(e_1,p_1^k)\sim(e_2,p_2^k)\ \forall\ k\in\N$ with $p_1^k\uparrow p_1$ then $p_2^k\uparrow p_2$ and $(e_1,p_1)\sim(e_2,p_2)$.
\\ 3. If $(e_1,p_1)\sim(e_2,p_2)$ and $(e_1,p_1')\sim(e_2,p_2')$ with $p_1'>p_1$ then $p_1'-p_2'=p_1-p_2$ and $(e_1,(p_1,p_1'])\sim(e_2,(p_2,p_2'])$.

Now let $A$ be the set of $p_1\in(0,\eta(e_1)]$ for which it is possible to find $p_2,p_3,\dots,p_n$ such that~(\ref{eq:pgammaell}) holds.
Suppose that $A\ne\emptyset$ and take $a_1=\inf A$, $b_1=\sup A$.
(When $A=\emptyset$ we take $J_i=\emptyset$ and the desired properties hold trivially.)
Consider a sequence $(p_1^k)\subseteq A$ with $p_1^k\uparrow b_1$.
By Property~2 above we have $(e_1,b_1)\sim\cdots\sim(e_n,b_n)$ and $b_1\in A$.
It follows from Property~1 that $a_1<b_1$ and we can take another sequence $(p_1^k)\subseteq A$ with $b_1>p_1^k\downarrow a_1$, $(e_1,p_1^k)\sim\cdots\sim(e_n,p_n^k)$.
By Property~3 $(e_1,(p_1^k,b_1])\sim\cdots\sim(e_n,(p_n^k,b_n])$ and $b_i>p_i^k\downarrow a_i$ for $i=1,\dots,n$; thus $(e_1,(a_1,b_1])\sim\cdots\sim(e_n,(a_n,b_n])$ and $b_1-a_1=\cdots=b_n-a_n$.
If it were the case that $a_1>0$ and $a_1\in A$, by Property~1 it would hold that $a_1-\delta\in A$ contradicting $a_1=\inf A$.
Therefore $A=(a_1,b_1]$. Suppose $p_1,\dots,p_n$ satisfy (\ref{eq:pgammaell}); by definition $p_1\in A=(a_1,b_1]$ and by Property~3 we have $b_i-p_i=b_1-p_1\in [0,b_1-a_1)=[0,b_i-a_i)$, that is, $p_i\in(a_i,b_i]$.
As a consequence we have that $J_i'\subseteq(a_i,b_i], i=1,\dots,n$ whenever $(e_1,J_1')\sim\cdots\sim(e_n,J_n')$, from which uniqueness follows.

\section{Proof of Lemma~\ref{lemma:lorder}}
\label{sec:lorder}

We start by proving Item~\ref{item:partorder}.
Relation $\succeq$ is obviously reflexive.
We now show that it is antisymmetric.
Assume $\ell\succeq\ell'\succeq\ell$, so that we are required to show that $\ell=\ell'$.
First suppose $D(\ell)\cap D(\ell')\ne\emptyset$ and take $x_0\in D(\ell)\cap D(\ell')$. Write $D(\ell)=\{x_{-n},\dots,x_0,\dots,x_m\}$, $D(\ell')=\{x_{-n'},\dots,x_0,\dots,x_{m'}\}$ and $D(\ell)\cap D(\ell')=\{x_{-\tilde n},\dots,x_0,\dots,x_{\tilde m}\}$ with $x_{i+1}=x_i+1$ and $\tilde m = m\wedge m', \tilde n = n\wedge n'$.
Since $\ell\succeq\ell'\succeq\ell$ we have $t(x_i)=t'(x_i)$ for $i=-\tilde n,\dots,\tilde m$. All we need to show is that $n=n'$ and $m=m'$.
Suppose $\tilde m=m'$. Then $(t'(x_{\tilde m}),x_{\tilde m})\in\partial_+\bar S\cup\partial^+\bar S$; and since $t(x_{\tilde m})=t'(x_{\tilde m})$ we cannot have $m>\tilde m$, thus $m=\tilde m$. By the same argument, if $\tilde m=m$ we conclude $m'=\tilde m$, therefore $m=\tilde m=m'$.
Similarly we show that $n=n'=\tilde n$.
It remains to consider the case $D(\ell)\cap D(\ell')\ne\emptyset$, which is ruled out by the following claim.
\begin{claim}
 \label{claim:empty}
 If $\ell,\ell'\in C(S)$ and $D(\ell)\cap D(\ell')=\emptyset$ then $I(\ell)\cap I(\ell')=\emptyset$
\end{claim}
We will show that $I(\ell)\cap I(\ell')\ne\emptyset$ implies $D(\ell)\cap D(\ell')\ne\emptyset$. Let $\tilde t\in I(\ell)\cap I(\ell')$ and take $x''\in D(\ell), x'\in D(\ell')$ such that $t(x'')=t'(x')=\tilde t$. Assume for simplicity $x''<x'$ and let $(\bar t,\bar x)$ denote the topmost site of $S$. Since $(\tilde t,x')\in\bar S$ we have $\bar x-\bar t\geqslant x'-\tilde t-2$ and $\bar x +\bar t\geqslant x'+\tilde t-2$.
Now, after passing by $(\tilde t,x'')$ when going upwards, $\ell$ must cross either of the lines
$\{(x,t):x-t=\bar x-\bar t\}$ or $\{(x,t):x+t=\bar x+\bar t\}$, because $\ell\in C(S)$. After crossing either of these lines there will be $(t_*,x_*)\in\ell$ with $(x_*-1)-(t_*+1)=\bar x - \bar t$ or $(x_*-1)+(t_*-1)=\bar x + \bar t$, respectively. Therefore $x^*\geqslant x'-|t'-t_*|$. But by~(\ref{eq:lbroktr})  we have $|t'-t_*|\leqslant|x_*-x''|$, so $x'-x_*\leqslant |x_*-x''|$. Assuming $x_*\leqslant x'$ (the other possibility trivially implies the desired result), one has $|x'-x_*|\leqslant|x''-x_*|$ thus $x_*\geqslant \frac{x'+x''}2$ and therefore $\frac{x'+x''}2\in D(\ell)$. Analogously we show that $\frac{x'+x''}2\in D(\ell')$ and the proof is done.

Finally let us see that $\succeq$ is transitive.
For a given point $y_*=(t_*,x_*)\in\bar S$, define
$A(t_*,x_*)=\{(t,x)\in \tilde\Z^2: t-x\geqslant t_*-x_*,t+x\geqslant t_*+x_*\}=\{(t,x)\in \tilde\Z^2: t\geqslant t_*+|x-x_*|\}$ 
and for $\ell\subseteq\bar S$, define $A(\ell)=\cup_{y\in\ell}A(y)$. It is immediate that $A(y)\subseteq A(\tilde y)$ iff $y\in A(\tilde y)$, therefore $A(\ell')\subseteq A(\ell)$ is equivalent to $\ell'\subseteq A(\ell)$. Now let $\ell\succeq\ell'\succeq\ell''$. It follows from these observations and from the claim below that $\ell\succeq\ell''$.
\begin{claim}
 Let $\ell,\ell'\in C(S)$. Then $\ell'\succeq\ell$ if and only if $\ell'\subseteq A(\ell)$.
\end{claim}
We start proving the `only if' part of the claim. Let $\ell'\succeq\ell\in C(S)$.
Suppose $D(\ell)\cap D(\ell')=\emptyset$. By Claim~\ref{claim:empty} $I(\ell)\cap I(\ell')=\emptyset$ and, since $t'(x')>t(x)$ for some $x'\in D(\ell')$, $x\in D(\ell)$, we have $t'(x')>t(x)$ for all $x'\in D(\ell')$, $x\in D(\ell)$. Consider the case $x>x'\ \forall\ x'\in D(\ell')$, $x\in D(\ell)$; the other situation is analogous.
Writing $D(\ell)=\{x_0<\cdots<x_n\}$ and $t_0=t(x_0)$, we must have $(t_0,x_0)\in\partial_-\bar S$ or $(t_0,x_0)\in\partial^-\bar S$ and the later is ruled out since there is $x'\in D(\ell')$ with $x'<x,t'(x')>t_0$ in $\bar S$. Now as $(t_0,x_0)\in\partial_-\bar S$ we have $t+x\geqslant t_0+x_0\ \forall\ (t,x)\in\bar S$ and, as $x'<x_0,t'(x')>t_0\ \forall\ x'\in D(\ell')$ we have $t'-x'>t_0-x_0$ for all $(t',x')\in\ell'$. Therefore $\ell'\subseteq A(t_0,x_0)\subseteq A(\ell)$.
Suppose on the other hand that $D(\ell)\cap D(\ell')\ne\emptyset$,
take $x_0\in D(\ell)\cap D(\ell')$ and write $D(\ell)=\{x_{-n},\dots,x_0,\dots,x_m\}$, $D(\ell')=\{x_{-n'},\dots,x_0,\dots,x_{m'}\}$ and $D(\ell)\cap D(\ell')=\{x_{-\tilde n},\dots,x_0,\dots,x_{\tilde m}\}$ with $x_{i+1}=x_i+1$ and $\tilde m = m\wedge m', \tilde n = n\wedge n'$.
For $x\in D(\ell)\cap D(\ell')$ we have $t'(x)\geqslant t(x)$ and of course $(t'(x),x)\in A(t(x),x)$.
Take $t_m=t(x_m)$, by definition $(t_m,x_m)\in\partial_+\bar S$ or $(t_m,x_m)\in\partial^+\bar S$.
In the later case it must be that $m\geqslant m'$, for if we suppose that $m'\geqslant m$, then as $t'(x_m)\geqslant t(x_m)$, we must have $t'(x_m)=t(x_m)$, thus $(t'(x_m),x_m)\in\partial^+\bar S$ and $m'=m$. In the former case we have $\{(t,x)\in\bar S:x\geqslant x_m\}\subseteq A(t_m,x_m)$. Therefore $(t'(x),x)\in A(\ell)$ for $x\in\{x_0,\dots,x_{m'}\}$. Analogous arguments show that $(t'(x),x)\in A(\ell)$ for $x\in\{x_{-n'},\dots,x_0\}$.

The `if' part is shorter. Suppose $\ell'\subseteq A(\ell)$. Take some $x'\in D(\ell')$, there is $(t,x)\in\ell$ such that $(t'(x'),x')\in A(t,x)$ and thus $t'(x')\geqslant t+|x-x'|\geqslant t$. Now for $x'\in D(\ell')\cap D(\ell)$, there is some $\tilde x\in D(\ell)$ such that $(t'(x'),x')\in A(t(\tilde x),\tilde x)$. We want to show that $t'(x')\geqslant t(x')$ and it follows from the fact that $\ell\subseteq A^-(t,x)$ for any $(t,x)\in\ell$, where
$A^-(t_*,x_*)=\{(t,x)\in \tilde\Z^2: t\geqslant t_*+|x-x_*|\}$.

Item~\ref{item:ellextreme} is easy.
If $\ell=(y_a,y_{a+1},\dots,y_b)\subseteq\ell'=(y_{a'},y_{a'+1},\dots,y_b')$ with $a'\leqslant a<b\leqslant b'$ for $\ell\in C(S)$, $\ell'\subseteq\bar S$, then $a=a'$, $b=b'$ and therefore $\ell=\ell'$; for if $a'<a$ we would have $y_a\in S$, contradicting $\ell\in C(S)$, same for $b>b'$.
Item~\ref{item:useless} is not used in this work and we omit its proof.

\section{Proof of Lemma~\ref{lemma:comparable}}
\label{sec:comparable}
Let $\ell,\ell'\subseteq\bar S$ such that $w(\ell),w(\ell')>0$.
If $t(x)=t'(x)\ \forall x\in D(\ell)\cap D(\ell')$ the result is trivial.
So suppose there is $x_0$ such that $t(x_0)\ne t'(x_0)$ and assume for simplicity that $t(x_0)<t'(x_0)$.
Write $D(\ell)\cap D(\ell') = \{x_{-\tilde m},\dots,x_{-1},x_0,x_1,x_{\tilde n}\}$,
$$\gamma(\ell)=\{y_{-m},e_{-m+1},J_{-m+1},y_{-m+1},\dots,e_0,J_0,y_0,\dots,e_m,J_m,y_m\},$$ 
and $$\gamma(\ell')=\{y'_{-m'},e'_{-m'+1}, J'_{-m'+1},y'_{-m'+1},\dots,e'_0,J'_0,y'_0,\dots, e'_{m'},J'_{m'},y'_{m'}\}.$$

We want to show that $t(x_i)\leqslant t'(x_i)$ for $i=0,1,\dots,\tilde n$; the proof for $i=0,-1,\dots,-\tilde m$ is analogous.

We claim that, for $i=0,1,\dots,\tilde n-1$, the following facts hold:
$t(x_i)\leqslant t'(x_i)$, $t(x_{i+1})\leqslant t'(x_{i+1})$ and, in case $t(x_i) = t'(x_i)$, $t(x_{i+1}) = t'(x_{i+1})$ we have $J_{i+1}\prec J'_{i+1}$, i.e., $p<p'$ for all $p\in J_{i+1}$, $p'\in J'_{i+1}$.
We proceed by induction.
For $i=0$ the result is obvious because of~(\ref{eq:lbroktr}) and $t(x_0)<t'(x_0)$.
Suppose the claim is true for $i=k-1$. We have three possibilities.
Case~A: $t(x_k)=t'(x_k)$ and $t(x_{k-1})=t'(x_{k-1})$.
Case~B: $t(x_k)=t'(x_k)$ and $t(x_{k-1})<t'(x_{k-1})$.
Case~C: $t(x_k)<t'(x_k)$.
In Case~C the claim holds for $i=k$ for the same reason it holds for $i=0$.
In Case~A $J_k\prec J'_k$ and, because of the association rules, we have $t(x_{k+1})<t'(x_{k+1})$ or $t(x_{k+1})=t'(x_{k+1})$ with $J_{k+1}\prec J'_{k+1}$; either way case the claim holds for $i=k$.
In case~B, since $t(x_{k-1})<t'(x_{k-1})$, by the association rules we cannot have $t(x_{k+1})>t'(x_{k+1})$ and if $t(x_{k+1})=t'(x_{k+1})$ we must have $J_{k+1}\prec J'_{k+1}$, so the claim holds.

The proof is complete.

\section{Formalization of the brick diagram}
\label{sec:etaell}
In this appendix we formalize the construction of the brick diagram described in the proof of
Theorem~\ref{theo:etaell}. Then we show that~(\ref{eq:lweightcalc}) holds for any flow field in the converse part of the theorem and that~(\ref{eq:lbtcarac}) holds for some flow field in its direct part.

As mentioned in the proof of the theorem, there are several equivalent representations of a flow field and this proof draws on each of them.
We describe how to obtain from each representation the next one, and we mention the properties of each representation that guarantee it is possible to come back to the previous setting.
Statements will be made without proof when their verification is a straight forward but tedious standard argument.

To simplify the presentation, assume that $S$ has the form $S=\{(t,x)\in\tilde\Z^2:|t|+|x|\leqslant N\}$ for some $N\in2 \Z_+$.

Let $S^* = \{y\in\Z^2:\exists y'\in S:|y-y'|=1\}$. Notice that $S^*\subseteq(\tilde\Z^2)^*=\Z^2\backslash\tilde\Z^2$.

The first alternative representation of a flow field is $p_{t,x}, (t,x)\in S^*$, defined below.

Take $W_{-N-1}=0$, for $i=-N,\dots,0$ take $W_i = W_{i-1} + \eta(e_{i,N+i}^\nwarrow)$ and for $i=1,\dots,N+1$ take $W_i = W_{i-1} + \eta(e_{i-1,N-i+1}^\nearrow)$.

Let $p_{0,N+1} = W_0$.

For each $x=N,N-1,\dots,1,0,-1,\dots,-N,-N-1$, set
\[
 p_{t,x} = p_{t+1,x+1} - \eta(e_{t+1,x}^\nwarrow) 
\]
for $t=-(N+1-x),-(N+1-x)+2,\dots,(N+1-x)-2$ and
\[
 p_{t,x} = p_{t-1,x+1} + \eta(e_{t-1,x}^\nearrow)
\]
for $t=N+1-x$.

With successive uses of (\ref{eq:lcons2}) it is not hard to see that $p_{t,x}$ has the following properties:
\begin{eqnarray*}
  && p_{-N-1,0}=0,
  \\ &&
  p_{t,x} \leqslant p_{t+1,x\pm1}\qquad \mbox{whenever they belong to }S^*
\end{eqnarray*}
and that it is possible to re-obtain $\bar\eta$ from the $p(t,x)$ by
\begin{equation}
\label{eq:etap}
\begin{array}{l}
  \eta(e_y^\swarrow) = p_{t,x-1} - p_{t-1,x}
\\
  \eta(e_y^\searrow) = p_{t+1,x} - p_{t,x-1}
\\
  \eta(e_y^\nwarrow) = p_{t,x+1} - p_{t-1,x}
\\
  \eta(e_y^\nearrow) = p_{t+1,x} - p_{t,x+1}.
\end{array}
\end{equation}

Now consider the set of all $p_{t,x}$ and reorder it by taking $A=\{q_0,\dots,q_M\}=\{p_{t,x}: (t,x)\in S^*\}$ with $q_0<q_1<\cdots<q_M$. Then $q_0=W_{-N-1}=p_{-N-1,0}=0$ and $q_M=W_{N+1}=p_{N+1,0}$.

For $(t,x)\in S^*$, take $k(t,x)\in\{0,\dots,M\}$ as the unique subindex $k$ that satisfies $p_{t,x}=q_{k}$. Then
\begin{eqnarray}
  && k(-N-1,0)=0 
  \nonumber \\ &&
  k(N+1,0)=M
  \nonumber \\ &&
  k(t,x) \leqslant k(t+1,x\pm1)\qquad \mbox{whenever they belong to }S^*
  \label{eq:kconsist} \\ &&
  \mbox{For all } k=0,\dots,M \mbox{ there is } y\in S^* \mbox{ such that } k(y)=k.
  \nonumber
\end{eqnarray}
Of course it is possible to recover $p_{t,x}$ from $A$ and $k(t,x)$:
\[
  p_{t,x} = q_{k(t,x)}.
\]

For $y=(t,x)\in \bar S$, take $k_\pm(y)\in\{0,\dots,M\}$ as
\[
 k_-(t,x) =
 \begin{cases}
   k(t-1,x),& y\in S\cup\partial^+\bar S\cup\partial^-\bar S \\ 
   k(t,x-1),& y\in\partial_+(S) \\
   k(t,x+1),& y\in\partial_-(S).
 \end{cases}
\]
\[
 k_-(t,x) = 
 \begin{cases}
   k(t+1,x),& y\in S\cup\partial_+\bar S\cup\partial_-\bar S \\ 
   k(t,x-1),& y\in\partial^+(S) \\
   k(t,x+1),& y\in\partial^-(S).
 \end{cases}
\]
Then
\begin{eqnarray*}
  && k_-(-N,0)=0
  \\ &&
  k_+(N,0) = M
  \\ &&
  k_+(t,x) \leqslant k_+(t+1,x\pm1)\qquad \mbox{whenever they belong to } \bar S
  \\ &&
  k_-(t,x) \leqslant k_-(t+1,x\pm1)\qquad \mbox{whenever they belong to } \bar S
  \\ &&
  k_+(t,x) = k_-(t+2,x)\qquad \mbox{whenever they belong to } \bar S
  \\ &&
  \mbox{For all } k=1,\dots,M \mbox{ there is } y\in S \mbox{ such that } k_+(y)=k
  \\ &&
  \mbox{For all } k=0,\dots,M-1 \mbox{ there is } y\in S \mbox{ such that } k_-(y)=k.
\end{eqnarray*}
To obtain $k$ from $k_\pm$, take, for $(t,x)\in S^*$,
\[
 k(t,x) = 
 \begin{cases}
   0 ,& t=-N-1 \\ 
   M ,& t=N+1 \\ 
   k_-(t+1,x)=k_+(t-1,x), & \mbox{otherwise}
 \end{cases}.
\]

We may also consider $t(j,x)$, $j \in \{1,\dots,M\}, x\in\{-N,\dots,N\}$.
For such $(j,x)$, take
\[
 t(j,x)=
 \begin{cases} 
 t, & (t,x)\in\bar S \mbox{ and } k_-(t,x)<j\leqslant k_+(t,x) \\ 
 -\infty,&  j\leqslant k_-(t,x) \mbox{ for all $t$ such that } (t,x)\in\bar S \\ 
 +\infty,&  j > k_+(t,x) \mbox{ for all $t$ such that } (t,x)\in\bar S.
 \end{cases}
\]
When $|t(j,x)|\ne\infty$, $(t(j,x),x)$ is an element of $\bar S$.
Given a fixed $j$, if $t(j,x_0)=+\infty$ for some $x_0>0$ (resp. $<0$) then $t(j,x)=+\infty$ for all $x\geqslant x_0$ (resp. $\leqslant x_0$); same for $-\infty$.
For $x=0$ it is always the case that $t(j,0)\in\{-N,\dots,N\}$.
Also, $t(j,x)\leqslant t(j+1,x)$ for $j\in\{1,\dots,M-1\}, x\in\{-N,\dots,N\}$.
For all $j\in\{1,\dots,M-1\}$ there is $x\in\{-N,\dots,N\}$ such that $t(j,x)< t(j+1,x)$.
Also $t(j,x+1) = t(j,x)\pm 1$ for all $j\in\{1,\dots,M\}, x\in\{-N,\dots,N-1\}$ such that $|t(j,x+1)|,|t(j,x)|\ne\infty$.
If $|t(j,x)|<\infty$ and $|t(j,x+1)|=\infty$ or $|t(j,x-1)|=\infty$ then $(t(j,x),x)\in\partial\bar S$.

To relate $t(j,x)$ with $k_\pm(t,x)$ take $t(0,x)=-\infty$, $t(M+1,x)=+\infty$ and
\[
  k_-(t,x) = \max\{k:0\leqslant k\leqslant M, t(k,x)<t\},
\]
\[
  k_+(t,x) = \min\{k:0\leqslant k\leqslant M, t(k+1,x)>t\}.
\]

Now for $j\in\{1,\dots,M\}$ take $a = \min\{x:t(j,x)\ne\pm\infty\}, b = \max\{x:t(j,x)\ne\pm\infty\}$ and define
$ y_i = (t(j,x),x) $
for $a\leqslant x\leqslant b$.
Take $\ell_j = (y_a,e_{a+1},y_{a+1},\dots,e_b,y_b)$, $e_j=\langle y_{j-1},y_j\rangle$.
Then $\ell_1,\dots,\ell_M\in C(S)$ and
\begin{equation}
 \label{eq:ellconstrord}
 \ell_1\prec\ell_2\prec\cdots\prec\ell_M.
\end{equation}

Given the set $\ell_1\prec\cdots\prec\ell_M$, write $D(\ell_j)=\{x^-_j,x^-_j+1,\dots,-1,0,1,\dots,x^+_j-1,x^+_j\}$
and $y=(t_j(x),x)$ for $y\in\ell_j$.
Then $(t_j(x^\pm_j),x^\pm_j)\in\partial_\pm\bar S\cup\partial^\pm\bar S$. Take
\[
 t(j,x)=
 \begin{cases} 
    t_j(x), & x_j^-\leqslant x\leqslant x_j^+ \\
    +\infty, & x>x_j^+ \mbox{ and } t_j(x_j^+)>0 \mbox{ or }
               x<x_j^- \mbox{ and } t_j(x_j^-)>0 \\
    -\infty, & x>x_j^+ \mbox{ and } t_j(x_j^+)<0 \mbox{ or }
               x<x_j^- \mbox{ and } t_j(x_j^-)<0.
  \end{cases}
\]

This completes the set of equivalences:
\[
 \bar\eta \leftrightarrow p \leftrightarrow (A,k), \quad
  k \leftrightarrow k_\pm \leftrightarrow t \leftrightarrow \ell.
\]

For $y=(t,x)\in \bar S$, take $K(y)=\{k_-(y)+1,\dots,k_+(y)\}$. It follows from the definition of $t(j,x)$ that
\begin{equation}
 \label{eq:tjequiv}
 \mbox{for all }(t,x)\in\bar S,j\in\{1,\dots,M\},
 \qquad j\in K(t,x) \Leftrightarrow t(j,x)=t.
\end{equation}

For $e\in\E(\bar S)$, take $k_-(e)=k_-(y) \vee k_-(y')$ and $k_+(e) = k_+(y) \wedge k_+(y')$, where $e=\langle y,y' \rangle$. Also take $K(e)=\{k_-(e)+1,\dots,k_+(e)\}=K(y)\cap K(y')$.
Yet for $\ell=\left( y_0,e_1,y_1,e_2,y_2,\dots,e_n,y_n \right)\subseteq\bar S$, take $K(\ell)=K(y_0)\cap\cdots\cap K(y_n)=K(e_1)\cap\cdots\cap K(e_n)$.

\begin{claim}
 \label{claim:1}
 Let $\ell\subseteq\bar S$ and $j\in\{1,\dots,M\}$. Then
 $\ell\subseteq\ell_j$ if and only if $j\in K(\ell)$.
\end{claim}

Proof: Take $\ell\subseteq\ell_j$. Write $\ell=(y_{a'},e_{a'+1},y_{a'+1},\dots,y_{b'})$ and $\ell_j=(y_{a},e_{a+1},y_{a+1},\dots,y_{b})$ with $a\leqslant a'<b'\leqslant b$. By construction of $\ell_j$, for all
$a\leqslant i\leqslant b$, $t(j,x_i)=t_i$ and, by~(\ref{eq:tjequiv}), $j\in K(y_i)$. Therefore $j\in K(a')\cap\cdots
\cap K(b') = K(\ell)$.
Now suppose $j\in K(\ell)$ for some $\ell\subseteq\bar S$. Write $\ell=(y_0,e_1,y_1,\dots,y_n)$ and, for $i=0,\dots,n$, $y_i=(t_i,x_i)$. By definition $j\in K(y_i)$ for all $i$, by~(\ref{eq:tjequiv}) this implies $t(j,x_i)=t_i$ and by construction of $\ell_j$ we have $y_i\in\ell_j\ \forall\ i$; it follows from this last fact that $\ell\subseteq\ell_j$, thus proving the claim.

Let $J^j=(q_{j-1},q_j],\ j\in\{1,\dots,M\}$. Notice that
\begin{equation}
 \label{eq:intrvldisjoint}
 J^j \cap J^{j'} = \emptyset\mbox{ when } j\ne j'.
\end{equation}
For $y\in S$, take $ J(y) = \bigcup_{j\in K(y)}J^j $ and for $\ell\subseteq\bar S$ take $J(\ell) = \bigcup_{j\in K(\ell)}J^j$. Given $e=\langle y,y'\rangle\in\E(\bar S)$, let $J(e)=(q_{k_-(e)},q_{k_+(e)}]$. Notice that $J(e)=J(y)\cap J(y')$, since $J(e)=\cup_{j\in K(e)}J^j$.
It follows from~(\ref{eq:etap})-(\ref{eq:kconsist}) that
\[
 |J(e)| = \eta(e).
\]

By the above fact there is one translation $T_e$ that relates $(0,\eta(e)]$ with $J(e)$.
We associate the atoms and subintervals of these translated intervals according to the following rules. We write $(e_1,p_1)\approx(e_2,p_2)$ when $(e_1,T_{e_1}^{-1}p_1)\sim(e_2,T_{e_2}^{-1}p_2)$ and $(e_1,J_1)\approx(e_2,J_2)$ when $(e_1,T_{e_1}^{-1}J_1)\sim(e_2,T_{e_2}^{-1}J_2)$.

\begin{claim}
\label{claim:ltransequiv}
 Let $e_1=\langle y_0,y_1\rangle,e_2=\langle y_1,y_2\rangle$ with $x_0<x_1<x_2$ be two adjacent edges in $\E(\bar S)$ and $J_1,J_2$ intervals.
 Then $(e_1,J_1)\approx(e_2,J_2)$ if and only if $J_1=J_2$ and $J_1\subseteq J(y_0)\cap J(y_1)\cap J(y_2)$.
\end{claim}

To prove the claim it suffices to show that $(e_1,p_1)\approx(e_2,p_2)$ iff $p_1=p_2\in J(e_1)\cap J(e_2)$, since $T_e:(0,\eta(e)]\to J(e)$ is just a translation and $J(e_1)\cap J(e_2)=J(y_0)\cap J(y_1)\cap J(y_2)$. 

So suppose $(e_1,p_1)\approx(e_2,p_2)$. Write $p_1=T_{e_1}\tilde p_1$ and $p_2=T_{e_2}\tilde p_2$. Of course $\tilde p_1\in(0,\eta(e_1)]$, thus $p_1\in J(e_1)$; also $p_2\in J(e_2)$. Writing $T_e$ more explicitly one gets $T_ep = p+q_{k_-(e)}$. We have 4 cases to consider. Case~1: $e_1=e_{y_1}^\swarrow$, $e_2=e_{y_1}^\nwarrow$.
In this case $\tilde p_1=\tilde p_2$, $k_-(y_0) \leqslant k_-(y_1)$ and $k_-(y_2) \leqslant k_-(y_1)$, thus $k_-(e_1)=k_-(e_2)$ and $p_1=\tilde p_1 + q_{k_-(e_1)} = \tilde p_2 + q_{k_-(e_2)} = p_2$, since by hypothesis $\tilde p_1=\tilde p_2$.
Case~2: $e_1=e_{y_1}^\searrow$, $e_2=e_{y_1}^\nwarrow$. In this case $\tilde p_1=\tilde p_2 - \eta(e_{y_1}^\swarrow)$ and $k_-(y_2)\leqslant k_-(y_1)$, thus $k_2(e_2) = k_-(y_1)$. By~(\ref{eq:etap}) $q_{k_-(y_1)} + \eta(e_{y_1}^\swarrow) = q_{k(t_1-1,x_1)}+\eta(e_{y_1}^\swarrow) = q_{k(t_1,x_1-1)} = q_{k_-(e_1)}$, so $p_1=\tilde p_1 + q_{k_-(e_1)} = \tilde p_1 + q_{k_-(y_1)} + \eta(e_{y_1}^\swarrow) = \tilde p_2 + q_{k_-(y_1)} = p_2$. Cases 3 and 4 are similar.

Conversely, suppose $p_1=p_2\in J(e_1)\cap J(e_2)$. Consider Case~3, i.e., $e_1=e_{y_1}^\swarrow$, $e_2=e_{y_1}^\nearrow$; the other cases are similar.
Take $\tilde p_1= T_{e_1}^{-1}p_1 = p_1-q_{k_-(e_1)}= p_1-q_{k_-(y_1)} \in (0,\eta(e_1)]$
and $\tilde p_2= T_{e_2}^{-1}p_2 = p_2-q_{k_-(e_2)}= p_2-q_{k_-(y_2)} = p_2-q_{k_-(t_1+1,x_1+1)} = p_2 - p_{x_1+1,t_1} = p_2 - [p_{x_1,t_1-1}+\eta(e_{y_1}^\nwarrow)] = p_2 - q_{k_-(y_1)} - \eta(e_{y_1}^\nwarrow) = \tilde p_1-\eta(e_{y_1}^\nwarrow) \in (0,\eta(e_2)]$. Since $\tilde p_2>0$ we have $\tilde p_1>\eta(e_{y_1}^\nwarrow)$. Therefore $(e_1,p_1)\approx(e_2,p_2)$ and the claim holds.

\index{broken line!translated}
It will be convenient to work with a different representation of a broken line, which we dub \emph{translated broken line}. The translated broken lines have a simpler representation that follows from Claim~\ref{claim:ltransequiv}.
Given a broken line $\gamma\in B(\bar\eta)$ in $\bar S$, we define the object $\psi=\psi(\gamma)$ by
\[
 \psi = \left((y_0,e_1,y_1,\dots,e_n,y_n),J\right),
\]
where $(y_0,e_1,y_1,\dots,e_n,y_n)=\ell(\gamma)$ and $J=T_{e_1}J_1=T_{e_2}J_2=\cdots=T_{e_n}J_n$.
Write $\gamma(\psi)$ for the unique broken line $\gamma$ such that $\psi=\psi(\gamma)$. For any $\psi$ of the above form we have $\gamma(\psi) = (y_0,e_1,J_1,y_1,\dots,e_n,J_n,y_n)$, where $J_i=T_{e_i}^{-1}J$.
The translated broken lines have all the properties analogous to those already discussed for the broken lines.
For $\psi=(\ell_1,J)$, define $\ell(\psi)=\ell_1$ and $w(\psi)=|J|=w(\gamma(\psi))$.

\begin{claim}
 \label{claim:2}
 Let $\ell=(y_0,e_1,y_1,\dots,e_n,y_n)\subseteq\bar S$ and $J\subseteq(0,q_M]$.
 Then $(e_1,J)\approx\cdots\approx(e_n,J)$ if and only if $J\subseteq J(\ell)$.
\end{claim}
The proof is short. Suppose $(e_1,J)\approx\cdots\approx(e_n,J)$. Since $(e_1,J)\approx(e_2,J)$, by Claim~\ref{claim:ltransequiv} we have $J\subseteq J(y_0)\cap J(y_1)\cap J(y_2)$. Also $(e_2,J)\approx(e_3,J)$, thus $J\subseteq J(y_3)$ as well. Similarly, for $i=2,\dots,n$, $(e_{i-1},J)\approx(e_i,J)$ and by Claim~\ref{claim:ltransequiv} $J\subseteq J(y_i)$. Therefore $J\subseteq\cap_{i=0}^n J(y_i) = \cap_{i=0}^n \cup_{j\in K(y_i)} J^j = \cup_{j\in [\cap_{i=0}^n K(y_i)]} J^j = 
\cup_{j\in K(\ell)} J^j = J(\ell)$; we have used~(\ref{eq:intrvldisjoint}) on the second equality. Conversely, if $J\subseteq J(\ell)$ we have $J\subseteq J(y_i)$ for $i=0,\dots,n$ and by Claim~\ref{claim:ltransequiv} we have $(e_{i-1},J)\approx(e_i,J)$ for all $i=2,\dots,n$, i.e., $(e_1,J)\approx\cdots\approx(e_n,J)$, which proves the claim.

We define $\tilde B(\bar\eta) = \{\psi(\gamma):\gamma\in B(\bar\eta)\}$. By definition of translated broken lines we have $\psi\in\tilde B(\bar\eta)$ if and only if $(e_1,J)\approx\cdots\approx(e_n,J)$. It follows from Claim~\ref{claim:2} that 
\begin{equation}
 \label{eq:tblequiv}
 \tilde B(\bar\eta) = \{\psi=(\ell,J):\ell\subseteq\bar S, J\subseteq J(\ell)\}.
\end{equation}
As for the broken lines, given an $\ell\subseteq\bar S$, one can define the maximal translated broken line $\psi(\ell)$ in $\tilde B(\bar\eta)$ that has trace $\ell$. It is clear that $\psi(\ell)=\psi(\gamma(\ell))$ and $\gamma(\ell)=\gamma(\psi(\ell))$. Now by~(\ref{eq:tblequiv}) one has $\psi(\ell) = (\ell,J(\ell))$ and therefore
\begin{equation}
 \label{eq:constrwcalc}
 w(\ell) = \sum_j w_j \I_{j\in K(\ell)}.
\end{equation}

Let $\ell\in C(S)$. By Claim~\ref{claim:1} and Item~\ref{item:ellextreme} of Lemma~\ref{lemma:lorder} we have $j\in K(\ell)$ if and only if $\ell=\ell_j$. Thus we have by (\ref{eq:ellconstrord})
\begin{equation}
 \label{eq:constrkell}
 K(\ell) =
 \begin{cases} 
    \{j\}, & \ell=\ell_j \\ 
    \emptyset, & \mbox{otherwise}
 \end{cases}
 \qquad \mbox{for any }\ell\in C(S).
\end{equation}

Now (\ref{eq:lbtcarac}) follows from (\ref{eq:constrkell}) and (\ref{eq:constrwcalc}).
Finally, combining (\ref{eq:constrwcalc}) and Claim~\ref{claim:1} gives
\[
 w(\ell) = \sum_j w_j \I_{\ell\subseteq\ell_j}
\]
and (\ref{eq:lweightcalc}) follows from (\ref{eq:lbtcarac}) and the above equation.

It remains to prove that, given any pair of sets $\{\ell_1\prec\ell_2\prec\cdots\prec\ell_M\}$ and $\{w_1,\dots,w_M>0\}$, equation~(\ref{eq:lbtcarac}) holds for some $\bar\eta$.
For $j=0,1,\dots,M$, let $q_j=\sum_{j'=1}^jw_{j'}$ and let $A=\{0=q_0<q_1<\cdots<q_M\}$.
For $x\in\{-N,\dots,N\}$, $j\in\{1,\dots,M\}$ define $t(j,x)$ from the $\{\ell_j\}$ as shown above.
It then follows that $t(j,x)$ has all the properties mentioned in the construction.
From $t(j,x)$, define $k_\pm(y),y\in\bar S$ and then $K(y),y\in S^*$.
With $A$ defined above and $k(\cdot)$ we can recover $p_{t,x},(t,x)\in S^*$ and from that we obtain the associated flow field $\bar\eta$.
Now for such $\bar\eta$~(\ref{eq:lbtcarac}) holds by the construction described in this appendix.

\section{Proof of~(\ref{eq:intersect})}
\label{sec:intersect}

As in the proof of Lemma~\ref{lemma:lorder}, take $A(t,x)=\{(\tilde t,\tilde x):\tilde t\geqslant t+|\tilde x-x|\}$, $A(V)=\cup_{y\in V}A(y)$. Consider also $\tilde A(t,x)=\{(\tilde t,\tilde x):\tilde t\leqslant t-|\tilde x-x|\}$. Then $y\in A(y')$ iff $y'\in\tilde A(y)$. It is clear that $y_m\in A(\ell)$, so take $n=\min\{i:y_i\in A(\ell)\}$.
There is $y_*\in\ell$ such that $y_n\in A(y_*)$ and thus $y_*\in\tilde A(y_n)$.
If $n=0$, $y_*\in\bar S\cap\tilde A(y_0)=\{(x_0,t_0),(x_0+1,t_0-1),(x_0-1,t_0-1)\}$.
It is clear that any of these possibilities for $y_*$ imply $y_0\in\ell$.
So suppose $n>0$.
By construction $x_{n-1}=x_n\pm1$, so assume for simplicity $x_{n-1}=x_n+1$.
Now $y_{n-1}\not\in A(\ell)$, which means $\ell\cap\tilde A(y_{n-1})=\emptyset$. So $y_*\in\tilde A(y_n)\backslash\tilde A(y_{n-1})$ and thus $x_*\leqslant x_n$, $t_*=t_n-(x_n-x_*)$.
Notice that $\ell$ must eventually reach $\{(t,x):t-(t_n+1)=(x_n+1)-x\}$ to cross $\bar S$. Let $y'=(t',x')$ be the point of the first time it happens, that is, the one with smallest $x'$. Since $\ell\cap\tilde A(y_{n-1})=\emptyset$, we must have $t'\geqslant t_n+1$, so $x'\leqslant x_n+1$.
Thus $t'-t_* = (t'-t_n)+(t_n-t_*) \geqslant 1+x_n-x_* \geqslant x'-x*$ and the equality holds if and only if $t'-t_n=1,x'-x_n=1$. But it must be the case that equality holds because (\ref{eq:lbroktr}) implies $|t'-t_*|\leqslant|x'-x_*|$. So $y'=(t_n+1,x_n+1)\in\ell$.
Now we only need to observe that $\ell$ must connect $y_*$ to $y'$ through $y_n$. Indeed, if $t(x_n)>t_n$ we would have $t(x_n)-t(x_*)>x_n-x_*$ and if $t(x_n)<t_n$ we would have $t(x')-t(x_n)>x'-x_n$ and either of them is absurd because of~(\ref{eq:lbroktr}).

\backmatter

\parskip 0pt
\addcontentsline{toc}{chapter}{Index}
\begin{theindex}

  \item abelian property, \see{commutativity}{11}
  \item annihilation, 30
    \subitem rule, 30
  \item association rules
    \subitem continuous case, 42
    \subitem discrete case, 38
  \item atoms, 42

  \indexspace

  \item $B(\bar\eta)$, \see{broken line associated to $\bar\eta$}{43}
  \item birth process, 29, 41
  \item boundary condition, 41
  \item brick diagram, 46
  \item broken line, 30, 36, 42
    \subitem associated to $\bar\eta$, 43
    \subitem domain of, 42
    \subitem maximal, 43
    \subitem trace of, 42
    \subitem translated, 66
    \subitem weight of, 42
  \item broken line process
    \subitem continuum, 30
    \subitem geometric, 29
  \item broken trace, 42
    \subitem domain of, 43

  \indexspace

  \item commutativity, 9, 11, 14
  \item consistency property, 35
  \item creation, \see{birth process}{29}
  \item $C(S)$, 43

  \indexspace

  \item domain, \see{hexagonal domain}{31}, 
		\see{rectangular domain}{41}, 
		\see{broken line domain}{42}, 
		\see{broken trace domain}{43}

  \indexspace

  \item envelopes, 9
  \item $\E(\bar S)$, 37
  \item $\bar\eta$, 33, 37
  \item $\eta$, 6
  \item $\bar\eta(\gamma)$, 43
  \item $\bar\eta(\ell)$, 43
  \item $\bar\eta(\zeta^\circ,\xi)$, 41

  \indexspace

  \item flow field, 37, 41
    \subitem decomposition, 44

  \indexspace

  \item $\gamma(\ell)$, \see{broken line, maximal}{43}
  \item generation, 38
  \item generator, 8

  \indexspace

  \item hexagonal domain, 31

  \indexspace

  \item intersection process, 30, 39

  \indexspace

  \item last passage percolation, 49
  \item $L(\gamma)$, 43
  \item $L(\ell)$, 43
  \item $\ell \succeq \ell'$, 44
  \item $\ell\subseteq\ell'$, 44
  \item $\ell\subseteq\bar S$, 43

  \indexspace

  \item monotonicity, 10, 14, 20

  \indexspace

  \item $\NN_m$, 12

  \indexspace

  \item $\tilde\Omega$, 12
  \item $\tilde\omega$, 12

  \indexspace

  \item percolation, \see{last passage percolation}{49}
  \item phase transition, 7

  \indexspace

  \item rectangular domain, 41
  \item relative age, 37
  \item restricted process, 35
  \item $R_x$, 14
  \item $R_x^t$, 20

  \indexspace

  \item $\Sigma_S$, 33
  \item stabilize, 14
  \item sweep, 22

  \indexspace

  \item trap, 22

  \indexspace

  \item $w(\ell)$, 44

  \indexspace

  \item $\xi(\ell)$, 43
  \item $\xi(\gamma)$, 43

\end{theindex}

\cleardoublepage
\addcontentsline{toc}{chapter}{List of Figures}
\listoffigures

\cleardoublepage

\end{document}